\documentclass{article}

\usepackage{amsmath}
\usepackage{amsfonts}
\usepackage{amssymb}
\usepackage{amsthm}
\usepackage{comment}
\usepackage{epsfig}
\usepackage{psfrag}
\usepackage{mathrsfs}
\usepackage{amscd}
\usepackage[all]{xy}
\usepackage{rotating}
\usepackage{lscape}
\usepackage{amsbsy}
\usepackage{verbatim}
\usepackage{moreverb}
\usepackage{fullpage}

\newtheorem{theorem}{Theorem}[section]
\newtheorem{prop}[theorem]{Proposition}
\newtheorem{lemma}[theorem]{Lemma}
\newtheorem{cor}[theorem]{Corollary}

\theoremstyle{definition}
\newtheorem{definition}[theorem]{Definition}

\newtheorem{notation}[theorem]{Notation}
\newtheorem{problem}[theorem]{Problem}

\theoremstyle{remark}
\newtheorem{remark}[theorem]{Remark}

\newtheorem{example}[theorem]{Example}

\DeclareMathOperator{\Aut}{Aut}

\DeclareMathOperator{\Inn}{Inn}
\DeclareMathOperator{\Out}{Out}

\DeclareMathOperator*{\hocolim}{hocolim}

\DeclareMathOperator{\Hom}{Hom}
\DeclareMathOperator{\Iso}{Iso}
\DeclareMathOperator{\Inj}{Inj}

\DeclareMathOperator{\Syl}{Syl}
\DeclareMathOperator{\res}{res}

\newcommand{\id}{\mathrm{id}}
\newcommand{\Ob}{\mathrm{Ob}}

\newcommand{\Mor}{\mathrm{Mor}}
\newcommand{\iso}{\mathrm{iso}}

\def\cA{\mathcal A}\def\cB{\mathcal B}\def\cC{\mathcal C}
\def\cE{\mathcal E}\def\cF{\mathcal F}\def\cG{\mathcal G}\def\cH{\mathcal H}
\def\cK{\mathcal K}\def\cL{\mathcal L}
\def\cO{\mathcal O}
\def\cS{\mathcal S}\def\cT{\mathcal T}
\def\cX{\mathcal X}
\def\cZ{\mathcal Z}

\def\cB{\mathcal B}\def\cC{\mathcal C}
\def\cE{\mathcal E}\def\cF{\mathcal F}\def\cH{\mathcal H}
\def\cK{\mathcal K}\def\cL{\mathcal L}
\def\cO{\mathcal O}
\def\cS{\mathcal S}\def\cT{\mathcal T}
\def\cX{\mathcal X}
\def\cZ{\mathcal Z}

\def\fN{\mathfrak N}\def\fX{\mathfrak X}

\def\fg{\mathfrak g}\def\fh{\mathfrak h}\def\fri{\mathfrak i}\def\fp{\mathfrak p}
\def\fU{\mathfrak U}

\def\FF{\mathbb F}

\def\XX{\mathbb X}

\def\TOP{\mathcal{TOP}}\def\GRP{\mathcal{GRP}}\def\GRPD{\mathcal{GRPD}} \def\CAT{\mathcal{CAT}} 

\def\id{\mathrm{id}}\def\Id{\mathrm{Id}}

\title{Fusion action systems}
\author{Matthew Gelvin}

\numberwithin{theorem}{subsection}

\begin{document}

\maketitle

\addtocounter{section}{-1}

\section{Introduction}

Let $p$ be a prime.

The theory of fusion systems is an abstraction of the $p$-local structure of finite groups.  First codified by Puig in language that can be found in \cite{PuigFrobeniusCategories}, fusion systems have since become a subject of interest to both group theorists and block theorists.  There are also applications to algebraic topology, originally in search of a solution to the Martino-Priddy Conjecture, which leads to the theory of $p$-local finite groups or classifying spaces of fusion systems, as introduced in \cite{BLO2}.  Much of the recent work on these subjects uses a different, though equivalent, formalism to describe fusion systems; a good source for this point of view is \cite{LinckelmannIntro}.

The work of this paper, which was mostly done toward completion of the author's doctoral thesis, began with an investigation of what it would mean for a fusion system to act on a finite set.  An initial answer took the form of a fusion system $\cF$ on a $p$-group $S$, together with a finite set $X$ with an $S$-action that ``respected fusion data.''  These \emph{$\cF$-stable $S$-sets} appear in \cite{BLO2}, where they are used to compute the cohomology ring of $\cF$.  From the topological point of view, a finite $G$-set gives rise to a covering space of $BG$; the hope that a similar situation would realize for $p$-local finite groups.  

Ultimately, it became apparent that simply imposing a condition on an $\cF$-set did not yield enough structure to do, for example, unstable homotopy theory.  To rigidify this flabbiness, we returned to ``ambient case'' of a finite group $G$ acting on the set $X$ with restricted action of $S\in\Syl_p(G)$ to see what structure had been lost.  By examining the simultaneous action of $G$ on $X$ and $S$ (where the latter is given by conjugation, and is the same sort of ``action'' that underlies fusion theory in general), we were led to define the notion of \emph{fusion action system}, and later  \emph{$p$-local finite group actions}.

\subsection{Results}

Explicitly, given a finite group $G$ acting on the finite set $X$ and $S\in\Syl_p(G)$, the fusion action system relative to these data is the category $\fX_G$ whose objects are the subgroups of $S$ and whose morphisms are given by
\[
\fX_G(P,Q)=\left\{(\varphi,\sigma)\big|\varphi\in\Inj(P,Q),\textrm{ }\sigma\in\Sigma_X,\textrm{ and }\exists g\in G\textrm{ s.t. }(\varphi,\sigma)=(c_g|_P,\ell_g)\right\}.
\]
Here, $\Sigma_X$ denotes the symmetric group on $X$, $c_g:P\to Q$ is the injective morphism given by conjugation by $g$, and $\ell_g$ is the permutation of $X$ induced by $g$.  The $p$-local finite group action arising from this situation is $\fX_G$ together with a more complicated category that can be roughly thought of as restoring to the fusion action system the data of which elements of $G$ induce a morphism in $\fX_G$.  Much of this paper is devoted to properly abstractifying these concepts, so that we may speak of fusion action systems or $p$-local finite group actions without actually referencing the finite group $G$.  For fusion action systems, as is the case with fusion systems, the abstraction process begins with with a naive description of a category similar to $\fX_G$ and proceeds by imposing further Saturation Axioms so as to more closely mimic the properties of the ambient case.  These Axioms are listed in Definition \ref{fusion action system saturation}, then simplified in the results of Subsection \ref{simplification subsection}.

With our main object of study defined, we go on to investigate the properties of fusion action systems.  These algebraic objects generalize fusion systems, but retain many important aspects of their structure.  For instance, it is not difficult to prove a variant of Alperin's Fusion Theorem, so that a saturated fusion action system is generated by inclusions and the automorphisms of subgroups (Theorem \ref{Alperin for fusion action systems}).  Topologically, a $p$-local finite group action arising from an actual finite $G$-action on $X$ allows us to reconstruct the Borel construction $EG\times_G X$ up to $p$-completion, just as a $p$-local finite group arising from $G$ determines the homotopy type of $|BG|^\wedge_p$ (Theorem \ref{Borel construction from p-local finite group actions}).  Even without reference to an ambient group, there is an obstruction theory to the existence and uniqueness of of a $p$-local finite group action associated to a fusion action system, which mimics that of \cite{BLO2} (Theorem \ref{fusion action system obstructions to p-local finite group actions}).  

The introduction of fusion action systems also leads us to a fusion-theoretic expression of some basic ideas from group theory.  If the underlying action of $S$ on $X$ is not faithful, we construct the \emph{core fusion subsystem} that mimics the fusion of the kernel of a nonfaithful ambient group action.  We further show that the core subsystem is saturated if the fusion action system is, and that the fusion action system can be viewed as a sort of extension of the core by a finite group (Theorems \ref{core subsystem is saturated} and \ref{fusion actions as extensions}).  Moreover, if the fusion action system is transitive, there is a notion of a stabilizer fusion action subsystem of a (well chosen) point, which plays the role of the stabilizer subgroup of that point (Theorem \ref{stabilizer saturation}).  If we further give ourselves the data of a $p$-local finite group action associated to our transitive fusion action system, the stabilizer subsystem has an associated $p$-local finite group whose classifying space is homotopic to the classifying space of the original $p$-local finite group action (Theorem \ref{stabilizer classifying spaces}); this is the fusion-theoretic analogue of the basic fact that if $X$ is a transitive $G$-set with $H\leq G$ the stabilizer of a point, $EG\times_G X\simeq BH$.

More generally, the purpose of this paper is to develop a framework in which to discuss permutations of finite sets in a fusion-theoretic context.  Many of the results are geared toward showing that the structure and technical results of fusion systems or $p$-local finite groups can be translated to and make sense in our new universe of discourse.  Consequently, many of the proofs and much of the development closely mirror the existing literature.  We view the possibility of making such a translation as evidence that the definitions of fusion action system and $p$-local finite group action are the ``correct'' ones, at least in a moral sense.

\subsection{Outline of the paper}

Each of this paper's three sections follow the same basic arc:  We are interested certain $p$-local data that arise from a finite group without reference to that group, so we begin by describing these data when we are actually given a finite group.  We then attempt to define what the abstract version of these data should be, only to find that (in the cases of fusion systems and fusion action systems) the naive definition is far to broad to be of much use to us.  The solution is to take basic technical properties from the ambient case and impose these as ``Saturation Axioms'' to define the abstract object of interest.  The sections then conclude with further investigations of the properties of these abstract $p$-local versions of groups and group actions, largely with the goal of connecting our definitions with the existing literature.  More specifically:

Section \ref{literature review} is a review of the notions fusion systems and $p$-local finite groups, told using the Martino-Priddy conjecture as motivation for the introduction of the more complicated structure of a centric linking system.  Although this material is already in the literature (cf. \cite{LinckelmannIntro} and \cite{BLO2} for a more detailed discussion) and seems to be gaining wider recognition, we include our discussion as a prelude to the story we tell about fusion action systems and $p$-local finite group actions, which we structure along along the same narrative lines.

Section \ref{fusion action system section} introduces the titular subject of this paper.  Subsection \ref{ambient fusion action systems} describes the $p$-local data we wish to emulate of a finite group $G$ with Sylow subgroup $S$ acting on the finite group $X$.  Subsection \ref{abstract fusion action system subsection} first gives the abstract definition of a fusion action system, then narrows the terms of discussion by turning the basic results of Subsection \ref{ambient fusion action systems} into the Saturation Axioms for fusion action systems.  As a demonstration of the validity of these Axioms, Subsection \ref{abstract fusion action system subsection} closes with a proof of a simple version of the Alperin Fusion Theorem in the context of fusion action systems.

Fusion action systems arise as an abstraction of the $p$-local structure of a finite group acting on a finite set, but they turn out to have other interesting interpretations as well.  First, a fusion action system $\fX$ has an underlying fusion system $\cF^\fX$; this observation allows us to see fusion action systems as generalizations of fusion systems, in the sense that if the set being acted on is trivial, $\fX=\cF^\fX$, and moreover we can recover the Saturation Axioms for fusion systems from those of fusion action systems (Subsection \ref{abstract fusion action system subsection}).  Second, as the underlying $S$-action becomes less trivial, a fusion action system looks less like a fusion system and more like the transporter system of a finite group, to the point where if $X$ is a faithful $S$-set there is a canonical finite group $G$ that realizes the underlying fusion system as $\cF_S(G)$ and the fusion action system $\fX$ as the transporter system of $S$ in $G$ (Subsection \ref{faithful fusion action systems}).  Finally, every fusion action system has a \emph{core fusion subsystem}, and the fusion action system can be interpreted as an extension of this fusion system by a finite group (Subsection \ref{core subsystems}).

In Subsection \ref{simplification subsection} we address the fact that the Saturation Axioms for fusion action systems are highly redundant by proving that two seemingly weaker sets of Axioms are in fact equivalent.  These simplified Saturation Axioms allow us to check saturation with considerably less annoyance, which convenience we use in \ref{Puig fusion action subsystem}.  In this final subsection of Section \ref{fusion action system section}, we reconnect fusion action systems with ideas from Puig's original definition of Frobenius categories in \cite{PuigFrobeniusCategories}.  As motivation, we are able to prove that the stabilizer fusion action subsystem of a fully stabilized point in a transitive saturated fusion action system is itself saturated.  In fact, we prove much more, showing in effect that Puig's theorem on $K$-normalizers holds in the context of fusion action systems, which will hopefully allow for the discovery of new saturated fusion action systems in the future.

Section \ref{p-local finite group actions section} turns the focus from pure algebra to algebraic topology by introducing the notion of a classifying space of a fusion action system.  The first three subsections mirror the development of $p$-local finite group theory in \cite{BLO2}:  Subsection \ref{ambient action systems} introduces the notion of the $X$-centric linking action system induced by an ambient group $G$ and shows that this augmentation of the fusion action system contains enough data to reconstruct the homotopy type of the Borel construction $(EG\times_G X)^\wedge_p$.  Subsection \ref{classifying spaces of fusion action systems} abstractifies this structure so that we are able to talk about $X$-centric linking action systems associated to the fusion action system $\fX$ without any ambient group.  Since it is not at first obvious what the classifying space of a fusion action system should be, we reformulate the notion of linking action system in terms of data contained entirely within the fusion action system itself, namely the associated orbit category $\cO^\fX$.  The linking action system is shown to provide a solution to a certain homotopy lifting problem indexed on $\cO^\fX$.  Finally, Subsection \ref{fusion action obstruction theory} sets up the obstruction theory to the existence and uniqueness of an $X$-centric linking action system associated to a given fusion action system, completing the generalization of $p$-local finite group theory.

In Subsection \ref{Linking action systems as transporter systems} we draw further connections between linking action systems and the existing literature, specifically the abstract transporter systems of \cite{OliverVenturaTransporterSystems}.  We show that the former is actually an example of the latter, equipped with the additional data of a homomorphism from the fundamental group of the classifying space to a symmetric group.  We further show that this process can be reversed:  Given an abstract transporter system $\cT$ associated to the fusion system $\cF$ and a homomorphism $\theta:\pi_1(|\cT|)\to\Sigma_X$, we can construct a fusion action system $\fX^\theta$ with underlying fusion system $\cF$ and show that $\cT$ contains a linking action system associated to $\fX^\theta$, at least if we assume that $\cT$ is defined on enough subgroups.  Moreover, we show that $\fX^\theta$ is saturated on the subgroups witnessed by $\cT$, suggesting the need for a theorem in the style of \cite{BCGLO1} describing when such a situation would imply that $\fX^\theta$ is saturated on all subgroups.

The final Subsection, \ref{stablizer p-local finite groups}, we return to the question of stabilizers of points in transitive fusion action systems introduced in Subsection \ref{Puig fusion action subsystem}, though now we approach the question from the topological perspective of trying to describe the classifying space of the stabilizer.  It turns out that, using the results of \cite{OliverVenturaTransporterSystems}, the added structure of a linking action systems makes it much easier to prove that the stabilizer fusion subsystem is saturated.   Moreover, we get an complete $p$-local finite group for the stabilizer basically for free, and show that the homotopy type the nerve of the stabilizer centric linking system is the same as the original classifying space for the fusion action system.  We should expect this to be true from basic group theory, so perhaps this final result is best interpreted as confirmation that the definitions made in this paper function the way we would hope they do.

\section{Background on fusion systems}\label{literature review}

	In this chapter we review the basics of the theory of fusion systems and related concepts.  We introduce the classical notion of the fusion system of a finite group and Puig's abstraction of this idea, the transporter systems of Oliver-Ventura, and the centric linking systems of Broto-Levi-Oliver.  As the final portion of this paper will deal with related concepts, we shall also give a brief discussion of the Martino-Priddy Conjecture and its proof by Oliver.
	
	Throughout this paper, let $G$ be a finite group and $S$ a finite $p$-group, thought of as a Sylow $p$-subgroup of $G$ when appropriate.
		
	\subsection{Transporter and fusion systems of finite groups}
	
	We are interested in the $p$-local structure of finite groups, and furthermore would like to encode these data in a reasonably combinatorial or even category manner.  Beyond the one-point category $\cB G$, which by definition has $\cB G(*,*)=G$, we may choose to focus on a single Sylow subgroup as an important part of the data.

	\begin{definition}
		Let $G$ be a finite group and $S\in\Syl_p(G)$.  The \emph{transporter system on $S$ relative to $G$} is the category $\cT_G=\cT_S(G)$ whose objects are all subgroups $P\leq S$ and whose morphisms are given by
		\[
			\cT_G(P,Q)=N_G(P,Q):=\left\{g\in G\big|{}^gP\leq Q\right\}
		\]
		$N_G(P,Q)$ is the \emph{transporter of $P$ to $Q$ in $G$}.
	\end{definition}
	
	This definition singles out a given Sylow subgroup and plays an important role throughout this document.  However, it contains too much information, especially $p'$-data.  In fact, $\cT_G$ can be easily seen to contain exactly the information of $\cB G$, together with a choice of Sylow subgroup, by noting that $\cT_G(1)\cong G$.\footnote{We here make note of our notational convention:  Just as for a category $\cC$ we denote by $\cC(a,b)$ the set of morphisms $\Hom_\cC(a,b)$, we shall write $\cC(a)$ for the automorphism group of the object $a$.}  Indeed, the natural functor $\cB G\to \cT_G$ sending $*$ to $1$ induces a homotopy inverse to the natural functor $\cT_G\to\cB G$; in the world of topology, $|\cB G|\simeq|\cT_G|$.
	
	One way of understanding the sense in which $\cT_G$ has too much information is to note that there may be distinct elements element $g,g'\in G$ that conjugate $P$ to $Q$ but are indistinguishable from the point of view of the conjugation action on $P$.  In other words,  $c_g|_P=c_{g'}|_P$ or $g^{-1}g'\in Z_G(P)$.  Let us suppose that the conjugation action is the truly important $p$-local data.
	
		\begin{definition}\label{ambient fusion system definition}
		For $G$ a finite group and $S\in\Syl_p(G)$, the \emph{fusion system on $S$ relative to $G$} is the category $\cF_G:=\cF_S(G)$ whose objects are all subgroups $P\leq S$ and whose morphisms are given by
		\[
			\cF_G(P,Q)=\Hom_G(P,Q):=\left\{\varphi\in\Inj(P,Q)\big|\exists g\in G\textrm{ s.t. }
				\varphi=c_g|_P\right\}
		\]
		Note that we can also write $\cF_G(P,Q)=\cT_G(P,Q)/Z_G(P)$, so $\cF_G$ can be thought of as a quotient of $\cT_G$ and we have a natural projection functor $\cT_G\to\cF_G$.
	\end{definition}

	\begin{example}
	The most basic example of a Sylow inclusion $S\leq G$ is the case that the supergroup $G$ is equal to $S$ itself.  We denote the resulting fusion system by $\cF_S$, the \emph{minimal fusion system on $S$}.  Minimality in this case means that if $H$ is any finite group with $S\in\Syl_p(H)$ then $\cF_S\subseteq\cF_S(H)=\cF_H$.  The importance of this minimal example will become clear with the introduction of abstract fusion systems.
	\end{example}
		
	One of the more important properties of fusion systems is that they exhibit ``local to global'' phenomena. For example, the global condition that $s,s'\in S$ are conjugate in $G$ can be realized by simply looking at automorphism groups in the fusion system itself.
	
	\begin{theorem}[Alperin's Fusion Theorem]\label{Alperin's fusion theorem}
	For $G$ a finite group and $S\in\Syl_p(G)$, the fusion system $\cF_G$ has the property that any morphism can be written as a composite of automorphisms in $\cF_G$ of subgroups together with subgroup inclusions.
	\begin{proof}
		See \cite{AlperinFusion} for a stronger version of this result.  This result was strengthened further in \cite{GoldschimdtConjugationFamily} to the Alperin-Goldschmidt fusion theorem, which described a particular class of subgroups whose automorphisms determine the fusion system.  Finally, Puig showed in \cite{PuigStructureLocale} that the class of groups identified by Goldschmidt is truly essential in order to generate the fusion system and proved in \cite{PuigFrobeniusCategories} an abstract analogue of the Fusion Theorem that makes no reference to the finite group $G$.
	\end{proof}
	\end{theorem}
	
	\subsection{Centric linking systems of finite groups}
	Definition \ref{ambient fusion system definition} shows that the fusion system $\cF_G$ can be thought of as the quotient of the transporter system $\cT_G$ obtained by killing the action of the centralizer of the source.  This quotienting process kills both $p$- and $p'$-information; if we wish to study the $p$-local structure of $G$, perhaps we should seek a less brutal quotient as an intermediary between $\cT_G$ and $\cF_G$.
	
	To find this intermediary category, technical considerations suggest that we restrict attention to a particular collection of subgroups of $S$.  The reasons will become clear in short order. Let us therefore introduce a seemingly ad hoc definition of the class of subgroups that will be central in the following discussion:
	
	\begin{definition}\label{p-centric in ambient case}
		A $p$-subgroup $P$ of $G$ is \emph{$p$-centric} if $Z(P)\in\Syl_p(Z_G(P))$.  Equivalently, $P$ is $p$-centric if there exists a (necessarily unique) $p'$-subgroup $Z_G'(P)\leq Z_G(P)$ such that $Z_G(P)=Z(P)\times Z_G'(P)$.
	\end{definition}
	
	\begin{notation}
		We shall reserve the notation $Z_G'(P)$ for $O^p(Z_G(P))$ in the case that $P$ is $p$-centric in $G$, in which case we also have $Z_G'(P)=O_{p'}(Z_G(P))$.
	\end{notation}
	
	\begin{notation}
		By $\cT_G^c$ we mean the full subcategory of the transporter system $\cT_G$ whose objects are the $p$-centric subgroups of $S$.  We use similar notation to denote full centric subcategories of fusion systems and other related categorical versions of groups we'll encounter.
	\end{notation}
	
	We are now in the position to introduce our intermediary between $\cT_G$ and $\cF_G$:
	
	\begin{definition}
		The \emph{centric linking system} of a finite group $G$ with Sylow $S$ is the category $\cL_G^c$ whose objects are the $p$-centric subgroups of $S$ and whose morphisms are the classes
		\[
			\cL_G^c(P,Q)=N_G(P,Q)/Z_G'(P)
		\]
	\end{definition}
	
	The quotient functors $\xymatrix{\cT_G^c\ar[r]&\cL_G^c\ar[r]& \cF_G^c}$ relate our three nontrivial notions of $G$ as a category and emphasize how some information is lost at each transition.  We shall make use of the relationship between these three players in the sequel.
	
	\subsection{The Martino-Priddy Conjecture}\label{Martino-Priddy section}
	
	The classifying space functor $B:\GRP\to\TOP$ is the primary tool we use in this document for studying groups in the context of algebraic topology.  We can see $B$ as the composition of $\cB:\GRP\to\CAT$ with the geometric realization functor, which suggests that perhaps our alternate categorical versions of finite groups should be viewed as topological spaces via geometric realization.
	
	For the transporter system $\cT_G$, this works:  As $|\cT_G|\simeq BG$, essentially all the algebraic information of $\cT_G$ is realized topologically in this manner.  Roughly speaking, the transporter system has a minimal object that includes in every other one.  One can see that in such a situation, all the topological data we could hope to extract from the transporter system is in some sense concentrated in the automorphisms of this minimal object; this is a situation that is explored in greater depth in \cite{OliverVenturaTransporterSystems}.
	
	However, simply taking the nerve of the fusion system $\cF_G$ will not yield an interesting space:  The object $1$ is initial in $\cF_G$, so $|\cF_G|$ is contractible.  Indeed, as fusion systems are naturally equipped with a collection of ``inclusion morphisms''---honest inclusions of subgroups in this case---this fact could be seen as a special case of the general reason why $|\cT_G|\simeq BG$.  We will have to be more clever about how we construct a topological space from a fusion system if we are to arrive at anything interesting.
	
	We therefore turn to the centric linking system $\cL_G^c$ as an intermediary between transporter systems and fusion systems.  The space $|\cL_G^c|$ should be related to $BG$ in some way, but as information is lost from the transition from transporter system to linking system it is unreasonable to expect that $|\cL_G^c|\simeq BG$.
	
	\begin{example}\label{stupid p'-data}
		Let $G$ be your favorite finite group, $S\in\Syl_p(G)$, and $H$ your favorite finite $p'$-group.  Then $H$ is a $p'$-subgroup of $Z_{G\times H}(S)$ so for any $P\leq S$ we have $H\leq Z_{G\times H}'(P)$.  We conclude $\cL_G^c\cong \cL_{G\times H}^c$ (actual isomorphism of categories).  It easily follows that the same result applies to fusion systems:  $\cF_G=~\cF_{G\times H}$ (equality of categories).
	\end{example}
	
	Since we construct the linking system by killing certain $\fp'$-primary data of the transporter system, we should look for an operation on topological spaces that ``isolates $p$-information'' in some appropriate sense.
	
		\begin{notation}
		Let $(-)^\wedge_p:\TOP\to\TOP$ denote the Bousfield-Kan $p$-completion functor of \cite{BousfieldKanBook}.  There is a natural transformation $\eta:\id_{\TOP}\Rightarrow(-)^\wedge_p$; for a space $\cX$, let $\eta_\cX:\cX\to\cX{}^\wedge_p$ denote the resulting $p$-completion map.
	\end{notation}
	
	\begin{definition}
		A space $\cX$ is \emph{$p$-complete} if the $p$-completion map $\eta_\cX:\cX\to\cX{}^\wedge_p$ is a homotopy equivalence.  $\cX$ is \emph{$p$-good} if $\cX{}^\wedge_p$ is $p$-complete.
	\end{definition}  

	The interested reader should refer to \cite{BousfieldKanBook} for further properties of the $p$-completion functor.
	
	So, even though it is unreasonable to expect that $|\cL_G^c|$ is homotopy equivalent to $BG$, it is conceivable that these spaces should be equivalent up to $p$-completion.  The following two theorems of \cite{BLO1} describe the relationship of the centric linking system of $G$ to the space $BG^\wedge_p$:
	
	\begin{theorem}
		For $G$ a finite group, the natural functors $\xymatrix{\cB G&\cT_G^c\ar[l]\ar[r]&\cL_G^c}$ induce mod-$p$ cohomology isomorphisms on realization.  In particular, $|\cL_G^c|^\wedge_p\simeq BG^\wedge_p$.
	\end{theorem}
	
	\begin{theorem}[Weak Martino-Priddy Conjecture]
		For finite groups $G$ and $H$, $BG^\wedge_p\simeq BH^\wedge_p$ if and only if the categories $\cL_G^c$ and $\cL_H^c$ are equivalent.
	\end{theorem}
	
	This still leaves the question of what space we should associate to the fusion system itself.  In fact, the surprising result is that the passage from transporter system to fusion system loses no topological information, and so in fact the space associated to $\cF_G$ should again be $BG^\wedge_p$.
	
	\begin{theorem}\label{Martino-Priddy conj}[Martino-Priddy Conjecture]
		The finite groups $G$ and $H$ have homotopic $p$-completed classifying spaces if and only if the $p$-fusion data of $G$ and $H$ are the same.	
		\begin{proof}
		The ``topology implies algebra'' direction is given in \cite{MartinoPriddyConjecture}.  The ``algebra implies topology'' direction was proved by Oliver in \cite{limZodd,limZ2}, using the machinery of \cite{BLO2} and the Classification Theorem of Finite Simple Groups.
		\end{proof}
	\end{theorem}
	
	That ``the $p$-fusion data of $G$ and $H$ are the same'' means that there is an isomorphism of fusion systems $\cF_G\cong\cF_H$.  The notion of isomorphism of fusion system is much stronger than saying that these categories are equivalent, or even isomorphic as categories:  Such a notion would record only the shape of the fusion system as a diagram without giving due deference to the structure of the objects of the fusion system.
	
	\begin{definition}
		Let $G$ and $H$ be finite groups with respective Sylows $S$ and $T$.  An isomorphism $\alpha:S\to T$ is a \emph{fusion preserving isomorphism} if for every $P,Q\leq S$ and $\beta\in\Hom(P,Q)$, $\beta\in\Hom_G(P,Q)$ if and only if $\alpha\beta\alpha^{-1}\in\Hom_H(\alpha P,\alpha Q)$.  In this case the fusion systems $\cF_G$ and $\cF_H$ are \emph{isomorphic as fusion systems}.
	\end{definition}

	So we claim that the data of the fusion system $\cF_G$ determine the $p$-completed homotopy type of $BG$, which returns us to the question of exactly how to associate a topological space to a fusion system.  We have already seen that simply taking the nerve of  $\cF_G$ yields nothing interesting.  This should not be surprising, as taking the geometric realization of the fusion system only records the shape of the category as a diagram without taking into account the fact that it is a diagram \emph{in $p$-groups}.  This is not a problem for either the transporter or linking systems of $G$, as for any $P\leq S$ we have $P\leq N_G(P)$ and therefore there is a natural way to identify $P$ with a subgroup of its automorphism group.  Such is not the case for fusion systems, so we must try a little harder to recover this information.  In so doing, we introduce yet another categorical version of the finite group $G$.
		
	\begin{definition}
	 The \emph{orbit category}\footnote{This notion is not to be confused with the category of orbits of $G$, whose objects are the transitive $G$-sets and where morphisms are maps of $G$-sets.  Although there is a relationship between these two notions, we will not make use of the category of $G$-orbits in this document.} of $\cF_G$ is the category $\cO_G:=\cO(\cF_G)$ whose objects are the subgroups of $S$ and whose morphisms are given by
	 \[
	 \cO_G(P,Q)=Q\backslash\cF_G(P,Q)
	 \]
	In other words, the hom-set from $P$ to $Q$ is the set orbits of the $Q$-action of $\cF_G(P,Q)$ given by postcomposition by $c_q$.
	
	$\cO_G^c$ will denote the full subcategory of $\cO_G$ whose objects are the $\cF$-centric subgroups of $S$.
	\end{definition}
	
	The functor $B-:\cF_G\to\TOP$ does not descend to a functor $\cO_G\to\TOP$, but because $\cO_G$ is defined by quotienting out inner automorphisms, it is easy to see that there is a homotopy functor $\overline B-:\cO_G\to ho\TOP$.  If we could find a homotopy lifting
	\[
	\xymatrix{
		&	\TOP\ar[d]\\
		\cO_G\ar@{-->}^{\widetilde B-}[ur]\ar[r]_-{\overline B-}&	ho\TOP
	}
	\]
	we could consider $\hocolim_{\cO_G}\widetilde B-$, and relate this space to $BG^\wedge_p$.  The following Proposition explains this relationship.
	
	\begin{prop}
		For $G$ a finite group, consider the diagram
		\[
		\xymatrix{
			\cL_G^c\ar[r]^-{*}\ar[d]_\pi&\TOP\\
			\cO_G^c\ar@{-->}[ru]_-L
		}
		\]
		where $\pi$ is the composite of the natural quotients $\cL_G^c\to\cF_G^c\to\cO_G^c$ and $L$ is the left homotopy Kan extension of the trivial functor $*$ over $\pi$.  Then $L$ is a homotopy lifting of $\overline B-:\cO_G^c\to ho\TOP$, and in particular we have
		\[
		\hocolim_{\cO_G^{c}} L\simeq\hocolim_{\cL_G^c}*=|\cL_G^c|\simeq_p BG
		\]
		\begin{proof}
			\cite{BLO2}.
		\end{proof}
	\end{prop}

	So far we have simply restated the original question of whether topological information is lost on the transition from linking system to fusion system:  If we have a linking system in mind for $\cF_G$, there is a homotopy lifting of $\overline B-$, which allows us to construct our desired space from the fusion system.  But what if there is another finite group $H$ such that $\cF_G=\cF_H$ and yet $\cL_G^c\neq\cL_H^c$: Is it possible there are two distinct homotopy liftings and thus two different spaces associated to $\cF_G$?  Or can this never happen?  How do we approach this problem?
	
	\subsection{Group theory without groups}\label{group theory without groups}
	
	The basic problem introduced at the end of Subsection \ref{Martino-Priddy section} is the need to think of fusion and linking systems as algebraic objects distinct from the finite groups from which they came. 
		
	Puig provided the necessary insight and abstraction to codify this generalization.  Here we introduce his idea of abstract fusion systems (or ``Frobenius categories'' in the terminology of \cite{PuigFrobeniusCategories}), though we shall use the language of Broto-Levi-Oliver.  We also review the abstraction of the notion of centric linking system, due to \cite{BLO2}.
		
	\begin{definition}
		Let $S$ be a $p$-group.  An \emph{abstract fusion system on S} is a category $\cF$ whose objects are all subgroups $P\leq S$ and whose morphisms are some collection of injective group maps:  $\cF(P,Q)\subseteq\Inj(P,Q)$.  We require that the following conditions be satisfied:
		\begin{itemize}
		\item{\bf ($S$-conjugacy)}  The minimal fusion system $\cF_S$ is a subcategory of $\cF$.  
		\item{\bf (Divisibility)}  Every morphism of $\cF$ factors as an isomorphism of groups followed by an inclusion.
		\end{itemize}
		Composition of morphisms is composition of group maps.
	\end{definition}
	
	This is a very simple definition.  In fact, it is perhaps too simple to be useful:  This mimics the situation where $S$ is a $p$-subgroup of some unnamed ambient group, but not where $S$ is a \emph{Sylow} $p$-subgroup.  There is a great deal of additional structure that comes from such a Sylow inclusion; the question is how to codify these interesting data without reference to an ambient group.  This will lead to the addition of \emph{Saturation Axioms} that must be imposed on a fusion system.
	
	\begin{definition}
	We will need the following terms to state the Saturation Axioms:	
		\begin{itemize}
		\item $P\leq S$ is \emph{fully normalized in $\cF$} if $|N_S(P)|\geq |N_S(Q)|$ for all $Q\cong_\cF P$.
		\item $P\leq S$ is \emph{fully centralized in $\cF$} if $|Z_S(P)|\geq |Z_S(Q)|$ for all $Q\cong_\cF P$.
		\item  For any $\varphi\in\cF(P,Q))_\iso$, let $N_\varphi\leq N_S(P)$ denote the group
		\begin{eqnarray*}
			N_\varphi&=&\left\{n\in N_S(P)\big|\varphi\circ c_n\circ\varphi^{-1}\in\Aut_S(Q)\right\}\\
				&=&\left\{n\in N_S(P)\big|\exists s\in S\textrm{ s.t. }\forall p\in P,\varphi({}^np)={}^s\varphi(p)\right\}
		\end{eqnarray*}
		\end{itemize}
		$N_\varphi$ will be called the \emph{extender} of $\varphi$.
	\end{definition}
	
	\begin{remark}
		Perhaps some motivation for these concepts is in order.  Each of these definitions comes from the idea that there are certain ``global'' phenomena that can be captured purely through local, fusion-theoretic data of a group.  For instance, if there is an ambient Sylow $G$ giving rise to the fusion system, then $P\leq S$ is fully normalized in $\cF$ if and only if $N_S(P)\in\Syl_p(N_G(P))$, and similarly for the concept of full centralization.
		
		The motivation for the extender $N_\varphi$ comes from Alperin's Fusion Theorem \ref{Alperin's fusion theorem}.  If we wish that the morphisms of a fusion system be generated by inclusions and automorphisms of subgroups, there must be some way of extending certain morphisms between different subgroups within the fusion system.  The extender is the maximal subgroup of $N_S(P)$ to which we could hope to extend $\varphi\in\cF(P,Q)$, so the question becomes when we can achieve this maximal extension.  
	\end{remark}
	
	\begin{definition}[Saturation Axioms for Abstract Fusion Systems]
		The fusion system $\cF$ is \emph{saturated} if
		\begin{itemize}
		\item Whenever $P$ is fully $\cF$-normalized, $P$ is fully $\cF$-centralized.
		\item Whenever $P$ is fully $\cF$-normalized, $\Aut_S(P)\in\Syl_p\left(\cF(P)\right)$.
		\item If $Q$ is fully $\cF$-centralized and $\varphi\in\Iso_\cF(P,Q)$, then there is some morphism $\widetilde\varphi\in\cF\left(N_\varphi,S\right)$ that extends $\varphi$:  $\widetilde\varphi|_P=\varphi$.
		\end{itemize}
		The first two conditions will be referred to as the \emph{Sylow Axioms} and the third the \emph{Extension Axiom}.
	\end{definition}
	
	Note that the Saturation Axioms are somewhat redundant, in that seemingly less restrictive Axioms turn out to be equivalent.  Examples of such simplified axiom sets can be found in \cite{OnofreiStancuCharacteristicSubgroup} and \cite{RobertsShpectorovSaturationAxiom}.  See also \cite{PuigFrobeniusCategories} for a different take on the nature of saturation.
	
	Saturated fusion systems form a category, though we shall not make great use of the notion of morphism of fusion system in this paper.
	
	Not only can we discuss fusion systems without reference to an ambient group, but there also exists a notion of abstract transporter systems that generalizes both the constructions $\cT_G$ and $\cL_G^c$.
	
	\begin{definition}[\cite{OliverVenturaTransporterSystems}]\label{abstract transporter system}
		Let $\cF$ be an abstract fusion system on $S$.  An \emph{abstract transporter system associated to $\cF$} is a category $\cT$ whose objects are some set of subgroups $P\leq S$ that is closed under $\cF$-conjugacy and overgroups, together with functors
		\[
		\xymatrix
		{
		\cT_S^{\Ob(\cT)}(S)\ar[rr]^-\delta&&	\cT\ar[rr]^\pi&&	\cF
		}
		\]
		For any $s\in N_S(P,Q)$, set $\widehat s=\delta_{P,Q}(s)\in\cT(P,Q)$. Similarly, for any $\fg\in\cT(P,Q)$, set $c_\fg=\pi_{P,Q}(\fg)\in\cF(P,Q)$.
		
		The following Axioms apply:
		\begin{itemize}
		\item[(A1)]  On objects, $\delta$ is the identity and $\pi$ is the inclusion.
		\item[(A2)]  For any $P\in\Ob(\cT)$, define
		\[
			E(P)=\ker\left[\pi_{P,P}:\cT(P)\to\cF(P)\right]
		\]
		Then for any $P,Q\in\Ob(\cT)$, the group $E(P)$ acts right-freely and $E(Q)$ acts left-freely on $\cT(P,Q)$.  Moreover, the map
		$
			\pi_{P,Q}:\cT(P,Q)\to\cF(P,Q)
		$
		is the orbit map of the $E(P)$-action.
		\item[(B)]  The functor $\delta$ is injective on morphisms, and for all $s\in N_S(P,Q)$ we have $c_{\widehat s}=c_s\in\cF(P,Q)$.
		\item[(C)]  For all $\fg\in\cT(P,Q)$ and all $p\in P$, the following diagram commutes in $\cT$:
		\[
		\xymatrix
		{
		P\ar[r]^\fg\ar[d]_{\widehat p}&	Q\ar[d]^{\widehat{c_{\fg}(p)}}\\
		Q\ar[r]_\fg&	Q
		}
		\]
		\end{itemize}
		{\bf Saturation axioms:}
		\begin{itemize}
		\item[(I)]  $\delta_{S,S}(S)\in\Syl_p(\cT(S))$.
		\item[(II)]  For all $\fg\in\Iso_\cT(P,Q)$ and supergroups $\widetilde P\trianglerighteq P$ and $\widetilde Q\trianglerighteq Q$ such that
		$
		\fg\circ\delta_{P,P}\left(\widetilde P\right)\circ\fg^{-1}\leq \delta_{Q,Q}\left(\widetilde Q\right)
		$,
		there is a morphism $\widetilde\fg\in\cT(\widetilde P,\widetilde Q)$ that satisfies $\widetilde\fg\circ\widehat 1_P^{\widetilde P}=\widehat 1_Q^{\widetilde Q}\circ\fg\in\cT(P,\widetilde Q)$.
		\end{itemize}
	\end{definition}
	
	\begin{notation}
		As was noted in the definition of transporter systems, for any $s\in N_S(P,Q)$ we set 
		\[
			\widehat s\big|_P^Q=\delta_{P,Q}(s)\in\cT(P,Q)
		\]
		and if the source and target are obvious from the context we shall simply write $\widehat s$.  Similarly, for any $R\leq N_S(P)$ we denote by $\widehat R\big|_P^P\leq \cT(P)$ the group $\delta_{P,P}(R)$, again writing simply $\widehat R$ if there is no chance of confusion.
	\end{notation}
	
	\begin{remark}
		$E(P)$ always contains $\delta_{P,P}(Z(P))$ by Axiom (B).  As $\delta$ is injective, we identify $Z(P)$ with its image in $\cT(P)$.  
	\end{remark}
	
	We say that $P\leq S$ is \emph{$\cF$-centric} if $Z_S(Q)=Z(Q)$ for all $Q$ that are $\cF$-conjugate to $P$.  The collection of $\cF$-centric subgroups is of central importance to the study of the homotopy theory of fusion systems.
	
	\begin{definition}\label{Linking systems as minimal transporter systems}
		Let $\cF$ be a saturated fusion system on $S$.  A transporter system $\cL$ associated to $\cF$ is an \emph{abstract centric linking system} if:
		\begin{itemize}
			\item $\Ob(\cL)$ is the collection of $\cF$-centric subgroups of $S$.
			\item For every $P\in\Ob(\cL)$, $E(P)=Z(P)$.
		\end{itemize}
		Thus an abstract linking system can be thought of as a minimal transporter system on the $\cF$-centric subgroups of $S$.
	\end{definition}
	
	The rationale for augmenting a fusion system with the more complicated structure of a centric linking system is that they lead to a proof the the Martino-Priddy Conjecture:  We already know that the centric linking system $\cL_G^c$ arising from a finite group $G$ determines the homotopy type of $BG^\wedge_p$, and using the obstruction theory of \cite{BLO2} and the Classification of Finite Simple Groups, Oliver was able to show in \cite{limZodd,limZ2} that there is a unique centric linking system associated to any fusion system of a finite group.

\section{Fusion action systems}\label{fusion action system section}

\subsection{The ambient case}\label{ambient fusion action systems}

	To motivate the following discussion, let us identify $ p$-local data of a finite group acting on a finite set that we are trying to abstract.  Let $G$ be a finite group, $X$ a $G$-set, and $S\in\Syl_ p(G)$.  We use the shorthand $\cF_G=\cF_S(G)$ and $\cF_S=\cF_S(S)$.  We additionally denote the action map defining the $G$-set structure on $X$ by $\ell$, and will write $\ell_g$ or $\ell_H$ for the permutation of $X$ represented by $g$ or the subgroup of $\Sigma_X$ that is the image of $H\leq G$, respectively.
	
	 \begin{definition}
	 The \emph{fusion action system on $S$ relative to the $G$-action on $X$} is the category  $\fX_G:=\fX_S(G)$ whose objects are the subgroups $P\leq S$ and whose morphisms are given by
		\[
			\fX_G(P,Q)=\left\{(\varphi,\sigma)\in\Hom(P,Q)\times\Sigma_X\big|\exists g\in G\textrm{ s.t. }\varphi=c_g|_P\in\Inj(P,Q)\textrm{ and } \sigma=\ell_g:{}_P X\cong{}_P^\varphi X\right\}.
		\]
	Any pair $(\varphi,\sigma)$ that appears as a morphism in $\fX_G(P,Q)$ has the property that for all $p\in P$ and $x\in X$, $\sigma(p\cdot x)=\varphi(p)\cdot\sigma(x)$, or equivalently that $\sigma\ell_p\sigma^{-1}=\ell_{\varphi(p)}$; such a pair of an injection of subgroups and permutation of $X$ is \emph{intertwined}.
	\end{definition}
	
	\begin{example}
	The restriction of the $G$-action to $S$ defines a fusion action system $\fX_S$, the \emph{minimal fusion action system} of a given $S$-action on $X$.
	\end{example}
	
	The fusion action system $\fX_G$ has an \emph{underlying fusion system} on $S$, obtained by simply ignoring the second coordinates of the morphisms of $\fX_G$.  This is of course $\cF_G$. 
	
	 If we ignore the first coordinates of $\fX_G$, we obtain another interesting algebraic structure, which turns out to be, effectively, a transporter system on a quotient of $S$.  Then $\ell_S\leq\Syl_ p(\ell_G)$, and we can talk about the transporter system on $\ell_S$ relative to this inclusion, $\cT_{\ell_G}:=\cT_{\ell_S}(\ell_G)$.  Then the functor $\fX_G\to\cT_{\ell_G}$ given on objects by $P\mapsto\ell_P$ and on morphisms by projection onto the second factor is surjective, in the sense that any morphism of $\cT_{\ell_G}$ lies in the image of a morphism of $\fX_G$.

	\begin{definition}
		For any pair of subgroups $P,Q\leq S$, let $\pi_\Sigma^{P,Q}:\fX_G(P,Q)\to \Sigma_X$ be the set-map projection onto the second coordinate.  We shall suppress the reference to the source and target groups and simply write $\pi_\Sigma$ for this map. When $P=Q$, the set-map $\pi_\Sigma$ is a homomorphism of groups.  In this case let $\Sigma_\fX^G(P)$ be the image of $\fX_G(P)$ under $\pi_\Sigma$,  and similarly set $\Sigma_\fX^S(P)$ to be the image of $\fX_S(P)$ under $\pi_\Sigma$.		
	\end{definition}
	
	To emphasize the connection with fusion systems, we introduce the following notation:
	
	\begin{notation}	
	Let $\Aut_G(P;X)$ denote the group $\fX_G(P)$.  This group is a simultaneous  generalization of $\Aut_G(P)=\cF_G(P)$ and the group $\Aut_G(X)$ of $G$-set automorphisms of $X$.  Similarly define $\Aut_S(P;X)$.
	\end{notation}
	
	There are several groups of automorphisms that arise from inspection of the category $\fX_G$, two of which are $\Aut_G(P;X)$ and $\Aut_S(P;X)$.  We need some notation to introduce the others.
	
	\begin{notation}
	Let $\widehat C$ be the \emph{core}, or kernel, of the $G$ action on $X$.  We define the \emph{$X$-normalizer} and \emph{$X$-centralizer in $G$} of a subgroup $H\leq G$ to be
	\[
		N_G(H;X):=N_G(H)\cap\widehat C\qquad\textrm{and}\qquad Z_G(H;X):=Z_G(H)\cap\widehat C.
	\]
	Similarly, let $C$ be the core of the $S$-action on $_SX$, so $C=\widehat C\cap S$.  The \emph{$X$-normalizer} and \emph{$X$-centralizer in S} of $P\leq S$ are then
	\[
	N_S(P;X):=N_S(P)\cap K\qquad\textrm{and}\qquad Z_S(P;X):=Z_S(P)\cap C.
	\]
	Note that $N_G(H;X)$ is just another name for $N_{\widehat C}(H)$.  We use this notation to emphasize the idea that $G$ is acting simultaneously on its subgroups (by conjugation) and on $X$ (by left multiplication).
	\end{notation}
	
	\begin{definition}\label{diamond group quotients}
	For any $P\leq S$, we have the inclusions of the groups $Z_G(P;X)$, $N_G(P;X)$, $Z_G(P)$, and $N_G(P)$ as depicted in Figure \ref{diamond of groups}.  All of these inclusions are normal, so we name to the respective quotients in that Figure as well.  The term \emph{$G$-automizers of $P$} will refer to any of these quotient groups.

	\begin{figure}[h!]
	\[
	\xymatrix{
		&			&N_G(P)\ar@{-}[ddd]|*+[F.]{\fX_G(P)}\ar@{-}[ddll]|*+[F.]{\cF_G(P)}\ar@{-}[drr]|*+[F.]{\Sigma_\fX^G(P)}\\
		&			&		&		&N_G(P;X)\ar@{-}[ddll]|*+[F.]{\cF_G(P)_0}\\
		Z_G(P)\ar@{-}[drr]|*+[F.]{\Sigma_\fX^G(P)_0}\\
		&		&Z_G(P;X)
	}
	\]
	\caption{Naming the $G$-automizer groups}\label{diamond of groups}
	\end{figure}
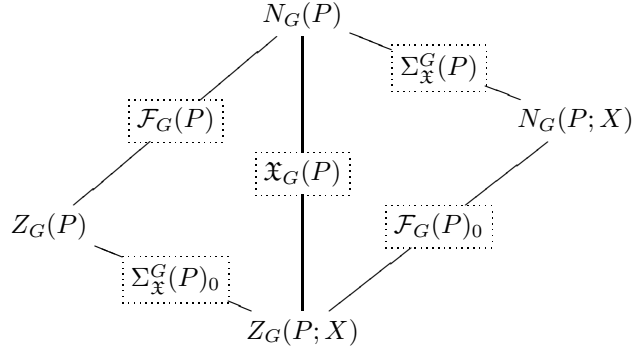
	
	We can expand on the definitions of these $G$-automizers as follows:
	\[
	\begin{array}{rlcll}
	\fX_G(P)&=&N_G(P)/Z_G(P;X)&=&\Aut_G(P;X)\\
	\cF_G(P)&=&N_G(P)/Z_G(P)&=&\left\{\varphi\in\Aut(P)\big|(\varphi,\sigma)\in\fX_G(P)\right\}\\
	\Sigma_\fX^G(P)&=&N_G(P)/N_G(P;X)&=&\left\{\sigma\in\Sigma_X\big|(\varphi,\sigma)\in\fX_G(P)\right\}\\
	\cF_G(P)_0&=&N_G(P;X)/Z_G(P;X)&=&\left\{\varphi\in\Aut(P)\big|(\varphi,\mathrm{id}_X)\in\fX_G(P)\right\}\\
	\Sigma_\fX^G(P)_0&=&Z_G(P)/Z_G(P;X)&=&\left\{\sigma\in\Sigma_X\big|(\mathrm{id}_P,\sigma)\in\fX_G(P)\right\}
	\end{array}
	\]
	We also have the short exact sequences
	\[
	\xymatrix@R=0pt{
	1\ar[r]&	\Sigma_\fX^G(P)_0\ar[r]&	\fX_G(P)\ar[r]&	\cF_G(P)\ar[r]&	1,\\
	1\ar[r]&	\cF_G(P)_0\ar[r]&	\fX_G(P)\ar[r]&	\Sigma_\fX^G(P)\ar[r]&	1.
	}
	\]	
	Clearly $\Sigma_\fX^G(P)_0$ can be identified with a subgroup of $\Aut_P(X)$ and $\cF_G(P)_0$ with a subgroup of $\cF_G(P)$.  If $\varphi\in\cF_G(P)_0$, the identity map defines an isomorphism of $P$-sets $\mathrm{id}_X:{}_PX\cong{}_P^\varphi X$, or $\ell_p=\ell_{\varphi(p)}$ for all $p\in P$.  Thus $\varphi(p)=p$ mod $\widehat C$, so we have $\cF_G(P)_0\leq\ker\left(\cF_G(P)\to\cF_{\ell_{G}}(\ell_{P})\right)$.
	\end{definition}
		
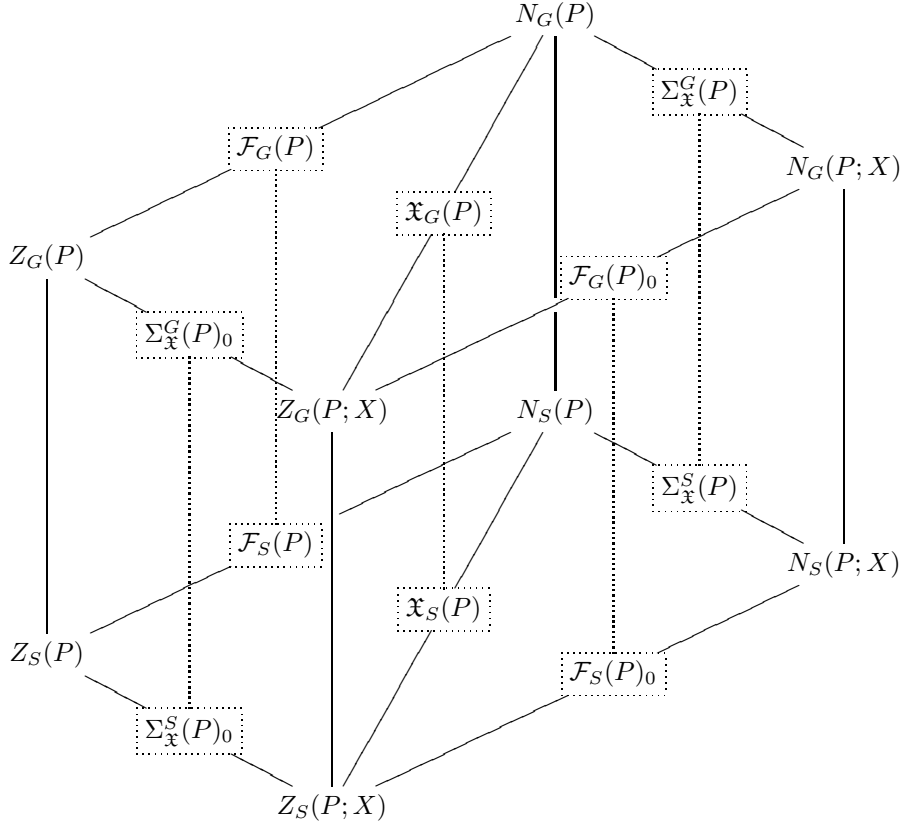
\begin{figure}[h]
	\[
	\xymatrix@=18pt{
	&&&&&N_G(P)\ar@{-}[ddddd]|!{[dddddll];[ddrrr]}\hole\\
	\\
	&&&&&&&&N_G(P;X)\ar@{-}[llluu]|*+[F.]{\Sigma_\fX^G(P)}="1a"\\
	Z_G(P)\ar@{-}[rrrrruuu]|(.45)*+[F.]{\cF_G(P)}="2a"\\
	\\
	&&&Z_G(P;X)\ar@{-}[llluu]|*+[F.]{\Sigma_\fX^G(P)_0}="3a"\ar@{-}[rrrrruuu]|(.55)*+[F.]{\cF_G(P)_0}="4a"\ar@{-}[rruuuuu]|*+[F.]{\fX_G(P)}="5a"&&N_S(P)\ar@{-}[uuuuu]|!{[ll];[rrruuu]}\hole\\
	\\
	&&&&&&&&N_S(P;X)\ar@{-}[uuuuu]\ar@{-}[llluu]|*+[F.]{\Sigma_\fX^S(P)}="1b"\\
	Z_S(P)\ar@{-}[uuuuu]\ar@{-}[rrrrruuu]|(.45)*+[F.]{\cF_S(P)}="2b"|!{[rrrdd];[rrruuu]}\hole\\
	\\
	&&&Z_S(P;X)\ar@{-}[uuuuu]\ar@{-}[llluu]|*+[F.]{\Sigma_\fX^S(P)_0}="3b"\ar@{-}[rrrrruuu]|(.55)*+[F.]{\cF_S(P)_0}="4b"\ar@{-}[rruuuuu]|*+[F.]{\fX_S(P)}="5b"
	\ar @{..}"1b";"1a"
	\ar @{..}"2b";"2a"
	\ar @{..}"3b";"3a"
	\ar @{..}"4b";"4a"
	\ar @{..}"5b";"5a"
	}
	\]
	\caption{Comparing $S$-automizers to $G$-automizers}\label{S vs. G}
	\end{figure}

	Relative to $S$ we have the same relationships amongst the groups $Z_S(P;X)$, $Z_S(P)$, $N_S(P;X)$, and $N_S(P)$.  We also have the inclusions $Z_S(P;X)\leq Z_G(P;X)$, etc.  These data give rise to the rather more complicated diagram of Figure \ref{S vs. G}.
	
	This diagram allows us to relate the ``minimal'' automizer groups---those that arise from the simultaneous $S$-action on its subgroups and $X$---with those that arise from $G$.  For certain subgroups $P\leq S$, there is a particularly nice relationship; to express it, we will need the following terminology:
	\\\\
	\begin{definition}
	Given the fusion action system $\fX_G$ and $P\leq S$, we say that
	\begin{itemize}
	\item $P$ is \emph{fully normalized relative to $G$} if $N_S(P)\in\Syl_ p\left(N_G(P)\right)$.
	\item $P$ is \emph{fully centralized relative to $G$} if $Z_S(P)\in\Syl_ p\left(Z_G(P)\right)$.
	\item $P$ is \emph{fully $X$-normalized relative to $G$} if $N_S(P;X)\in\Syl_ p\left(N_G(P;X)\right)$.
	\item $P$ is \emph{fully $X$-centralized relative to $G$} if $Z_S(P;X)\in\Syl_ p\left(Z_G(P;X)\right)$.
	\end{itemize}
	The reference to $G$ will be omitted if it is clear from the context.	
	\end{definition}
	
	Just as in the case of ordinary fusion systems, we can think of these terms in the following manner:  Instead of starting out with a chosen Sylow subgroup $S\leq G$ and a $ p$-subgroup $P\leq S$, we want to investigate the $ p$-subgroup $P$ on its own terms.  For instance, to understand the $ p$-part of $N_G(P)$ we must pick a ``right'' Sylow of $G$:  Such a Sylow must contain not just $P$, but also a Sylow subgroup of $N_G(P)$.  Saying that $P$ is fully normalized means that we have made this choice correctly within the $G$-conjugacy class of $P$ inside $S$.  Moreover, Sylow's theorems show that for a fixed Sylow $S$ and $P\leq S$, there is a subgroup $Q\leq S$ such that $Q$ is $G$-conjugate to $P$ and $Q$ is fully normalized, etc.
	
	We now use this terminology to describe the relationships between the $S$- and $G$-automizers of $P$.
	
	\begin{prop}\label{ambient Sylow conditions}
	Fix $\fX_G$ and $P\leq S$.
	\begin{enumerate}
		\item[(1)]  If $P$ is fully normalized relative to $G$, then 
		\begin{itemize}
			\item $\Aut_S(P)\in\Syl_ p\left(\cF_G(P)\right)$
			\item $\Aut_S(P;X)\in\Syl_ p\left(\fX_G(P)\right)$
			\item $\Sigma_\fX^S(P)\in\Syl_ p\left(\Sigma_\fX^G(P)\right)$
		\end{itemize}
		and furthermore $P$ is fully centralized, $X$-normalized, and $X$-centralized relative to $G$.
		\item[(2)]  If $P$ is fully $X$-normalized relative to $G$, then $\cF_S(P)_0\in\Syl_ p\left(\cF_G(P)_0\right)$ and $P$ is fully $X$-centralized relative to $G$.
		\item[(3)]  If $P$ is fully centralized relative to $G$, then $\Sigma_\fX^S(P)_0\in\Syl_ p\left(\Sigma_\fX^G(P)_0\right)$ and $P$ is fully $X$-centralized relative to $G$.
		\end{enumerate}
	\begin{proof}
		If $P$ is fully normalized, consider the diagram of short exact sequences:
		\[
		\xymatrix
		{
		1\ar[r]&	Z_G(P;X)\ar[r]&		N_G(P;X)\ar[r]&		\fX_G(P)\ar[r]&	1\\
		1\ar[r]&	Z_S(P;X)\ar[r]\ar@{^(->}[u]&		N_S(P;X)\ar[r]\ar@{^(->}[u]_{\textrm{Sylow}}&		
				\Aut_S(P;X)\ar[r]\ar@{^(->}[u]&	1
		}
		\]
		The middle inclusion's being Sylow implies that the other two are as well, and thus $\Aut_S(P;X)\in\Syl_ p\left(\fX_G(P)\right)$ and $P$ is fully $X$-centralized.  The other claims are proved in a similar manner.
	\end{proof}
	\end{prop} 
	
	Lemma \ref{ambient Sylow conditions} tells us certain properties of $\fX_G$ that follow from a subgroup's having the property that the normalizer, $X$-normalizer, or centralizer of that subgroup in $S$ is Sylow in the respective group of $G$.  It gives no information, however, about the case when $P$ is fully $X$-centralized relative to $G$, which turns out to be very important for understanding extensions of morphisms in the fusion action system.  Recall that if $g\in G$ is such that $^gP\leq S$, then $g$ determines the morphism $(c_g,\ell_g)\in\fX_G(P,{}^gP)$.
	
	\begin{definition}
		For $P\leq S$ and $g\in G$ such that $^gP\leq S$, define
		\[
		N_{(c_g,\ell_g)}=\left\{n\in N_S(P)\big|\exists s\in N_S({}^gP)\textrm{ s.t. } (c_{gng^{-1}}, \ell_{gng^{-1}})=(c_s,\ell_s)\right\}.
		\]
		The group $N_{(c_g,\ell_g)}$ is the \emph{extender} of $(c_g,\ell_g)$ in the fusion action system.
	\end{definition}
	
	\begin{remark}\label{obvious extension facts}
	Consider the following easy observations:
	\begin{itemize}
	\item $P\cdot Z_S(P;X)\leq N_{(c_g,\ell_g)}\leq N_S(P)$.  
	
	\item We define a notion of ``translation'' along morphisms:  Given $g\in G$ and $P\leq S$ such that $^gP\leq S$, set
	\[
		t_{(c_g,\ell_g)}:\fX_G(P)\to \fX_G({}^gP):
		(\varphi,\sigma)\mapsto(c_g\varphi c_g^{-1},\ell_g\sigma\ell_g^{-1}).
	\]
	We could then define $N_{(c_g,\ell_g)}$ to be the subgroup of $n\in N_S(P)$ such that
	\[
	(c_n,\ell_n)\in t_{(c_g,\ell_g)}^{-1}\left(\Aut_S({}^gP;X)\right)
	\]
	or even more confusingly, the preimage in $N_S(P)$ of the preimage in $\fX_G(P)$ of $\Aut_S({}^gP;X)$.
		
	\item The extender $N_{(c_n,\ell_n)}$ is the largest subgroup of $N_S(P)$ to which we could hope to extend $(c_g,\ell_g)$.  Here, extension of a morphism means that we would find a $g'\in G$ such that $^{g'}N_{(c_g,\ell_g)}\leq S$, $\ell_{g'}=\ell_g$, and $c_{g'}=c_g$ on $P$.  The last two conditions are equivalent to requiring that $g'Z_G(P;X)=gZ_G(P;X)$.  In this case we say that $(c_{g'},\ell_{g'})=(c_{g'},\ell_g)\in\fX_G\left(N_{(c_g,\ell_g)},S\right)$ is an extension of $(c_g,\ell_g)\in\fX_G(P,S)$.
	
	The reason why $N_{(c_g,\ell_g)}$ is the domain of the largest possible extension is as follows:  Pick some $n\in N_S(P)$ and imagine there some $g'\in G$ such that $^{g'}\langle P,n\rangle\leq S$ and $(c_g,\ell_g)=(c_{g'}|_P,\ell_{g'})$.  Then $^{g'}n\in N_S({}^gP)$, so there is some $s\in N_S({}^gP)$ such that $(c_{g'n(g')^{-1}},\ell_{g'n(g')^{-1}})=(c_s,\ell_s)$.  On the other hand, the assumption that $c_g=c_{g'}$ on $P$ implies that, on restrition, we have $c_{g'n(g')^{-1}}=c_{gng^{-1}}$, and similarly we have $\ell_{g'n(g')^{-1}}=\ell_{gng^{-1}}$, which tells us that $n\in N_{(c_g,\ell_g)}$, as desired.
	\end{itemize}
	\end{remark}
	
	\begin{lemma}\label{ambient extension for actions}
	Let $P\leq S$ and $g\in G$ be such that $^gP\leq S$ is fully $X$-centralized in $G$.  Then there is a $g'\in G$ such that $^{g'}N_{(c_g,\ell_g)}\leq S$ and $(c_{g'}|_P,\ell_{g'})=(c_g,\ell_g)$.  In other words, $g'$ defines an extension of $(c_g,\ell_g)$ to a morphism in $\fX_G\left(N_{(c_g,\ell_g)}, S\right)$.
	\begin{proof}	
	If $^gP$ is fully $X$-centralized, then
	\begin{eqnarray*}
	[N_S({}^gP):N_S({}^gP)\cdot Z_G({}^gP;X)]&=&\frac{\left|N_S({}^gP)\right|\left|Z_G({}^gP;X)\right|/|Z_S({}^gP;X)|}{|N_S({}^gP)|}\\
	&=&[Z_S({}^gP;X):Z_G({}^gP;X)]
	\end{eqnarray*}
	is prime to $ p$, so $N_S({}^gP)\in\Syl_ p\left(N_S({}^gP)\cdot Z_G({}^gP;X)\right)$.  From the definition of the extender we have $^gN_{(c_g,\ell_g)}$ is a $ p$-subgroup of $N_S({}^gP)\cdot Z_G({}^gP;X)$, so we can choose some $z\in Z_G({}^gP;X)$ so that $^{zg}N_{(c_g,\ell_g)}\leq N_S({}^gP)$.  Then setting $g'=zg$ gives the desired extension $(c_{g'},\ell_{g'})$.
	\end{proof}
	\end{lemma}

\subsection{Abstract fusion action systems}\label{abstract fusion action system subsection}
	
	In this section we describe an abstraction of the fusion action systems of Section \ref{ambient fusion action systems} that makes no reference to an ambient group.  Let us start by naming our universe of discourse:
	
	\begin{definition}
	Let $\fU:=\fU(S;X)$ be the category whose objects are all subgroups $P\leq S$ and whose morphisms are given by 
	\[
	\fU(P,Q)=\left\{(\varphi,\sigma)\big|\varphi\in\Inj(P,Q),\sigma\in\Sigma_X,\textrm{ and }\varphi\textrm{ intertwines }\sigma\right\}.
	\]
	Recall that the condition on the pair $(\varphi,\sigma)$ means that for all $p\in P$ we have $\sigma\ell_p\sigma^{-1}=\ell_{\varphi(p)}$.
	\end{definition}

	By identifying $\fU$ as the ``universe of discourse,'' we mean that this is the category in which all of our abstract fusion action systems will live.
	
	\begin{definition}
		An \emph{abstract fusion action system} of $S$ acting on $X$ is a category $\fX$ such that
		$
		\fX_S\subseteq\fX\subseteq\fU
		$
		and that satisfies the \emph{Divisibility Axiom}:  Every morphism of $\fX$ factors as an isomorphism followed by an inclusion.
		
		The \emph{underlying fusion system} of  $\fX$ is the fusion system $\cF^\fX$ on $S$, where
		\[
		\cF^\fX(P,Q)=\left\{\varphi\in\Inj(P,Q)\big|\exists\sigma\in\Sigma_X\textrm{ s.t. }(\varphi,\sigma)\in\fX(P,Q)\right\}.
		\]
	Thus $\cF^\fX$ is the fusion system obtained by projecting the morphisms of $\fX$ onto the first coordinate.  
	\end{definition}
	
	By referring to this underlying fusion system we are able to speak of $P\leq S$ being fully normalized or fully centralized (by definition, in $\cF^\fX$), as well as fully $X$-normalized and $X$-centralized.  The latter two are again formed by considering the normalizers or centralizers of subgroups intersected with $C$, the core of the action of $S$ on $X$.
	
	\begin{lemma}\label{core strongly closed}
	The core $C$ is strongly closed in $\cF^\fX$, i.e., if $\varphi$ is a morphism of $\cF^\fX$ and $c\in C$ is in the source of $\varphi$, then $\varphi(c)\in C$.
	\begin{proof}
	If $c\in C$ and $\varphi\in\cF^\fX(\langle c\rangle, S)$, then there is some $(\varphi,\sigma)\in\fX(\langle c\rangle,S)$.  Then  $\varphi$ and $\sigma$'s being intertwined implies that $\ell_{\varphi(c)}=\sigma\ell_c\sigma^{-1}=\id_X$ and thus $\varphi(c)\in c$.
	\end{proof}
	\end{lemma}
	
	Of course, we need further axioms to restrict ourselves to the interesting cases.  The goal is to turn the observations of Lemmas \ref{ambient Sylow conditions} and \ref{ambient extension for actions} into the \emph{Saturation Axioms} for abstract fusion action systems.  
	
	Note that all the terms of Definition \ref{diamond group quotients} have obvious analogues in this context, so we may talk about the group the projection maps $\pi_\cF:\fX(P,Q)\to\cF^\fX(P,Q)$ and $\pi_\Sigma:\fX(P,Q)\to\Sigma_X$.  In particular, we have $\Sigma_\fX(P)$ is the group of permutations of $X$ that appear as a second coordinate of some automorphism of $P$ in $\fX$, while $\cF^\fX(P)_0$ and $\Sigma_\fX(P)_0$ are the kernels of the natural projections $\fX(P)\to\Sigma_\fX(P)$ and $\fX(P)\to\cF^\fX(P)$, respectively.  All of these groups will be referred to as the \emph{$\fX$-automizers} of the subgroup $P$.
	
	We thus find ourselves in the situation where'd like to compare the various $\fX$-automizers to their $S$-automizer subgroups, and require that these inclusions be Sylow in certain circumstances.  Unfortunately, without reference to an ambient group $G$, we do not have the entirety of Figure \ref{S vs. G}, but only the truncated Figure \ref{S vs. fX}.
	
	\begin{figure}[h]
		\[
	\xymatrix@=18pt{
	&&&&&\ar@{}[ddddd]|!{[dddddll];[ddrrr]}\hole\\
	\\
	&&&&&&&&\ar@{}[llluu]|*+[F.]{\Sigma_\fX(P)}="1a"\\
	\ar@{}[rrrrruuu]|(.45)*+[F.]{\cF^\fX(P)}="2a"\\
	\\
	&&&\ar@{}[llluu]|*+[F.]{\Sigma_\fX(P)_0}="3a"\ar@{}[rrrrruuu]|(.55)*+[F.]{\cF^\fX(P)_0}="4a"\ar@{}[rruuuuu]|*+[F.]{\fX(P)}="5a"&&N_S(P)\ar@{}[uuuuu]|!{[ll];[rrruuu]}\hole\\
	\\
	&&&&&&&&N_S(P;X)\ar@{}[uuuuu]\ar@{-}[llluu]|*+[F.]{\Sigma_\fX^S( P)}="1b"\\
	Z_S(P)\ar@{}[uuuuu]\ar@{-}[rrrrruuu]|(.45)*+[F.]{\cF_S(P)}="2b"\\
	\\
	&&&Z_S(P;X)\ar@{}[uuuuu]\ar@{-}[llluu]|*+[F.]{\Sigma_\fX^S(P)_0}="3b"\ar@{-}[rrrrruuu]|(.55)*+[F.]{\cF_S(P)_0}="4b"\ar@{-}[rruuuuu]|*+[F.]{\fX_S(P)}="5b"
	\ar @{..}"1b";"1a"
	\ar @{..}"2b";"2a"
	\ar @{..}"3b";"3a"
	\ar @{..}"4b";"4a"
	\ar @{..}"5b";"5a"
	}
	\]
	\caption{Comparing $S$-automorphisms to $\fX$-automorphisms}\label{S vs. fX}
	\end{figure}
	
	What we need is an analogue of $P$'s being fully normalized, etc., that makes no reference to an ambient group.
	
	\begin{definition} Fix a fusion action system $\fX$ and a subgroup $P\leq Q$
	\begin{itemize}
	\item $P$ is \emph{fully normalized in $\fX$} if $\left|N_S(P)\right|\geq\left| N_S(Q)\right|$ for all $Q\cong_{\cF^\fX} P$.
	\item $P$ is \emph{fully centralized in $\fX$} if $\left|Z_S(P)\right|\geq\left| Z_S(Q)\right|$ for all $Q\cong_{\cF^\fX} P$.
	\item $P$ is \emph{fully $X$-normalized in $\fX$} if $\left|N_S(P;X)\right|\geq\left| N_S(Q:X)\right|$ for all $Q\cong_{\cF^\fX} P$.
	\item $P$ is \emph{fully $X$-centralized in $\fX$} if $\left|Z_S(P;X)\right|\geq\left| Z_S(Q:X)\right|$ for all $Q\cong_{\cF_\fX} P$.
	\end{itemize}
	\end{definition}
	
	If $\fX=\fX_G$ for some finite group $G$, it is well known (cf. \cite{BLO2}, Proposition 1.3) that $P$ is fully normalized in $\cF^\fX=\cF_G$ if and only if $P$ is fully normalized with respect to $G$, and similarly for both definitions of full centralization.  The following Lemma shows that the same is true for the two new terms we have introduced
	
	\begin{lemma}
		If the saturated fusion action system $\fX$ on $S$ is realized by $G$, $P\leq S$ is fully $X$-normalized if and only if
		$
			N_S(P;X)\in\Syl_p(N_G(P;X)).
		$
		
		Similarly, $P$ is fully $X$-centralized if and only if
		$
			Z_S(P;X)\in\Syl_p(Z_G(P;X)).
		$
		\begin{proof}
			We prove the result for $X$-normalizers and note that the same argument works for $X$-centralizers.  
			
			If $N_S(P;X)=S\cap N_G(P;X)$ is Sylow in $N_G(P;X)$, then $|N_S(P;X)|\geq |N_S({}^gP;X)|$ for all $g\in G$ such that $^gP\leq S$.
			
			Suppose now that $|N_S(P;X)|\geq |N_S({}^gP;X)|$ for all $g\in G$ such that $^gP\leq S$.  Pick $T\in\Syl_p(N_G(P;X))$, and $g\in G$ such that $T\cdot P\leq {}^gS$.  Then we have $P^g\leq S$, and the fact that $N_{^gS}(P;X)\in\Syl_p(N_G(P;X))$ implies that $N_S(P^g;X)\in\Syl_p(N_G(P^g;X))$ (this uses the fact that the core of the action is normal in $G$).  The assumption that $P$ is fully $X$-normalized now implies that the orders of the $X$-normalizers in $S$ of $P$ and $P^g={}^{g^{-1}}P$ are equal, which in turn forces $N_S(P;X)\in\Syl_p(N_G(P;X))$, as desired.
		\end{proof}
	\end{lemma}
	
	The final ingredient to define saturation is the notion of the extender of an isomorphism.
	
	\begin{definition}
	The \emph{extender} of $(\varphi,\sigma)\in\fX(P,Q)_\iso$ is the subgroup of $N_S(P)$ given by
	\[
	N_{(\varphi,\sigma)}=\left\{n\in N_S(P)\big|(\varphi c_n\varphi^{-1},\sigma\ell_n\sigma^{-1})\in\Aut_S(Q;X)\right\}.
	\]
	\end{definition}
	
	\begin{definition}[Saturation Axioms for Fusion Action Systems]\label{fusion action system saturation}
	The abstract action $\fX$ is \emph{saturated} if
	\begin{enumerate}
	\item  For any $P\leq S$, the following implications hold:
		\[
		\xymatrix{
		P\textrm{ fully normalized}\ar@{=>}[r]\ar@{=>}[d]&	P\textrm{ fully }X\textrm{-normalized}\ar@{=>}[d]\\
		P\textrm{ fully centralized}\ar@{=>}[r]&			P\textrm{ fully }X\textrm{-centralized}
		}
		\]
	\item If $P$ is fully normalized, then
		\begin{itemize}
			\item $\Aut_S(P)\in\Syl_p\left(\cF^\fX(P)\right)$.
			\item $\Aut_S(P;X)\in\Syl_p\left(\fX(P)\right)$.
			\item $\Sigma_\fX^S(P)\in\Syl_p\left(\Sigma_\fX(P)\right)$.
		\end{itemize}
	\item If $P$ is fully $X$-normalized, then $\cF_S(P)_0\in\Syl_p\left(\cF^\fX(P)_0\right)$.
	\item If $P$ is fully centralized, then $\Sigma_\fX^S(P)_0\in\Syl_p\left(\Sigma_\fX(P)_0\right)$.
	\item If $Q$ is fully $X$-centralized and $(\varphi,\sigma)\in\fX(P,Q)_\iso$, then there is some
		$
		\left(\widetilde\varphi,\sigma\right)\in\fX\left(N_{(\varphi,\sigma)},S\right)
		$
		that extends $(\varphi,\sigma)$, i.e., such that $\varphi=\widetilde\varphi\big|_P$.
	\end{enumerate}
	We shall refer to Points 2-4 collectively as the \emph{Sylow Axioms} and Point 5 as the \emph{Extension Axiom} for fusion action systems.
	\end{definition}

	\begin{notation}
		Given a fusion action system $\fX$, $P\leq S$, and $Q\leq N_S(P)$, we shall denote by $\overline Q$ the image of $Q$ in $\fX(P)$; that is, $\overline Q=\left\{(c_q|_P,\ell_q)\big|q\in Q\right\}$.
	\end{notation}	

	It is obvious that the Saturation Axioms for fusion action systems are very similar to those of fusion systems.  Let us note a further connection between the two.
	
	\begin{prop}\label{underlying fusion system is saturated}
	If the abstract fusion action system $\fX$ is saturated, the underlying fusion system $\cF^\fX$ is as well.
	\begin{proof}
		All of the saturation conditions are clear except for the Extension Axiom.  Pick $\varphi\in\cF^\fX(P,Q)_\iso$ with $Q$ fully centralized in $\cF^\fX$.  Recall that
		\[
			N_\varphi=\left\{n\in N_S(P)\big|\varphi c_n\varphi^{-1}\in\Aut_S(Q)\right\}.
		\]
		If $(\varphi,\sigma)\in\fX(P,Q)$ is a morphism that lies over $\varphi$, we have
		\[
			N_{(\varphi,\sigma)}=\left\{n\in N_S(P)\big|\left(\varphi c_n\varphi{-1},\sigma\ell_n\sigma^{-1}\right)\in\Aut_S(Q;X)\right\}.
		\]
		Clearly $N_{(\varphi,\sigma)}\leq N_\varphi$ for all possible choices of $\sigma$.  If we can show that we have equality for some particular choice of $\sigma$, the Extension Axiom for the fusion action system $\fX$ will imply the desired Axiom for the fusion system $\cF^\fX$.
		
		We can restate the problem as follows:  Fix some $\sigma$ such that $(\varphi,\sigma)\in\fX(P,Q)$, and let $N\leq\fX(Q)$ denote the group
		$
			(\varphi,\sigma)\circ \overline{N_\varphi}\circ(\varphi,\sigma)^{-1}
		$.  Then $N$ is a $ p$-subgroup of $\fX(Q)$.  If $N\leq\Aut_S(Q;X)$, by definition $N_\varphi=N_{(\varphi,\sigma)}$ and we are done.  The goal then becomes to show that there is some $\tau\in\Sigma_\fX(Q)_0$ such that $^\tau N\leq\Aut_S(Q;X)$, for this will imply that $N_\varphi=N_{(\varphi,\tau\sigma)}$ and complete the proof.  (Here and elsewhere we shall ignore the distinction between an element $\sigma\in\Sigma_S$ and the corresponding morphism $(\id_Q,\sigma)$ that lies in $\fX$.)
		
		For $n\in N_\varphi$, we have $\varphi c_n\varphi^{-1}\in\Aut_S(Q)$, so $(\varphi c_n\varphi^{-1},\sigma\ell_n \sigma^{-1})$ differs from an element of $\Aut_S(Q;X)$ by an element of $\Sigma_\fX(Q)_0$.  The claim is that $\Aut_S(Q;X)$ is Sylow in $\Sigma_\fX(Q)_0\cdot\Aut_S(Q;x)$.  We compute
		\begin{eqnarray*}
			[\Aut_S(Q;X):\Sigma_\fX(Q)_0\cdot\Aut_S(Q;X)]
			&=&\frac{\left|\Sigma_\fX(Q)_0\right|\left|\Aut_S(Q;X)\right|}{\left|\Sigma_\fX(Q)_0\cap\Aut_S(Q;X)\right|\left|\Aut_S(Q;X)\right|}\\
			&=&[\Sigma_\fX(Q)_0\cap\Aut_S(Q):\Sigma_\fX(Q)_0]
		\end{eqnarray*}
		By Axiom (4) for saturated fusion action systems, $\Sigma_\fX(Q)_0\cap\Aut_S(Q;X)=\Sigma_\fX^S(Q)_0$ is Sylow in $\Sigma_\fX(Q)_0$.  Thus there is some $\tau\in\Sigma_\fX(Q)_0$ such that $^\tau N\leq\Aut_S(Q;X)$, and the result is proved.
	\end{proof}
	\end{prop}
	
	It is possible to define a morphism of fusion action systems.  Basically, the obvious thing to try works, but as we shall make no use of the category of fusion action systems in this article, we shall note define morphisms here.
	
	We conclude this section with some general properties of abstract fusion action systems, starting with a sort of converse to the Sylow Axioms of the Saturation Axioms for fusion action systems.
		
	\begin{prop}
	Fix a saturated fusion action system $\fX$ and a subgroup $P\leq S$.
	\begin{itemize}
	\item[(a)]  The following are equivalent:
		\begin{enumerate}
		\item $P$ is fully normalized in $\fX$.
		\item $P$ is fully centralized in $\fX$ and $\Aut_S(P)\in\Syl_ p(\cF(P))$.
		\item $P$ is fully $X$-normalized in $\fX$ and $\Sigma_\fX^S(P)\in\Syl_ p\left(\Sigma_\fX(P)\right)$.
		\item $P$ is fully $X$-centralized in $\fX$ and $\Aut_S(P;X)\in\Syl_ p\left(\fX(P)\right)$.
		\end{enumerate}
	\item[(b)]  $P$ is fully centralized iff $P$ is fully $X$-centralized and $\Sigma_\fX^S(P)_0\in\Syl_ p\left(\Sigma_\fX(P)_0\right)$.
	\item[(c)] $P$ is fully $X$-normalized iff $P$ is fully $X$-centralized and $\cF_S(P)_0\in\Syl_ p\left(\cF(P)_0\right)$.
	\end{itemize}
	\begin{proof}
	Half of these implications are part of the Saturation Axioms, so prove only the remaining ones.
	
	\begin{itemize}
		\item[(a)] $P\leq S$ is fully normalized if:
		\begin{itemize}
		\item $P$ is fully centralized and $\Aut_S(P)\in\Syl_ p(\cF(P))$.  In this case, the assumptions on $P$ together with the inclusion of short exact sequences
		\[
		\xymatrix
		{
		1\ar[r]&	\Sigma_\fX(P)_0\ar[r]&	\fX(P)\ar[r]&	\cF(P)\ar[r]	&	1\\
		1\ar[r]&	\Sigma_\fX^S(P)_0\ar@{^(->}[u]^{\mathrm{Sylow}}\ar[r]&	\Aut_S(P;X)\ar@{^(->}[u]\ar[r]&
			\Aut_S(P)\ar@{^(->}[u]_{\mathrm{Sylow}}\ar[r]&	1
		}
		\]
		force $\Aut_S(P;X)\in\Syl_ p\left(\fX(P)\right)$.  Then if $Q$ is fully normalized and $\cF$-conjugate to $P$ via $(\varphi,\sigma)\in\fX(Q,P)_\iso$, we have that
		$
			(\varphi,\sigma)\circ\Aut_S(Q;X)\circ(\varphi,\sigma)^{-1}\leq\fX(P)
		$
		is an inclusion of a $ p$-subgroup.  Thus there is $(\psi,\tau)\in\fX(P)$ such that
		$
			(\psi\varphi,\tau\sigma)\circ\Aut_S(Q;X)\circ(\psi\varphi,\tau\sigma)^{-1}\leq\Aut_S(P;X)
		$.
		
		Note that $\Aut_S(Q;X)$ is the image in $\fX(Q)$ of $N_S(Q)$, and that $P$'s being fully centralized implies that it is fully $X$-centralized.  The Extension Axiom now gives the existence of $(\widetilde{\psi\varphi},\tau\sigma)\in\fX(N_S(Q),N_S(P))$, from which it follows that $\big|N_S(Q)\big|=\big|N_S(P)\big|$ and $P$ is fully normalized as well.
		\item
		$P$ is fully $X$-normalized and $\Sigma_\fX^S(P)\in\Syl_ p\left(\Sigma_\fX(P)\right)$.  The inclusion of short exact sequences
		\[
		\xymatrix
		{
		1\ar[r]&	\cF(P)_0\ar[r]&	\fX(P)\ar[r]&	\Sigma_\fX(P)\ar[r]&1\\
		1\ar[r]&	\cF_S(P)_0\ar@{^(->}[u]^{\mathrm{Sylow}}\ar[r]&	\Aut_S(P;X)\ar@{^(->}[u]\ar[r]&
			\Sigma_\fX^S(P)\ar@{^(->}[u]_{\mathrm{Sylow}}\ar[r]&	1
		}
		\]
		shows that $\Aut_S(P;X)\in\Syl_ p\left(\fX(P)\right)$.  The rest of the proof is the same as the first point.
		\item  $P$ is fully $X$-centralized and $\Aut_S(P;X)\in\Syl_ p\left(\fX(P)\right)$.  The argument is the same as the end of the previous two.
		\end{itemize}		
		
		\item[(b)] $P$ is fully centralized if $P$ is fully $X$-centralized and  $\Sigma_\fX^S(P)_0\in\Syl_ p\left(\Sigma_\fX(P)_0\right)$:  Observe that $\big|Z_S(P)\big|=\big|Z_S(P;X)\big|\big|\Sigma_\fX^S(P)_0\big|$.  If $P\cong_{\cF^\fX} Q$, then we have $\Sigma_\fX(P)_0\cong\Sigma_\fX(Q)_0$, so it is easy to see that the order of (the $ p$-group) $Z_S(P)$ is maximized precisely when $\Sigma_\fX^S(P)_0\in\Syl_ p\left(\Sigma_\fX(P)_0\right)$ and the order of $Z_S(P;X)$ is maximized.  The result follows.
		
		\item[(c)] $P$ is fully $X$-normalized if $P$ is fully $X$-centralized and $\cF_S(P)_0\in\Syl_ p\left(\cF(P)_0\right)$:  The same argument as the previous paragraph applies, with the observation that $\big| N_S(P)\big|=\big|N_S(P;X)\big|\big|\cF_S(P)_0\big|$.
		\end{itemize}
	\end{proof}
	\end{prop}
			
	Let us close with a basic reality-check for our Axioms:  If these Saturation Axioms are ``right,'' the least we would expect is for some sort of Alperin Fusion Theorem to hold.
	
	\begin{prop}\label{fully X-normalized target}
	If $(\varphi,\sigma)\in\fX(P,Q)$ is an isomorphism such that $Q$ is fully normalized, then there are morphisms $(\widetilde\varphi,\sigma')\in\fX(N_S(P),S)$ and $(\psi,\tau)\in\fX(P)$ such that $\widetilde\varphi\big|_P=\varphi\circ\psi$ and $\sigma'=\sigma\circ\tau$.
	\begin{proof}
		By the assumption that $Q$ is fully normalized, the Saturation Axioms state that $\Aut_S(Q;X)\in\Syl_ p(\fX(Q))$, so the $ p$-group $(\varphi,\sigma)\circ\Aut_S(P;X)\circ(\varphi^{-1},\sigma^{-1})\leq\fX(Q)$ is subconjugate to $\Aut_S(Q;X)$.  Say
		\[
			(\chi,\upsilon)\circ(\varphi,\sigma)\circ\Aut_S(P;X)
			\circ(\varphi,\sigma)^{-1}\circ(\chi,\upsilon)^{-1}\leq\Aut_S(Q;X)
		\]
		and set $(\psi,\tau)=(\varphi^{-1}\chi\varphi,\sigma^{-1}\upsilon\sigma)\in\fX(P)$, so that
		\[
			(\varphi,\sigma)\circ(\psi,\tau)\circ\Aut_S(P;X)\circ
			(\psi,\tau)^{-1}\circ(\varphi,\sigma)^{-1}\leq \Aut_S(Q;X).
		\]
		Now the fact that $Q$ is fully normalized and thus fully $X$-centralized implies that $(\varphi\psi,\sigma\tau)\in\fX(P,Q)$ has an extension $(\widetilde\varphi,\sigma')\in\fX\left(N_{(\varphi\chi,\sigma\tau)}(P),S\right)$ by the Extension Axiom.  Recall that
		\[
			N_{(\varphi\chi,\sigma\tau)}(P)=
			\left\{s\in N_S(P)\big|(\varphi\chi\circ c_s\circ 
			\chi^{-1}\varphi^{-1},\sigma\tau\circ\ell_s\circ\tau^{-1}\sigma^{-1}\in\Aut_S(Q;X)\right\}
		\]
		so that by construction $N_S(P)\leq N_{(\varphi\chi,\sigma\tau)}(P)$.  The result follows.
		\end{proof}
	\end{prop}
	
	One could also prove analogous statements in the cases that the target is fully centralized or $X$-normalized, but we will not have need for such results.
	
	\begin{theorem}[Alperin Fusion Theorem for Fusion Action Systems]\label{Alperin for fusion action systems}
		If the fusion action system $\fX$ is saturated, every morphism of $\fX$ can be written as a composite of automorphisms of subgroups of $S$ and inclusions.
		\begin{proof}
			The proof goes by downward induction on the order of the source $P$.  If $P=S$, then $(\varphi,\sigma)\in\fX(S)$, and there's nothing to prove.
			
			Therefore suppose that $P\lneq S$.  Without loss of generality, we may assume that $Q=\varphi(P)$ is fully normalized:  Otherwise, pick $Q'$ in the $\cF$-conjugacy class of $P$ that is fully normalized and an isomorphism $(\psi,\tau)\in\fX(P,Q')$.  If the result is true for $(\psi,\tau)$ and $(\varphi,\sigma)\circ(\psi,\tau)^{-1}$, it clearly is for $(\varphi,\sigma)$ as well.
			
			Proposition \ref{fully X-normalized target} shows that if the target of $(\varphi,\sigma)$ is fully normalized, then $(\varphi,\sigma)$ can be composed with an element of $\fX(P)$ so that the resulting morphism extends to $\fX(N_S(P),S)$.  But $N_S(P)\gneq P$ as $P$ is a proper subgroup, and the inductive hypothesis gives the rest.
		\end{proof}
	\end{theorem}

	\subsection{Faithful fusion action systems}\label{faithful fusion action systems}
	
	Note that a fusion system $\cF$ is actually a special case of the notion of fusion action system, where the set $X$ being acted on is a one-point set.  In particular, fusion systems could be thought of as fusion action systems where the $S$-action is trivial.
	
	In this section we examine in closer detail a fusion action system $\fX$ at the opposite extreme to $\cF$:  We assume that $X$ is faithful as an $S$-set.   We will see that in this situation, the underlying saturated fusion system $\cF$ is always realizable by a finite group $G$, and moreover that $G$ acts on $X$ in such a way that $\fX=\fX_G$.

	Let $G=\fX(1)$ be the group of permutations of $X$ that appear in some morphism of $\fX$.  The map $S\to G:s\mapsto\ell_s$ is an injection if we assume that $X$ is faithful as an $S$-set, so we will identify $S$ with its image in $G$.   The Saturation Axioms for fusion actions imply that $S\in\Syl_ p (G)$, so the fusion system $\cF_G$ is saturated.
	
	The following proposition shows that exotic fusion systems cannot be the underlying fusion system of a $S$-faithful fusion action system.

	\begin{prop}\label{faithful actions realize fusion systems}
		If $\fX$ is a saturated fusion action system and $S$ acts faithfully on $X$, then $\cF^\fX=\cF_G$.
		\begin{proof}
			First we show that $\cF^\fX\subseteq\cF_G$.  It suffices to show that any $\varphi\in\cF^\fX(P,S)$ is realized by conjugation by an element of $G$.  Pick some $(\varphi,\sigma)\in\fX(P,S)$, so that $\varphi$ and $\sigma$ are intertwined, so $\sigma\ell_p\sigma^{-1}=\ell_{\varphi(p)}\in\Sigma_X$ for all $p\in P$.  This is exactly to say that $\sigma$ conjugates $P$ via $\varphi$, proving the first inclusion.
			
			For the reverse inclusion, suppose that for some $\sigma\in G$ and $P\leq S$ we have $^\sigma P\leq S$; we must show that $c_\sigma\in\cF^\fX(P,S)$.  By the Extension Axiom, the morphism $(\mathrm{id_1},\sigma)\in\fX(1)$ extends to some
			$
				(\varphi,\sigma)\in\fX(N_{(\mathrm{id_1},\sigma)},S)$, where 
				$
				N_{(\mathrm{id_1},\sigma)}=\left\{n\in S\big|\sigma\circ\ell_n\circ\sigma^{-1}\in\Aut_S(X)\right\}			$.
			In other words, for all $n\in N_{(\mathrm{id_1},\sigma)}$, there is some (necessarily unique by faithfulness of the $X$-action) $s\in S$ so that $\sigma\circ\ell_n\circ\sigma^{-1}=\ell_s$.  Therefore $N_{(\mathrm{id_1},\sigma)}$ is the maximal subgroup of $S$ that is conjugated into $S$ by $\sigma$.  Thus $P\leq N_{(\mathrm{id_1},\sigma)}$, and if we can show that the assignment $n\mapsto s$ is equal to $\varphi$, the result will follow, as $\varphi$ is by assumption a morphism of $\cF^\fX$.
			
			This final assertion follows again from the fact that $\varphi$ and $\sigma$ are intertwined: $\sigma\ell_n\sigma^{-1}=\ell_{\varphi(n)}$ for $n\in N_{(\mathrm{id_1},\sigma)}$ forces $s=\varphi(n)$ in the above notation.  The result is proved.
		\end{proof}
	\end{prop}
	
	In this situation, there is a natural action of $G=\fX(1)\leq\Sigma_X$ on $X$, and the following is immediate:
	
	\begin{cor}\label{faithful actions realize fusion action systems}
		The fusion action system $\fX$ is realized by the $G$-action on $X$; i.e., $\fX=\fX_G$.
		\begin{proof}
			The key point is that if the $S$-action on $X$ is faithful, then any morphism $(\varphi,\sigma)\in\fX(P,Q)$ is actually determined by $\sigma$ alone.  This follows from the assumption that $\varphi$ and $\sigma$ are intertwined, so that for all $p\in P$, we have $\ell_{\varphi(p)}=\sigma\ell_p\sigma^{-1}$ along with the identification $p\leftrightarrow \ell_p$.  In particular, $\fX\subseteq\fX_G$, since $(\varphi,\sigma)\in\fX(P,Q)$ implies that $^\sigma P\leq Q$ and $c_\sigma=\varphi$.
			
			On the other hand, for any $\sigma\in N_G(P,Q)$ we have $(c_\sigma,\sigma)\in\fX_G(P,Q)$, and the claim is that this is also a morphism of $\fX(P,Q)$.  This follows from the Extension Axiom of saturated fusion actions as in Proposition \ref{faithful actions realize fusion systems} and the observation that $^\sigma P\leq S$, implies $P\leq N_{(\id_1,\sigma)}$.  Thus $\fX_G\subseteq\fX$.  
		\end{proof}
	\end{cor}
	
	\subsection{Core subsystems}\label{core subsystems}

	In this section let $\fX$ be a saturated fusion action system.  Recall that $C\leq S$ denotes the core of the $S$ action on $X$.
	
	\begin{definition}
		The \emph{core fusion system} associated to $\fX$ is the fusion system $\cC$ on $C$ with morphisms $\cC(P,Q)=\left\{\varphi\in\Inj(P,Q)\big|(\varphi,\id_X)\in\fX(P,Q)\right\}$.		
	\end{definition}
	
	\begin{remark}
	There is a natural functor $\fX\to\cB\Sigma_X$ that sends inclusions to the identity; this is precisely what the various maps $\pi_\Sigma$ piece together to give.  The fusion system $\cC$ can be thought of the kernel fusion action system of this functor, once we note that, as the action of $\cC$ on $X$ is trivial, there is no distinction between a fusion system and a fusion action system.
	\end{remark}
	
	The first goal of this section is that $\cC$ is saturated as a fusion system.
	
	\begin{prop}\label{core almost normal}
	$\cC$ satisfies the following properties:
	\begin{enumerate}
		\item $C$ is strongly closed in $\cF^\fX$.
		\item For all $P,Q\leq C$, $\varphi\in\cC(P,Q)$ and $\psi\in\cF^\fX(C)$, we have $\psi\varphi\psi^{-1}\in\cC(\psi P,\psi Q)$.
		\item For all $P,Q\leq C$ and $\chi\in\cF^\fX(P,Q)$, there are $\psi\in\cF^\fX(C)$ and $\varphi\in\cC(\psi (P),Q)$ such that $\chi=\varphi\circ\psi\big|_P^{\psi(P)}$.
	\end{enumerate}
	\begin{proof}$\left.\right.$
	\begin{enumerate}
		\item Cf. Lemma \ref{core strongly closed}.
		\item Pick some $(\psi,\tau)\in\fX(C)$ that lies over $\psi\in\cF^\fX(C)$.  Then
		$
			(\psi,\tau)(\varphi,\id_X)(\psi,\tau)^{-1}=(\psi\varphi\psi^{-1},\id_X)
		$,
		and the result follows.
		\item  Let $(\chi,\upsilon)\in\fX(P,Q)$ lie over $\chi\in\cF^\fX(P,Q)$, and consider $(\id_1,\upsilon)\in\fX(1)$.  By the Extension Axiom for saturated fusion action systems, $(\id_1,\upsilon)$ extends to a map $(\widetilde\psi,\upsilon)\in\fX(N_{(\id_1,\upsilon)},S)$, and it follows easily from the definition that $C\leq N_{(\id_1,\upsilon)}$.  Let $(\psi,\upsilon)\in\fX(C)$ be the restriction to $C$ (which is necessarily an automorphism of $C$), and set 
		$
			(\varphi,\id_X)=(\chi,\upsilon)\circ(\psi,\upsilon)^{-1}\big|_{\psi(P)}^P\in\fX(\psi(P),Q)
		$.
		Projection to $\cF^\fX$ shows now that $\psi$ and $\varphi$ have the desired property.
	\end{enumerate}
	\end{proof}
	\end{prop}
	
	The content of Proposition \ref{core almost normal} is that $\cC$ is almost a normal subsystem of $\cF$ in the terminology of Puig, which we shall call an \emph{invariant subsystem}, following Aschbacher.  We say ``almost'' because in the definition of invariance, we would also need that $\cC$ is saturated as a fusion system, which is our current goal.  We need the third point of Proposition \ref{core almost normal} to prove saturation, but luckily the logic is not circular.
	
	\begin{cor}\label{core whatever vs. fusion X-whatever}
	Every $P\leq C$ is fully normalized (resp. centralized) in $\cC$ if and only if $P$ is fully $X$-normalized (resp. $X$-centralized) in $\fX$.
	\begin{proof}
		Note that $N_S(P;X)=N_C(P)$ and $Z_S(P;X)=Z_C(P)$ by definition, which makes the ``if'' implications obvious.  We prove the ``only if'' implication for normalizers; the same argument works for centralizers.
		
		Suppose that $P$ is fully normalized in $\cC$, and let $Q\leq C$ be fully $X$-normalized and $\cF^\fX$-conjugate to $P$ (such a $Q$ must exist as $C$ is strongly closed in $\cF^\fX$).  For any isomorphism $\chi\in\cF^\fX(P,Q)_\iso$, we can find $\psi\in\cF^\fX(C)$ and $\varphi\in\cC(\psi(P),Q)$ such that $\chi=\varphi\circ\psi\big|_P^{\psi(P)}$ by the third point of Proposition \ref{core almost normal}.  As $\psi$ extends to an automorphism of $C$, it sends the $X$-normalizer of $P$ to the $X$-normalizer of $\psi(P)$.  As $\varphi$ is necessarily an isomorphism in $\cC$, the assumption that $|N_C(P)|\geq|N_C(\psi(Q))|$ forces that $|N_S(P;X)|\geq|N_S(Q;X)|$.  The assumption that $Q$ is fully $X$-normalized forces equality, and the result is proved.
	\end{proof}
	\end{cor}
	
	\begin{theorem}\label{core subsystem is saturated}
	$\cC$ is a saturated fusion system.
	\begin{proof}
		Using the Saturation Axioms of \cite{OnofreiStancuCharacteristicSubgroup}, it suffices to show that $\Aut_C(C)\in\Syl_ p(\cC(C))$ and that fully normalized subgroups satisfy the Extension Axiom.
		
		To see that $\Aut_C(C)\in\Syl_ p(\cC(C))$:  We've already noted that $C$ is strongly closed in $\cF^\fX$, and therefore is fully $X$-normalized.  By the Saturation Axioms for fusion action systems, we have $\cF^\fX_S(C)_0\in\Syl_ p(\cF^\fX(C)_0)$.  It is easy to see from the definitions that $\cF^\fX(C)_0=\cC(C)$, and we can describe $\cF^\fX_S(C)_0$ as the subgroup of $\Aut_S(C)$ consisting of those automorphisms that are induced by an element of $S$ that acts trivially on $X$, so that $\cF^\fX_S(C)_0=\Aut_C(C)$.
		
		To check the Extension Axiom:  If $P\leq C$ is fully normalized in $\cC$, by Corollary \ref{core whatever vs. fusion X-whatever} we see that $P$ is fully $X$-normalized in $\cF^\fX$, and thus fully $X$-centralized.  Given an isomorphism $\varphi\in\cC(Q,P)_\iso$, we have $(\varphi,\id_X)\in\fX(Q,P)$, and the extension condition for fusion actions gives a morphism
		$
			(\widetilde\varphi,\id_X)\in\fX(N_{(\varphi,\id_X)}^{\fX},C)
		$
		where the target may be assumed to be $C$ instead of $S$ because $C$ is strongly closed in $\cF^\fX$.  We wish to show that 
		$
			N_\varphi^\cC=\left\{n\in N_C(Q)\big| \varphi c_n\varphi^{-1}\in\Aut_K(P)\right\}
		$
		is contained in $N_{(\varphi,\id_X)}^{\fX}$, as then $\widetilde\varphi|_{N_\varphi^\cC}$ will give the desired extension of $\varphi$.  We have
		\[
			N_{(\varphi,\id_X)}^{\fX}=\left\{n\in N_S(Q)\big|(\varphi\circ c_n\circ\varphi^{-1},\ell_n)\in\Aut_S(P;X)\right\}.
		\]
		Since $C$ acts trivially on $S$, for $n\in N_C(Q)$ such that $\varphi\circ c_n\circ\varphi^{-1}=c_k\in\Aut_C(P)$ we have 
		\[
		(\varphi\circ c_n\circ\varphi^{-1},\ell_n)=(c_k,\id_X)=(c_k,\ell_k)\in\Aut_C(P;X)\leq\Aut_S(P;X)
		\]
		and the result is proved.		
	\end{proof}
	\end{theorem}

	\begin{cor}
		$\cC$ is a normal subsystem of $\cF^\fX$ in the sense of \cite{AschbacherNormalSubsystem}.
		\begin{proof}
		Recall that in Aschbacher's definition of normality, in addition to the results of Proposition \ref{core almost normal} and Theorem \ref{core subsystem is saturated}, we must also verify the technical condition that for all $\varphi\in\cC(C)$ there is an extension $\widetilde\varphi\in\cF^\fX(C\cdot Z_S(C))$ such that $[\widetilde\varphi,Z_S(C)]\leq Z(C)$, or equivalently $[\widetilde\varphi,Z_S(C)]\leq C$.
		
		If $\varphi\in\cC(C)$ we have by definition $(\varphi,\id_X)\in\fX(C)$.  Since $C$ is strongly closed in $\cF^\fX$, it is in particular fully $X$-centralized, and so $(\varphi,\id_X)$ extends to some morphism in $\fX(N_{(\varphi,\id_X)},S)$.  We claim that $Z_S(C)\leq N_{(\varphi,\id_X)}$: Indeed, for $z\in Z_S(C)$, we have
		$
		(\varphi,\id_X)\circ(c_z,\ell_z)\circ(\varphi,\id_X)^{-1}=(\id_C,\ell_z)=(c_z,\ell_z)
		$
		since $z$ acts trivially on $C$.  Let $(\widetilde\varphi,\id_X)$ be the extension of $(\varphi,\id_X)$ to $C\cdot Z_S(C)$.  As the strongly $\cF^\fX$-closed $C$ is in the domain of $\widetilde\varphi$, we in fact have $(\varphi,\id_X)\in\fX(C\cdot Z_S(C))$.
		
		Now we just need to show that $[\widetilde\varphi,Z_S(C)]\leq C$:  For all $z\in Z_S(C)$ we want $\widetilde\varphi(z)z^{-1}\leq C$, or equivalently $\ell_{\widetilde\varphi(z)}\ell_z^{-1}=\id_X$.  But since the pair $(\widetilde\varphi,\id_X)$ is intertwined, we have have $\ell_{\widetilde\varphi(z)}=\ell_z$, and the result follows.
		\end{proof}
	\end{cor}
	
	We close this section by describing how the fusion action system $\fX$ can be thought of as an extension of $\cC$ by the finite group $G=\fX(1)$.  First note that we have the functors
	$
	\xymatrix
	{
	\cC\ar@{^{(}->}[r]^-\iota&
	\fX\ar@{->>}[r]^-{\pi_\cT}&
	\cT_G
	},
	$
	where $\iota$ is ``injective'' in the sense that $\cC$  naturally sits as a subcategory of $\fX$, and $\pi_\cT$ is ``surjective'' in the sense that every morphism of $\cT_G$ is in the image of $\pi_\cT$.  Just as extensions of the finite group $N$ by $H$ determines a morphism $H\to\Out(N)$, would like to say that a fusion action system gives rise to an outer action of the ``quotient'' group $G$ on the core subsystem $\cC$.  First we must understand what is meant by an ``outer action.''  We begin by recalling the definition of morphism of fusion system.
		
			\begin{definition}\label{fusion morphism}
			  For fusion systems $(S,\cF)$ and $(S',\cF')$, a group map $\alpha:S\to S'$ is \emph{fusion-preserving} if there exists a functor $F_\alpha:\cF\to\cF'$ such that $F_\alpha(P)=\alpha P$ for all $P\leq S$ and the following diagram of groups commutes for all $\varphi\in\cF(P,Q)$:
	\[
	\xymatrix{
		P\ar[r]^-\alpha\ar[d]_\varphi	&	\alpha P\ar[d]^{F_\alpha(\varphi)}\\
		Q\ar[r]_-\alpha				&	\alpha Q
	}
	\]
	Such a functor $F_\alpha$ is necessarily unique if it exists.
	
	A \emph{morphism of fusion systems} $(S,\cF)\to(S',\cF')$ is a fusion-preserving morphism $\alpha:S\to S'$.  The set of all such morphisms is denoted $\Hom(\cF,\cF')$, and in the case that $(S,\cF)=(S',\cF')$ the resulting group of fusion preserving automorphisms of $S$ is written $\Aut(\cF)$.
	\end{definition}
	
	\begin{remark}
	Implicit in Definition \ref{fusion morphism} is the easy result that if $\alpha\in\Hom(\cF,\cF')$, $\alpha'\in\Hom(\cF',\cF'')$, and $F_\alpha, F_{\alpha'}$ are the corresponding functors, then $F_{\alpha'\alpha}=F_{\alpha'}F_\alpha$ is the functor that shows that $\alpha'\alpha\in\Hom(\cF,\cF'')$.
	\end{remark}
	
	There is an easy way to check which automorphisms of the group $S$ actually lie in $\Aut(\cF)$.  To express this, let us recall a piece of terminology from Remark \ref{obvious extension facts}:
	
	\begin{definition}\label{translation morphisms}
	Given a $p$-group $S$, subgroups $P,Q\leq R\leq S$, and an injective map $\gamma:R\to S$, the \emph{translation along $\gamma$ from $P$ to $Q$} is the map
	\[
	\begin{array}{rl}t_\gamma\big|_P^Q:&
	\xymatrixrowsep{1ex}
	\xymatrix{
	\Hom(P,Q)\ar[r]&\Hom(\gamma P,\gamma Q)\\
	\eta\ar@{|->}[r]&\gamma\eta\gamma^{-1}
	}
	\end{array}
	\]
	In cases where there will be no confusion, we simply write $t_\gamma$ for $t_\gamma\big|_P^Q$.
	\end{definition}
	
	\begin{prop}\label{fusion automorphism condition}
	Given a fusion system $\cF$ on the $p$-group $S$ and $\alpha\in\Aut(S)$, we have $\alpha\in\Aut(\cF)$ if and only if $t_\alpha\left(\cF(P,Q)\right)=\cF(\alpha P,\alpha Q)$ for all $P,Q\leq S$.
	\begin{proof}
		If there is a functor $F_\alpha$ making $(\alpha,F_\alpha)$ a morphism of $\cF$, for all $\varphi\in\cF(P,Q)$ we must have 
		$
		F_\alpha(\varphi)=\alpha\varphi\alpha^{-1}=t_\alpha(\varphi)\in\cF(\alpha P,\alpha Q)
		$.
		Thus if $\alpha$ extends to a functor of $\cF$, it follows that $t_\alpha(\cF(P,Q))\subseteq\cF(\alpha P,\alpha Q)$ for all $P,Q\leq S$.  The fact that $\alpha$ is an automorphism of $S$ forces equality.
		
		Conversely, if we have $t_\alpha(\cF(P,Q))=\cF(\alpha P,\alpha Q)$, the assignment $\varphi\mapsto t_\alpha\varphi$ can easily be seen to give the action of the desired $F_\alpha$ on morphisms.
	\end{proof}
	\end{prop}
	
	\begin{example}\label{motivation of Out(F)}
		$\cF(S)\subseteq\Aut(\cF)$.  In other words, if $\varphi\in\Aut(S)$ is a morphism in $\cF$, $\varphi$ is actually fusion preserving.  This follows immediately from Proposition \ref{fusion automorphism condition} together with the Divisibility Axiom of fusion systems, which implies that the restriction of a morphism in $\cF$ to a subgroup also lies in $\cF$.
		
		Moreover, it is easy to see that $\cF(S)\trianglelefteq\Aut(\cF)$:  If $\alpha\in\Aut(\cF)$ and $\varphi\in\cF(S)$, $\alpha\varphi\alpha^{-1}=t_\alpha(\varphi)\in\cF(S)$ again by Proposition \ref{fusion automorphism condition}.
	\end{example}
	
	\begin{definition}
		The \emph{inner automorphism group of $\cF$} is $\Inn(\cF):=\cF(S)$.  The \emph{outer automorphism group of $\cF$} is the quotient $\Out(\cF):=\Aut(\cF)/\Inn(\cF)$.
	\end{definition}

	\begin{theorem}\label{fusion actions as extensions}
		The fusion action system $\fX$ determines a unique homomorphism $\kappa:G\to\Out(\cC)$.
		\begin{proof}
			An easy application of the Extension Axiom for saturated fusion action systems implies that each $\sigma\in G$ appears in the second coordinate of some $(\varphi,\sigma)\in\fX(C)$.   As $\cC$ is normal in $\cF^\fX$, $\varphi\in\Aut(\cC)$. The indeterminacy of the assignment $\sigma\mapsto\varphi$ is measured by $\cC(C)$, so the assignment $\sigma\mapsto[\varphi]\in\Out(\cC)$ is the desired map.
		\end{proof}
	\end{theorem}

\subsection{Simplification of the Saturation Axioms}\label{simplification subsection}

	The Saturation Axioms of Definition \ref{fusion action system saturation} took their form from an analysis of the situation of a finite group acting on a finite set, just as the first set of Saturation Axioms for fusion systems in \cite{BLO2} could be derived.  The problem with the original list of Axioms is that it is both redundant and annoying to check in specific cases.  The goal of this Subsection is to give a much shorter list of Saturation Axioms, which will be very similar in form to the simplified Axioms for saturated fusion systems introduced in \cite{RobertsShpectorovSaturationAxiom}.  The form of the proof follows an unpublished manuscript by Bob Oliver.

	Let $\fX$ be a fusion action system associated to a given $S$ action on $X$.

	\begin{definition}
		A subgroup $P\leq S$ is \emph{fully automized in $\fX$} if $\Aut_S(P;X)\in\Syl_p(\fX(P))$. $P$ is \emph{receiving} if for all $(\varphi,\sigma)\in\fX(Q,P)_\iso$ there exists an extension $(\widetilde\varphi,\sigma)\in\fX(N_{(\varphi,\sigma)},S)$.
	\end{definition}
	
	\begin{theorem}\label{RS axioms}
		The fusion action system $\fX$ is saturated if and only if every $\fX$-conjugacy class of subgroups contains a fully automized, receiving subgroup.
	\end{theorem}
	
	The proof of Theorem \ref{RS axioms} will take the form of most of the lemmas in this Subsection.  If $\fX$ is saturated in the original definition then every subgroup that is fully normalized in $\fX$ is both fully automized and receiving.  We must therefore just show the other implication.
	
	\begin{lemma}
		If $P\leq S$ is fully automized then $\Sigma_\fX^S(P)_0\in\Syl_p(\Sigma_\fX(P)_0)$ and $\cF^\fX_S(P)_0\in\Syl_p(\cF^\fX(P)_0)$.
	\begin{proof}
	Consider the inclusion of exact sequences
	\[
	\xymatrix{1\ar[r]&\Sigma_\fX(P)_0\ar[r]&\fX(P)\ar[r]&\cF^\fX(P)\ar[r]&1\\
	1\ar[r]&\Sigma_\fX^S(P)_0\ar[r]\ar[u]&\Aut_S(P;X)\ar[r]\ar[u]&\cF^\fX_S(P)\ar[r]\ar[u]&1}
	\]
	The middle inclusion is Sylow if and only if the outer two are, so we have the first implication.  The second follows from considering the exact sequences of the form $\xymatrix{1\ar[r]&\cF^\fX(P)_0\ar[r]&\ar[r]\fX(P)\ar[r]&\Sigma_\fX(P)\ar[r]&1}$.
	\end{proof}
	\end{lemma}
	
	\begin{lemma}\label{RS axiom extension to normalizer}
	If $P$ is fully automized and receiving then $P$ is fully normalized.
	\begin{proof}
		Pick $Q\leq S$ fully normalized and some isomorphism $(\varphi,\sigma)\in\fX(Q,P)_\iso$.  We wish to show that, perhaps after modifying the choice of isomorphism $(\varphi,\sigma)$, we can extend this isomorphism to the whole of the normalizer of $Q$; as any such map will send $N_S(Q)$ into $N_S(P)$, the choice of $Q$ to be fully normalized will complete the proof.
		
		Let $\overline{N_S(Q)}$ denote the image of $N_S(Q)$ in $\fX(Q)$.  Then $^{(\varphi,\sigma)}\overline {N_S(Q)}\leq\fX(P)$ is a $p$-subgroup.  As $P$ is fully automized, there is some $(\psi,\tau)\in\fX(P)$ such that $^{(\psi,\tau)}\left({}^{(\varphi,\sigma)}\overline{N_S(Q)}\right)\leq\Aut_S(P;X)$, so $N_S(Q)=N_{(\psi\varphi,\tau\sigma)}$.  The fact that $P$ is receiving implies the existence of an extension of $(\psi\varphi,\tau\sigma)$ to a morphism in $\fX(N_S(Q),N_S(P))_\iso$, and the result is proved.
	\end{proof}
	\end{lemma}
	
	\begin{lemma}
	Suppose that $P$ is receiving.
	\begin{enumerate}
	\item If $\Sigma_\fX^S(P)_0\in\Syl_p(\Sigma_\fX(P)_0)$, then $P$ is fully centralized in $\fX$.
	\item If $\cF^\fX_S(P)_0\in\Syl_p(\cF^\fX(P)_0)$, then $P$ is fully $X$-normalized in $\fX$.
	\end{enumerate}
	\begin{proof}
	The same proof works in both cases, so we shall just present the first.  Pick $Q\leq S$ fully centralized and $(\varphi,\sigma)\in\fX(Q,P)_\iso$.  Then $^{(\varphi,\sigma)}\overline{Z_S(Q)}$ is a $p$-subgroup of $\Sigma_\fX(P)_0$.  The condition on $P$ forces there is some $(\id_P,\tau)\in\Sigma_\fX(P)_0$ such that $^{(\id_P,\tau)}\left({}^{(\varphi,\sigma)}\overline{Z_S(Q)}\right)\leq\Sigma_\fX^S(P)_0\leq\Aut_S(P;X)$.  Thus $Z_S(Q)\leq N_{(\varphi,\tau\sigma)}$, so the fact that $P$ is receiving implies that there is a morphism of $\fX$ that sends $Z_S(Q)$ into $Z_S(P)$.  The choice of $Q$ to be fully centralized implies that $P$ is as well.
	\end{proof}
	\end{lemma}
	
	\begin{lemma}
	If $P$ is receiving then $P$ is fully $X$-centralized in $\fX$.
	\begin{proof}
		For any isomorphism $(\varphi,\sigma)\in\fX(Q,P)_\iso$, it follows from the definition that $Z_S(Q;X)\leq N_{(\varphi,\sigma)}$, so any such extension of $(\varphi,\sigma)$ will sends $Z_S(Q;X)$ into $Z_S(P;X)$.  In particular, if $Q$ is fully $X$-centralized we get the result.
	\end{proof}
	\end{lemma}
	
	\begin{prop}
	Suppose that $Q\leq S$ is fully automized and receiving and we are given $(\varphi,\sigma)\in\fX(P,Q)_\iso$.
	\begin{enumerate}
	\item If $P$ is fully $X$-centralized then $P$ is receiving.
	\item If $P$ is fully $X$-normalized then $P$ is receiving and $\cF^\fX_S(P)_0\in\Syl_p(\cF^\fX(P)_0)$.
	\item If $P$ is fully centralized then $P$ is receiving and $\Sigma_\fX^S(P)_0\in\Syl_p(\Sigma_\fX(P)_0)$.
	\item If $P$ is fully normalized then $P$ is fully automized and receiving.
	\end{enumerate}
	\begin{proof}
	By the proof of Lemma \ref{RS axiom extension to normalizer}, we may choose $(\varphi,\sigma)$ such that there is an extension $(\widetilde\varphi,\sigma)\in\fX(N_S(P),N_S(Q))$.
	\begin{enumerate}
	\item  Pick any subgroup $R\leq S$ and isomorphism $(\psi,\tau)\in\fX(R,P)_\iso$.  As $Q$ is receiving, the morphism $(\varphi\psi,\sigma\tau)\in\fX(R,Q)$ has an extension $(\widetilde{\varphi\psi},\sigma\tau)\in\fX(N_{(\varphi\psi,\sigma\tau)},S)$.  We claim that $N_{(\psi,\tau)}\leq N_{(\varphi\psi,\sigma\tau)}$ and $\widetilde{\varphi\psi}(N_{(\psi,\tau)})\leq\widetilde\varphi(N_S(P))$:  If we can show these to be true, $(\widetilde\varphi,\sigma)^{-1}\circ(\widetilde{\varphi\psi},\sigma\tau)\big|_{N_{(\psi,\tau)}}$ will define our desired extension of $(\psi,\tau)$.
	
	If $n\in N_{(\psi,\tau)}$, by definition we have $(\psi c_n\psi^{-1},\tau\ell_n\tau^{-1})=(c_s,\ell_s)$ for some $s\in N_S(P)$.  Then $\left((\varphi\psi) c_n(\varphi\psi)^{-1},(\sigma\tau)\ell_n(\sigma\tau)^{-1}\right)=(\varphi c_s\varphi^{-1},\sigma\ell_s\sigma^{-1})$, and we want to claim this is equal to $(c_{s'},\ell_{s'})$ for some $s'\in N_S(Q)$.  This is true if and only if $s\in N_{(\varphi,\sigma)}$, but we have chosen $(\varphi,\sigma)$ such that the extender is all of $N_S(P)$, and thus we have proved the first claim.
	
	To prove the second claim, note that $(\widetilde\varphi,\sigma)$ sends $Z_S(P;X)$ into $Z_S(Q;X)$, and as $P$ is fully $X$-centralized, we conclude that this is an isomorphism and therefore $Z_S(Q;X)\leq\widetilde\varphi(N_S(P))$.  Now for $n\in N_{(\psi,\tau)}$, we again have $(\psi c_n\psi^{-1},\tau\ell_n\tau^{-1})=(c_s,\ell_s)$.  Then we have
	\[
		(c_{\widetilde{\varphi\psi}(n)},\ell_{\widetilde{\varphi\psi}(n)})=(\varphi c_s\varphi^{-1},\sigma\ell_s \sigma^{-1})=(c_{\widetilde\varphi(s)},\ell_{\widetilde\varphi(s)})
	\]
	which implies that $\widetilde{\varphi\psi}(n)$ differs from  the element $\widetilde\varphi(s)\in\widetilde\varphi(N_S(P))$ by an element of $Z_S(Q;X)$, which we've already noted lies in the image of $\widetilde\varphi$.  Therefore $\widetilde{\varphi\psi}(N_{(\psi,\tau)})\leq\widetilde\varphi(N_S(P))$, completing the proof.
	\item  To see that $P$ is receiving, it suffices to show that $P$ is fully $X$-centralized.  As $\widetilde\psi$ takes $N_S(P)$ into $N_S(Q)$, it restricts to an injection of $N_S(P;X)$ into $N_S(Q;X)$ and the assumption that $P$ is fully $X$-centralized forces this to be an isomorphism.  As $Z_S(P;X)\leq N_S(P;X)$, and similarly for $Q$, the inverse to the isomorphisms of $X$-normalizers takes the $X$-centralizer of $Q$ into the $X$-centralizer of $P$.  But $Q$ is receiving, and hence fully $X$-centralized, so the $X$-centralizers are isomorphic, and $P$ is fully $X$-centralized as well.
	
	To verify the Sylow condition, note that we've already seen that $Q$ is fully $X$-normalized and $X$-centralized, so we have
	\[
	\left|\cF^\fX_S(P)_0\right|=|N_S(P;X)|/|Z_S(P;X)|=|N_S(Q;X)|/|Z_S(Q;X)|=\left|\cF^\fX_S(Q)_0\right|
	\]
	But we've also seen that $\cF^\fX_S(Q)_0\in\Syl_p\left(\cF^\fX(Q)_0\right)$, so the isomorphism $\cF^\fX (Q)\cong\cF^\fX(P)$ gives the result.	
	\end{enumerate}
	The proofs of Items 3 and 4 are the same as that of Item 2, replacing every mention of $X$-normalizers with centralizers and normalizers, respectively.
	\end{proof}
	\end{prop}
	
	We finish this section with another a seeming further simplification of the Saturation Axioms, though of course it will turn out to be equivalent.  This version is the analogue of the definition in, for example, \cite{OnofreiStancuCharacteristicSubgroup} (I first became aware of this simplification of the Saturation Axioms for fusion systems by reading an unpublished paper of Stancu).
	
	\begin{theorem}\label{more saturation simplification}
	The fusion action system $\fX$ is saturated if and only if $S$ is fully automized and every fully normalized subgroup of $S$ is receiving.
	\begin{proof}
		If $\fX$ is saturated than the original Axioms imply these weakened conditions, so we must just prove the opposite implication.  Using the simplified Axioms of  Theorem \ref{RS axioms}, it suffices to show that every fully normalized subgroup is fully automized.
		
		Suppose otherwise, and let $P\leq S$ be of maximal order with the property that $P$ is fully normalized and  $\Aut_S(P;X)\notin\Syl_p(\fX(P))$.  By the assumption that $S$ itself is fully automized, we must have that $P$ is a proper subgroup.  Then there is some some $p$-power element $(\varphi,\sigma)\in \fX(P)-\Aut_S(P;X)$ that normalizes $\Aut_S(P;X)$.   We claim that $N_{(\varphi,\sigma)}=N_S(P)$.  Indeed, if $n\in N_S(P)$, we have $(\varphi c_n\varphi^{-1},\sigma\ell_n\sigma^{-1})=(c_{n'},\ell_{n'})$ for some other $n'\in N_S(P)$ by since $(\varphi,\sigma)$ normalizes $\Aut_S(P;X)$, which proves the claim.
		
		Now by assumption $P$ is receiving, so there is some extension $(\widetilde\varphi,\sigma)\in\fX(N_S(P))$ of $(\varphi,\sigma)$.  We may assume without loss of generality that $(\widetilde\varphi,\sigma)$ has $p$-power order, as $(\varphi,\sigma)$ does.  Pick $(\psi,\tau)\in\fX(N_S(P),Q)_\iso$ such that $Q$ is fully normalized in $\fX$.  Then $^{(\psi,\tau)}(\widetilde\varphi,\sigma)\in\fX(Q)$ is a $p$-element, and as $Q$ has order strictly larger than $P$, we conclude that (by adjusting $(\psi,\tau)$ as necessary), $^{(\psi,\tau)}(\widetilde\varphi,\sigma)=(c_n,\ell_n)\in\Aut_S(Q;X)$ for some $n\in N_S(Q)$.  Then $c_n\psi(P)=\psi\widetilde\varphi(P)=\psi(P)$ since $\widetilde\varphi(P)=\varphi(P)=P$, so $n\in N_S(\psi(P))$.  Since $N_S(\psi(P))\leq\psi\left(N_S(P)\right)$ (because $Q$ was chosen to be fully normalized), we can therefore consider the element $m=\psi^{-1}(n)\in N_S(P)$, and we claim that $(\widetilde\varphi,\sigma)=(c_m,\ell_m)\in\fX(N_S(P))$.  We calculate for $a\in N_S(P)$
		\[
		mam^{-1}=\psi^{-1}(n\psi(a)n^{-1})=\psi^{-1}\left(\psi\widetilde\varphi\psi^{-1}\psi(a)\right)= \widetilde\varphi(a)
		\]
		The fact that $(\psi,\tau)$ is an intertwined pair implies that $\tau\ell_m\tau^{-1}=\ell_{\psi(m)}=\ell_n$.  On the other hand, we've assumed that $\tau\sigma\tau^{-1}=\ell_n$, so the claim is verified.
		
		But this is a contradiction, for we have realized $(\varphi,\sigma)$ (after restriction to $P$) as a morphism induced by an element of $S$.  We conclude that $\Aut_S(P;X)\in\Syl_p(\fX(P))$, and thus $\fX$ is saturated.
	\end{proof}
	\end{theorem}
	
	It should be noted that the different between the proof of Theorem \ref{more saturation simplification} and the analogous result for fusion systems is almost purely formal; the above proof looks nearly identical to that found in \cite{LinckelmannIntro}.  The same phenomenon will reappear in the next Subsection, as well as later in the paper.

\subsection{Stabilizer subsystems and $K$-normalizers}\label{Puig fusion action subsystem}

The motivation for this Subsection is to describe the ``stabilizer'' fusion action subsystem of a well-chosen point of $X$; we will return to the question of why this is a worthwhile goal in Subsection \ref{stablizer p-local finite groups}.  In particular, we shall show that these stabilizer subsystems are saturated, and thus the underlying fusion systems are as well.  On the way to proving this fairly narrow result, we shall have to introduce an analogue of Puig's notion of $K$-normalizers, and end up proving a much more general Theorem.

\begin{definition}
	The fusion action system $\fX$ is \emph{transitive} if the natural group action of $\cG=\fX(1)$ on $X$ is transitive.
\end{definition}

\begin{definition}
	Given a transitive fusion action system $\fX$ and a point $x\in X$, the \emph{stabilizer fusion action system} is the fusion action system $\fX_x$ on $S_x$ whose action on $X$ is given by restriction of the original $S$-action, and whose morphisms are given by
	\[
	\fX_x(P,Q)=\left\{(\varphi,\sigma)\in\fX(P,Q)\big|\sigma(x)=x\right\}.
	\]
\end{definition}

\begin{definition}
	If $\fX$ is a transitive fusion action system, then $x\in X$ is \emph{fully stabilized} if $|S_x|\geq|S_y|$ for all $y\in X$.
\end{definition}

\begin{theorem}\label{stabilizer saturation}
	If $x\in X$ is fully stablized in the saturated fusion action system $\fX$, the stabilizer subsystem $\fX_x$ is also saturated.
\end{theorem}

	Instead of trying to prove Theorem \ref{stabilizer saturation} directly, we shall instead generalize the result twice.
	
\begin{definition}
	Given a fusion action system $\fX$ and a subgroup $H\leq G=\fX(1)$, the \emph{preimage fusion action system of $H$ in $\fX$} is the fusion action system $\fX_H$ on the subgroup $T=\ell^{-1}(S\cap H)\leq S$, with action on $X$ the restriction of the $S$-action, and with morphisms given by
	\[
		\fX_H(P,Q)\left\{(\varphi,\sigma)\in\fX(P,Q)\big|\sigma\in H\right\}
	\]
\end{definition}

	\begin{theorem}\label{preimage saturation}
		If $\fX$ is a saturated fusion action system and $H\leq \fX(1)$ is chosen such that $\ell_T\in\Syl_p(H)$, the preimage fusion action system $\fX_H$ is saturated.
	\end{theorem}
	
	If we denote by $\Sigma_{X-x}$ the subgroup of $\Sigma$ that fixes $x$, it is immediate from the definition that that the preimage fusion action system $\fX_{\Sigma_{X-x}}$ is just the stabilizer subsystem $\fX_x$.  Thus Theorem \ref{stabilizer saturation} will follow from Theorem \ref{preimage saturation} with the aid of the following lemma.
	
	\begin{lemma}\label{full stabilization in fusion action systems}
		If the fusion action system $\fX$ is transitive then $x\in X$ is fully stabilized if and only $\ell_{S_x}\in\Syl_p(\fX(1)\cap\Sigma_{X-x})$.
		\begin{proof}
		As the core $C$ is contained in all stabilizers $S_y$ for $y\in X$ and $\left|\ell_{S_y}\right|=|S_y|/|C|$, it follows immediately that $x$ if fully stabilized if and only if $\left|\ell_{S_x}\right|$ is maximal among the orders of the groups $\ell_{S_y}$.  Moreover, the Sylow Axioms for saturation imply that $\ell_S\in\Syl_ p\left(\fX(1)\right)$, and clearly $\ell_{S_y}=\ell_S\cap\left(\fX(1)_y\right)$, so we find ourselves in the following situation:

			Let $G$ be a finite group that acts transitively on the finite set $X$.  For $S\in\Syl_ p(G)$ and $x\in X$, we have $|S_x|\geq |S_y|$ for all $y\in X$ if and only if $S_x\in\Syl_ p(G_x)$.  This result is easy to see, but we include the proof for the sake of completeness.
			
			If $x\in X$ is such that the order of $S_x$ is maximal, let $T$ be a Sylow subgroup of $G_x$ that contains $S_x$.  As $T$ is a $ p$-subgroup of $G$ and $S\in\Syl_ p(G)$, there is some $g\in G$ such that $^gT\leq S$.  We have $^gG_x=G_{g\cdot x}$, so $^gT\leq S_{g\cdot x}$.  The assumption on the maximal order of the stabilizer of $x$ then implies that $\big|{}^gT\big|\leq\big|S_x\big|$, from which the assumption that $S_x\leq T$ implies that $S_x=T\in\Syl_ p(G_x)$.
			
			Conversely, if $S_x\in\Syl_ p(G_x)$, for any $y\in X$, pick $g\in G$ such that $g\cdot y=x$.  Then the fact that $^g G_y=G_x$ implies that $^g S_y$ is a $ p$-subgroup of $G_x$, and hence subconjugate to $S_x$ by the Sylow assumption.  Thus $\big|S_y\big|\leq\big|S_x\big|$, as desired.
		\end{proof}
	\end{lemma}
	
	We digress briefly from the flow of this section to note a corollary that intuitively ``should'' be true.
	
	\begin{cor}\label{stabilizer subconjugacy}
		Let $\fX$ be a transitive saturated fusion action system and $x$ a fully stabilized point of $X$.  Then for every $y\in X$, the stabilizer $S_y$ is $\cF^\fX$-subconjugate to $S_x$.  In particular, the stabilizers of distinct fully stabilized points of $X$ are isomorphic.
		\begin{proof}
		First note that if $(\varphi,\sigma)$ is a morphism of $\fX$ such that $\sigma(y)=x$ and $s\in S_y$, then $\varphi(s)$ lies in $S_x$ when defined.  This is simply a restatement of the fact that $(\varphi,\sigma)$ is an intertwined pair:
		$
		\varphi(s)\cdot x=\varphi(s)\cdot\sigma(y)=\sigma(s\cdot y)=\sigma(y)=x
		$.
		
		Thus the result will follow from the Extension Axiom if we can find some $\sigma\in\fX(1)$ such that $\sigma(y)=x$ and $S_y\leq N_{(1,\sigma)}$.  
		As $\fX$ is assumed to be transitive, there exists a $\sigma\in\fX(1)$ such that $\sigma(y)=x$.  Then the group
		$
		H:=\left\{\sigma\circ\ell_s\circ\sigma^{-1}\big| s\in S_y\right\}
		$
		is a $ p$-subgroup of $\fX(1)_x$.  By Proposition \ref{fully stabilized characterization} (a), $\ell_{S_x}\in\Syl_ p\left(\fX(1)_x\right)$, so without loss of generality we may assume that we have chosen $\sigma$ such that that $H\leq \fX(1)_x$ and $\sigma(y)=x$.  But then $S_y\leq N_{(1,\sigma)}$, as desired.
		\end{proof}
	\end{cor}
	
	Returning to our overall goal for this section, we must generalize some more.
	
	\begin{definition}
		For $X$ an $S$-set and $P\leq S$, let $\Aut(P;X)$ denote the group of all pairs
		\[
		\left\{(\varphi,\sigma)\big|\varphi\in\Aut(P),\sigma\in\Sigma_X,\varphi\textrm{ and }\sigma\textrm{ are intertwined}\right\}.
		\]
		Similarly, let $\Aut_S(P;X)$ denote the subgroup consisting of all pairs of the form $(c_s,\ell_s)$ for $s\in N_S(P)$.
	\end{definition}
	
	\begin{definition}
		Let $K$ be a subgroup of $\Aut(P;X)$.  The \emph{$K$-normalizer of $P$ in $S$} is the group
		\[
		N_S^K(P)=\left\{n\in N_S(P)\big|(c_n,\ell_n)\in K\right\}
		\]
		We denote by $\Aut_S^K(P)$ the group $\Aut_S(P;X)\cap K$, so that $\Aut_S^K(P)\cong N_S^K(P)/Z_S(P;X)$.  Similarly set $\Aut_\fX^K(P):=\fX(P)\cap K$.
	\end{definition}

	\begin{definition}
		Given a fusion action system $\fX$, a subgroup $P\leq S$, and some $K\leq\Aut(P;X)$, the \emph{$K$-normalizer subsystem of $P$ in $\fX$} is the fusion action subsystem $\fN:=N_\fX^K(P)$ on $N_S^K(P)$ whose morphisms are given by
		\[
		\fN(Q,R)=\left\{(\varphi,\sigma)\in\fX(P,Q)\big|\exists(\widetilde\varphi,\sigma)\in\fX(PQ,PR)\textrm{ s.t. }\widetilde\varphi|_Q=\varphi\textrm{ and }(\widetilde\varphi|_P,\sigma)\in K\right\}
		\]
	\end{definition}
	
	The goal is to relate this notion to preimage fusion action systems, and thus back to stabilizer subsystems.  Of course, we must still find conditions that will make $K$-normalizer subsystems saturated.
	
	Note that if $(\varphi,\sigma)$ is any intertwined morphism pair from $P$ to $Q$, and $K\leq\Aut(P;X)$, then $^{(\varphi,\sigma)}K:=(\varphi,\sigma)\circ K\circ (\varphi,\sigma)^{-1}$ is a subgroup of $\Aut(Q;X)$.  We may thus compare the respective normalizer groups $N_S^K(P)$ and $N_S^{^{(\varphi,\psi)}K}(Q)$, or at least their orders.

	\begin{definition}
		For a fusion action system $\fX$ on $S$, $P\leq S$, and $K\leq\Aut(P;X)$, we say that $P$ is \emph{fully $K$-normalized in $\fX$} if $\left|N_S^K(P)\right|\geq\left|N_S^{^{(\varphi,\psi)}K}(Q)\right|$ for all $(\varphi,\sigma)\in\fX(P,Q)_\iso$.
	\end{definition}
	
	Let us record a basic fact about $K$-normalizers to get a better sense of how their use will play out.
	
	\begin{lemma}\label{K-normalizers to K-normalizers}
		If $(\varphi,\sigma)\in\fX(P,Q)_\iso$ has some extension $(\widetilde\varphi,\sigma)\in\fX(\widetilde P,S)$, then for any $n\in\widetilde P\cap N_S^K(P)$ we have $\widetilde\varphi(n)\in N_S^{^{(\varphi,\sigma)}K}(Q)$.
		\begin{proof}
		 If $n\in N_S^K(P)$ then $(c_n,\ell_n)\in K$, so it will suffice to show that if $\widetilde\varphi(n)$ is defined then $(c_{\widetilde\varphi (n)},\ell_{\widetilde\varphi (n)})=(\varphi,\sigma)(c_n,\ell_n)(\varphi,\sigma)^{-1}$ on $Q$.  Any $q\in Q$ can be written uniquely as $\varphi(p)$ for $p\in P$, and the fact that $\widetilde\varphi$ is an extension of $\varphi$ implies that 
		\[
		c_{\widetilde\varphi (n)}(q)=\widetilde\varphi(n)\widetilde\varphi(p)\widetilde\varphi(n)^{-1}=\widetilde\varphi(npn^{-1})=\varphi(npn^{-1})=(\varphi\circ c_n\circ\varphi^{-1})(q)
		\]
		while the fact that $\ell_{\widetilde\varphi(n)}=\sigma\circ\ell_n\circ\sigma^{-1}$ follows directly from the assumption that $(\widetilde\varphi,\sigma)$ is an intertwined pair.
		\end{proof}
	\end{lemma}

	\begin{cor}
		If $(\varphi,\sigma)\in\fX(P,Q)_\iso$ extends to $(\widetilde\varphi,\sigma)\in\fX(P\cdot N_S^K(P),S)$, then if $P$ is fully $K$-normalized we have that $Q$ is fully $^{(\varphi,\sigma)}K$-normalized.
		\begin{proof}
		Lemma \ref{K-normalizers to K-normalizers} implies that in this situation $N_S^K(P)$ injects into $N_S^{^{(\varphi,\sigma)}K}(Q)$; the assumption that $P$ is fully $K$-normalized implies that this is an isomorphism and the result follows.
		\end{proof}
	\end{cor}

	\begin{theorem}\label{K-normalizer saturation}
		With the notation as above, if $P$ is fully $K$-normalized in the saturated fusion action system $\fX$, the $K$-normalizer subsystem $\fN$ is saturated as well.
	\end{theorem}

	Before we begin the proof of Theorem \ref{K-normalizer saturation} let us see how this implies Theorem \ref{preimage saturation}.  Given $H\leq\fX(1)$, we can view $H$ as a subgroup of $\Aut(1;X)$, and thus consider the $H$-normalizer $N_S^H(1)$.  By definition, this is the subgroup of $X$ whose elements $n$ are characterized by $\ell_n\in H$.  Therefore the $H$-normalizer fusion action subsystem $\fN=N_\fX^H(1)$ is a fusion action system on the same group $T$ on which the preimage fusion action subsystem $\fX_H$ is based.  Moreover, $(\varphi,\sigma)\in\fN(P,Q)$ if and only if $\sigma\in H$, so we see that in fact these two subsystems are equal.  We are left to understand in what circumstances the identity subgroup is fully $H$-normalized in $\fX$, which the following Lemma will demonstrate is exactly the condition placed on $H$ in Theorem \ref{preimage saturation}:
	
	\begin{lemma}\label{full K-normalization conditions}
		Given a saturated fusion action system $\fX$, $P\leq S$, and $K\leq\Aut(P;X)$, then $P$ is fully $K$-normalized in $\fX$ if and only if $P$ is fully $X$-centralized in $\cF^\fX$ and $\Aut_S^K(P;X)\in\Syl_p(\Aut_\fX^K(P;X))$.
	\begin{proof}
		First suppose that $P$ is fully $X$-centralized in $\cF^\fX$ and $\Aut_S^K(P;X)\in\Syl_p(\Aut_\fX^K(P;X))$.  For any $(\varphi,\sigma)\in\fX(P,Q)_\iso$ it is immediate that $\Aut_\fX^K(P;X)\cong\Aut_\fX^{^{(\varphi,\sigma)}K}(Q;X)$.  The second assumption on $P$ then implies that $\left|\Aut_S^K(P;X)\right|\geq\left|\Aut_\fX^{^{(\varphi,\sigma)}K}(Q;X)\right|$, and the assumption that $P$ is fully $X$-centralized lets us write
		\[
		\left|N_S^K(P)\right|=\left|Z_S(P;X)\right|\cdot\left|\Aut_S^K(P;X)\right|\geq\left|Z_S(Q;X)\right|\cdot\left|\Aut_\fX^{^{(\varphi,\sigma)}K}(Q;X)\right|=\left|N_S^{^{(\varphi,\sigma)}K}(Q;X)\right|
		\]
		so that $P$ is fully $K$-normalized.
		
		Now suppose that $P$ is fully $K$-normalized, and pick an $\cF^\fX$-conjugate subgroup $Q$ that is fully normalized.  Without loss of generality [cite result] we may assume that we have an isomorphism $(\varphi,\sigma)\in\fX(P,Q)_\iso$ that extends to $(\widetilde\varphi,\sigma)\in\fX(N_S(P),N_S(Q))$.  Lemma \ref{K-normalizers to K-normalizers} implies that $\widetilde\varphi(N_S^K(P))\leq N_S^{^{(\varphi,\sigma)}K}(Q)$, so by the assumption that  $P$ is fully $K$-normalized we must have that $\widetilde\varphi|_{N_S^K(P)}$ is actually an isomorphism.  The $X$-centralizer of any subgroup is contained in any $K$-normalizers (indeed, $Z_S(P;X)=N_S^{\{\id_P\}}(P)$ is the minimal $K$-normalizer group), we see that $\widetilde\varphi$ induces an isomorphism on $X$-centralizers.  The fact that $Q$ is fully $X$-centralized implies then that $P$ is as well.
		
		To verify the Sylow condition, note that $(\varphi,\sigma)\circ\Aut_S^K(Q;X)\circ(\varphi,\sigma)^{-1}$ is a $p$-subgroup of $\Aut_\fX^{^{(\varphi,\sigma)}K}(Q;X)\leq\fX(P)$.  As $Q$ was taken to be fully normalized, $\Aut_S(Q;X)\in\Syl_p(\fX(P))$, so we conclude that there is some $(\psi,\tau)\in\fX(Q)$ such that
		\[
		(\varphi,\sigma)\circ\Aut_S^K(Q;X)\circ(\varphi,\sigma)^{-1}\leq(\psi,\tau)\circ\Aut_S(Q;X)\circ(\psi,\tau)^{-1} \cap\Aut_\fX^{^{(\varphi,\sigma)}K}(Q;X)\in\Syl_p(\Aut_\fX^{^{(\varphi,\sigma)}K}(Q;X))
		\]
		so that
		\[
		\Aut_S(Q;X)\cap\Aut_\fX^{\left(^{(\varphi,\sigma)}K\right)^{(\psi,\tau)}}(Q;X)=\Aut_S^{\left(^{(\varphi,\sigma)}K\right)^{(\psi,\tau)}}(Q;X)\in\Syl_p\left(\Aut_\fX^{\left(^{(\varphi,\sigma)}K\right)^{(\psi,\tau)}}(Q;X)\right)
		\]
		Now the assumption that $P$ is fully $K$-normalized implies that $\left|\Aut_S^K(Q;X)\right|\geq\left|\Aut_S^{\left(^{(\varphi,\sigma)}K\right)^{(\psi,\tau)}}(Q;X)\right|$, from which it follows that $\Aut_S^K(Q;X)\in\Syl_p(\Aut_\fX^K(Q;X))$.
	\end{proof}
	\end{lemma}

	As the identity subgroup is always fully $X$-centralized, we see that the primage fusion action subsystem $\fX_H$ is saturated precisely when the condition of Theorem \ref{preimage saturation} is met, i.e., when $\ell_T\in\Syl_p(H)$.
	
	\begin{cor}\label{extension to K-normalizers}
		If $(\varphi,\sigma)\in\fX(P,Q)_\iso$ and $K\leq\Aut(P;X)$ are given such that $Q$ is fully $^{(\varphi,\sigma)}K$-normalized.  If $\fX$ is saturated, then there exist $(\psi,\tau)\in\fX(P\cdot N_S^K(P),S)$ and $(\chi,\nu)\in\fX(P)$ such that $(\psi|_P,\tau)=(\varphi\chi,\sigma\nu)$.
		\begin{proof}
		Let $N=\overline{N_S^K(P)}$ denote the image of $N_S^K(P)$ in $\fX(P)$, so that $(\varphi,\sigma)\circ N\circ(\varphi,\sigma)^{-1}$ is a $p$-subgroup of $\Aut_\fX^{^{(\varphi,\sigma)}K}(Q;X)$.  By Lemma \ref{full K-normalization conditions}, the assumption on $Q$ implies that $\Aut_S^{^{(\varphi,\sigma)}K}(Q;X)\in\Syl_p\left(\Aut_\fX^{^{(\varphi,\sigma)}K}(Q;X)\right)$, so there is some $(\chi',\nu')\in\Aut_\fX^{^{(\varphi,\sigma)}K}(Q;X)$ such that $(\chi'\varphi,\nu'\sigma)\circ N\circ(\chi'\varphi,\nu'\sigma)^{-1}\leq\Aut_S^{^{(\varphi,\sigma)}K}(Q;X)$.  Thus $N\leq N_{(\chi'\varphi,\nu'\sigma)}$, so using the fact that $Q$ is fully $X$-centralized by Lemma \ref{full K-normalization conditions} and thus receiving, we have the existence of an extension $(\psi,\tau)\in\fX(P\cdot N_S^K(P),S)$ of $(\chi'\varphi,\nu'\sigma)$.  Finally, define $(\chi,\nu):=(\varphi^{-1}\chi'\varphi,\sigma^{-1}\nu'\sigma)$ to get the final result.
		\end{proof}
	\end{cor}

	We now turn to the proof of Theorem \ref{K-normalizer saturation}, by which we mean that one can now observe that the setup of $K$-normalizers of fusion action systems is so closely analogous to that of fusion systems that in fact Puig's original proof of the analogous result in \cite{PuigFrobeniusCategories}, or, for example, Linckelmann's version of the same result in \cite{LinckelmannIntro} goes through mutatis mutandis (meaning basically appending a permutation of $X$ to every morphism of the fusion system).  Note that in order to appeal to Linckelmann's proof we make use of the simplified Saturations Axioms of Theorem \ref{more saturation simplification}.  As reproducing the proofs almost word for word is not particularly enlightening in this context, the interested reader should refer to the original sources.

\section{$ p$-local finite group actions}\label{p-local finite group actions section}
	\subsection{The ambient case}\label{ambient action systems}
	
	Let us return for this Subsection to the situation where we are given a finite group $G$, a Sylow $S\in\Syl_ p(G)$, the fusion system $\cF_G$, and $X$ a $G$-set.  We wish to reconstruct the homotopy type of the Borel construction $B_GX:=EG\times_G X$, at least up to $ p$-completion, with a minimum of $ p'$-data.  This section reproduces some results of \cite{BLO1} in the context of fusion action systems arising from ambient groups.  
	
		The game we want to play in reconstructing the $ p$-completed homotopy type of $B_GX$ is to look for a new category that will both allow us to construct $B_GX^\wedge_ p$ but that does not contain too much extra information.  In some sense, the transporter system $\cT_G=\cT_S(G)$ contains all the information of $G$ that we care about (for instance, the natural functor $\cB G\to\cT_G$ that sends the unique object of $\cB G$ to the subgroup $\left\{1\right\}$ induces a homotopy equivalence $\left|\cT_G\right|\simeq BG$).  The problem of course is that $\cT_G$ contains too much information.  The goal then becomes figuring out the right amount of data of $\cT_G$ to forget and still be able to understand the Borel construction.
	
	The first way of forgetting information of $\cT_G$ is to consider full subcategories whose sets of objects are closed under $G$-conjugacy and overgroups.  In other words, we simply throw out all sufficiently small subgroups and the information of their $\Hom$-sets.  Exactly which subgroups we will allow must depend somehow on the fusion data of $G$ and the action of $G$ on $X$:

	\begin{definition}
	A $ p$-subgroup $P\leq S$ is \emph{$ p$-centric at $X$} if $Z(P;X)\in\Syl_ p(Z_G(P;X))$.  Equivalently, $P$ is $ p$-centric at $X$ if $Z_G(P;X)=Z(P;X)\times Z_G'(P;X)$ for some (necessarily unique) $ p'$-group $Z_G'(P;X)$.
	\end{definition}

	We shall generally call such a subgroup \emph{$X$-centric} and omit mentioning the prime $ p$.  This frees up our nomenclature so we can recall
	
	\begin{definition}
		A $ p$-subgroup $P\leq S$ is \emph{$ p$-centric} if $Z(P)\in\Syl_ p(Z_G(P))$, or equivalently if $Z_G(P)=Z(P)\times Z_G'(P)$ for some (again, unique) $ p'$-group $Z_G'(P)$.  This is equivalent to saying that $P$ is $X$-centric for $X=*$ the trivial $S$-set.
	\end{definition}
	
	Note that the condition of being $X$-centric is determined purely by fusion data and makes no reference to the fusion action system $\fX_G$.  This motivates the following definition:

	\begin{definition}
		Given a saturated fusion system $\cF$ on $S$ and an $S$-set $X$, a subgroup $P\leq S$ is \emph{$\cF$-centric at $X$} if $Z(Q;X)=Z_S(Q:X)$ for all $Q\cong_\cF P$.
	\end{definition}
	
	\begin{remark}
	Actually, if we're fixing the fusion system on $S$ without going straight to a fusion action system, we really should require that the $S$-set $X$ is \emph{$\cF$-stable}.  By this we mean that for all $\cF$-conjugate subgroups $P$ and $Q$ we have $\left|X^P\right|=\left|X^Q\right|$, or a couple of other equivalent conditions.  If we give ourselves an ambient group $G$, or even a fusion action system, the $\cF$-stability condition is automatically realized, so we shall not concern ourselves with it too much.
	\end{remark}

	In the presence of an ambient group $G$, the two notions of $X$-centricity coincide:

	\begin{prop}
		If $G$ is a finite group that acts on $X$ and $S\in\Syl_ p(G)$, then a subgroup $P\leq S$ is $ p$-centric at $X$ if and only if it is $\cF_G$-centric at $X$.
		\begin{proof}
		First suppose that $P$ is $ p$-centric at $X$.  As $Z_G(P;X)=Z(P;X)\times Z_G'(P;X)$ for some uniquely defined $ p'$-subgroup $Z_G'(P;X)$, it is clear that $Z_S(P;X)=Z(P;X)$.  If $g\in G$ is such that $^gP\leq S$, the fact that $Z_G({}^gP;X) ={}^gZ_G(P;X)$ immediately shows that $Z_G({}^gP;X)=Z({}^gP;X)\times{}^gZ_G'(P;X)$, and we have the desired conclusion for $^gP$.
		
		Conversely, suppose that $P$ is $\cF_G$-centric at $X$, and pick $T\in\Syl_ p(Z_G(P;X))$ such that $Z(P;X)\leq T$.  As $S\in\Syl_ p (G)$, there is some $g\in G$ such that $^gP\leq{}^gT\leq S$.  Therefore $^gT\leq Z_S({}^gP;X)=Z({}^gP;X)$ by the $\cF$-centricity at $X$ of $P$  and we conclude that $T=Z(P;X)$, as desired.
		\end{proof}
	\end{prop}

	The main purpose of introducing the abstract fusion-centric way of thinking of $X$-centricity at this point is that it makes certain results cleaner to prove:

	\begin{prop}
	If $X$ and $Y$ are $\cF$-stable $S$-sets and there is a surjective map of $S$-sets $f:X\to Y$, then every $Y$-centric subgroup is also $X$-centric.
	\begin{proof}
	Let $C$ and $D$ be the cores of $X$ and $Y$, respectively.  The existence of such a surjection forces $C\leq D$. Thus $Z_S(P)\cap C\leq Z_S(P)\cap D$, so if $Z_S(P)\cap D=Z(P)\cap D$ then $Z_S(P)\cap C=Z(P)\cap C$.  The result follows.
	\end{proof}
\end{prop}

	\begin{cor}\label{centric implies X-centric}
	$ p$-centricity implies $X$-centricity for all $\cF$-stable $X$.
	\end{cor}

	In particular, every $ p$- or $\cF$-centric (both henceforth ``centric'') subgroup of $S$ is automatically $X$-centric for every $\cF$-stable $X$.  Moreover, the more faithful $X$ is (so the smaller the core $C$ is), the easier it is for subgroups of $S$ to be $X$-centric, to the point where if $X$ is a faithful $S$-set, \emph{every} $P\leq S$ is $X$-centric.
	
	\begin{definition}
		Let $\cT_G^{cX}$ denote the full subcategory of $\cT_G$ whose objects are the $X$-centric subgroups.  Similarly, for $\fX$ an abstract action with finite set $X$, let $\fX^{cX}$ be the full subcategory on the $X$-centric subgroups.
	\end{definition}
	
	The first thing we must show is that we have not lost too much information by this restriction.  Let $\XX:\cT_G^{cX}\to\TOP$ be the functor
	\[
	\xymatrix{
	P\ar@{|..>}[r]\ar[d]_g&		X\ar[d]^{\ell_g}\\
	Q\ar@{|..>}[r]&				X
	}
	\]
	
	\begin{prop}\label{transporter Borel}
	There exists a mod-$ p$ equivalence $\hocolim_{\cT_G^{cX}}\XX\simeq_ p B_GX$.
	\begin{proof}
			Let $\cE G$ be the contractible groupoid (meaning there is precisely one morphism from any object to any other object) with objects $g\in G$, $\cK_{cX}$ the category associated to the poset of $X$-centric subgroups, and $K_{cX}$ the associated $G$-simplicial complex.  We view the $G$-set $X$ as a discrete $G$-category.  If $\cG(\XX)$ is the Grothendieck construction associated to the functor $\XX$, we have $\left|\cG(\XX)\right|\simeq\hocolim_{\cT_G^{cX}}\XX$ by \cite{ThomasonHomotopyColim}, so we will use the Grothendieck construction as our model for the homotopy colimit.
			
	There is a functor $\cE G\times_G(X\times\cK_{cX})\to\cG(\XX)$ given by
	\[
	\xymatrix{
   		[g,(x,P)]\ar@{|..>}[r]\ar[d]_{[g\mapsto hg,(id_x,P\leq Q)]}&
   		({}^gP,g\cdot x)\ar[d]^h\\
   		[hg,(x,Q)]\ar@{|..>}[r]&
   		(^{hg}Q,hg\cdot x)
	}
	\]
	that is easily seen to have an inverse isomorphism of categories
	\[
	\xymatrix{
		(P,x)\ar@{|..>}[rr]\ar[d]_h&&
		[1,(x,P)]\ar[d]^{[1\mapsto h,(id_x,P\leq{}^{h^{-1}}Q)]}\\
		(Q,h\cdot x)\ar@{|..>}[r]&
		[1,(h\cdot x,Q)]\ar@{=}[r]&
		[h,(x,{}^{h^{-1}}Q)]
	}
	\]
	Thus, on taking realizations, we get a homeomorphism of spaces
	\[
		EG\times_G(X\times K_{cX})\simeq \hocolim_{\cT_G^{cX}}\XX.
	\]

	Finally, consider the natural projection map $EG\times_G(X\times K_{cX})\to EG\times_GX$; if we can show that this is a mod-$ p$ homology isomorphism, the result will follow.  
 
	By Corollary \ref{centric implies X-centric}, the collection of $X$-centric subgroups contains all centric subgroups, and it is easy to check that it is closed under $ p$-overgroups.  Thus \cite[Theorem 8.3]{dwyerhomdecomp} applies to show that $K_{cX}$ is $\FF_p$-acyclic, and therefore the natural map $EG\times_G(X\times K_{cX})\to EG\times_G X$, being a fibration with fiber $K_{cX}$, is a mod-$ p$ homology isomorphism by the Serre spectral sequence.
 	\end{proof}
	\end{prop}

	There is a more drastic way of reducing information in $\cT_G$ than simply restricting our attention to various subgroups, in which we quotient out $ p'$-information directly.  
	
	Note that there is a free right action of $Z(P;X)$ on $\cT_G(P,Q)$, so we can consider categories whose morphisms are orbits of subgroups of $Z(P;X)$ in the transporter system.  The need to lose only $ p'$-information is one reason why we have restricted our attention to the $X$-centric subgroups of $S$, as we shall now see.

	Recall that for $H$ a finite group, $O^ p(H)$ is the smallest normal subgroup of $H$ of $ p$-power index.  If $P$ is $X$-centric, $Z_G(P;X)=Z(P;X)\times Z_G'(P;X)$ for $Z_G'(P;X)$ a $ p'$-group, and we have $O^ p(Z_G(P;X))=Z'_G(P;X)$.  

 	\begin{definition}  For a $G$-set $X$, define the \emph{$ p$-centric linking action system of $S$ at $X$} to be the category $\cL_G^{cX}=\cL^{cX}_S(G)$ whose objects are $X$-centric subgroups of $S$, and where $\cL^{cX}_G(P,Q)=N_G(P,Q)/O^ p(Z_G(P;X))$.
 \end{definition}

	Just as we formed the Borel construction $B_GX$ up to $ p$-completion by considering the homotopy colimit of a functor $\cT_G^{cX}\to\TOP$, there is a similarly defined functor whose homotopy colimit is of particular interest to us in the context of linking action systems.  Let $\overline\XX:\cL_G^{cX}\to\TOP$ be the functor
	\[
		\xymatrix{
		P\ar@{|..>}[r]\ar[d]_{[g]}&
		X\ar[d]^{\ell_g}\\
		Q\ar@{|..>}[r]&
		X
	}
	\]

	By construction, the natural quotient $\pi:\cT_G^{cX}\to\cL_G^{cX}$ satisfies the conditions of \cite[Lemma 1.3]{BLO1}.  In particular, this implies

	\begin{prop}\label{linking Borel}
	The quotient $\pi:\cT_G^{cX}\to\cL_G^{cX}$ induces a mod-$ p$ equivalence:
	\[
		\hocolim_{\cL_G^{cX}}\overline\XX\simeq_ p\hocolim_{\cT_G^{cX}}\overline\XX\circ \pi.
	\]
	\end{prop}

	\begin{theorem}\label{Borel construction from p-local finite group actions}
	\[
		\hocolim_{\cL_G^{cX}}\XX\simeq_ p B_GX
	\]
	\begin{proof}
		Combine Propositions \ref{transporter Borel} and \ref{linking Borel} with the fact that $\XX=\overline\XX\circ\pi$.
	\end{proof}
	\end{theorem}

\subsection{Classifying spaces of abstract fusion actions}\label{classifying spaces of fusion action systems}
	
	\subsubsection{Stating the problem}
	
	In this section let $\fX$ be a fixed abstract fusion action system with underlying fusion system $\cF$.  The goal is to describe what a ``classifying space'' for $\fX$ should look like.  The heuristic is that $B\fX$ should recover the homotopy type of the Borel construction of $X$, up to $ p$-completion, just as in the ambient case.  
	
	The problem with the heuristic is that, without an ambient group, we lack a Borel construction to aim for.  Section \ref{ambient action systems} tells us that we could perhaps instead try to define $B\fX$ to be the homotopy colimit of some functor into $\TOP$, but again without an ambient group we do not know what the source category for such a functor should look like.  We could try to rectify the situation by developing some notion of an \emph{abstract linking action system}.  This would be a category $\cL^\fX$ associated to $\fX$ that would have the ``right'' properties so that we could construct a functor $\XX:\cL^\fX\to\TOP$, all of which would be a generalization of the ambient case of Section \ref{ambient action systems}.
	
	This will in fact be the plan of attack take, but first a detour:  In order to properly understand what is meant  by $\cL^\fX$ having the ``right'' properties, we should try to understand the space we are looking for purely in terms that can be described by the fusion action system itself.  The question then becomes what spaces we can make from the category $\fX$.  The first guess of $\left|\fX\right|$ will not give us what we want, as this relies only on the shape of the category of $\fX$ and does not take into account the fact that it should be thought of as a combination of a diagram in groups together with permutations of $X$.
	
	More accurately, we should think of $\fX$ as a diagram in groupoids via the functor $\cB_-X:\fX\to\GRPD$. For $P\leq S$, let $\cB_PX$ denote the translation category of the $P$-set $X$.  Recall that the objects of $\cB_PX$ are the points of $X$ and the morphisms are given by 
	$
		\cB_PX(x,x')=\left\{\check p\big|p\in P\textrm{ and }p\cdot x=x'\right\}
	$.
	We then define $\cB_-X$ to be the functor
	\[
	\xymatrix{
		P\ar@{|..>}[r]\ar[d]_{(\varphi,\sigma)}&	\cB_PX\ar[d]^{\cB(\varphi,\sigma)}\\
		Q\ar@{|..>}[r]&						\cB_QX
	}
	\]
	where $\cB(\varphi,\sigma)$ is the functor that acts on objects by $\sigma$ and morphisms by $\varphi$.
	
	If we then set $B_-X:\fX\to\TOP$ to be $\left|\cB_-X\right|$, we can consider the space $\hocolim_{\fX}B_-X$, which is defined with information contained in $\fX$.  
	
	Unfortunately, this space is not what we want if $S$ is nonabelian or acts nontrivially on $X$, as the following example illustrates: 
	
	 Consider the case that $X=*$ is the trivial $S$-set. In this case $\fX=\cF$ and the space in question is $\hocolim_\cF BP$, while the space we want to recover up to $ p$-completion is $EG\times_G*=BG$.  In particular, the natural map $BS\to BG$ is Sylow, which is to say if $R$ is any $ p$-group and we have a map $BR\to BG$, there is a factorization
	\[
	\xymatrix{
	BS\ar[rr]&		&	BG\\
	&	BR\ar@{-->}[ul]\ar[ur]
	}
	\]
	up to homotopy.  If $S$ is nonabelian, $\Inn(S)$ is nontrivial, and therefore $R=\Inn(S)\ltimes S$ contains $S$ as a proper subgroup.  Let $\cG=\cG(\cB-)$ be the Grothendieck category for the functor $\cB:\cF\to\GRP$, so that
	\[
		\hocolim_\cF \cB-=\left|\cG\right|
	\]
	By definition $\cG(S)=\Inn(S)\ltimes S$, so there is an obvious map of spaces $BR\to\left|\cG\right|$ that does not factor through $BS$.
	
	The problem with simply taking $\hocolim_\fX\cB_-X$ is that this construction is in some sense double-counting the noncentral elements of $S$, as well as the elements not in $K$.  In particular, any non-$X$-central $s\in S$ defines both a morphism in $\check s\in\cB S$ and a morphism $(c_s,\ell_s)\in\fX(S)$; by simply taking the homotopy colimit of $B_-X$ these separate morphisms both contribute even though they come from the same element of $S$.  In other words, $\hocolim_\fX B_-X$ is too big, in that it has too many arrows.
	
	So let's kill the offending morphisms.
	
	\begin{definition}
	 The \emph{orbit category} of $\fX$ is the category $\cO^\fX:=\cO(\fX)$ whose objects are the subgroups of $S$ and whose morphisms are given by
	 $
	 \cO^\fX(P,Q)=Q\backslash\fX(P,Q)
	 $.
	In other words, the $\Hom$-set from $P$ to $Q$ is the set orbits of the $Q$-action of $\fX(P,Q)$ given by postcomposition by $(c_q,\ell_q)$.
	
	The full subcategory of $\cO^\fX$ whose objects are the $\cF$-centric subgroups of $S$ will be denoted $\cO^{c\fX}$.
	\end{definition}
	
	The functor $ B_-X$ does not descend to $\cO^{c\fX}\to\TOP$, but because $\cO^{c\fX}$ is defined by quotienting out inner automorphisms, it is easy to see that there is a homotopy functor $\overline B_-X:\cO^{c\fX}\to ho\TOP$.  If we could find a homotopy lifting $\widetilde B$ as in Figure \ref{homotopy lifting problem} we could consider the space $\hocolim_{\cO^{c\fX}}\widetilde B_-X$ as a possibility for the classifying space of $\cF^X$.
	\begin{figure}[h!]
	\[
	\xymatrix{
		&	\TOP\ar[d]\\
		\cO^{c\fX}\ar@{-->}^{\widetilde B_-X}[ur]\ar[r]_-{\overline B_-X}&	ho\TOP
	}
	\]
	\caption{The homotopy lifting problem}\label{homotopy lifting problem}
	\end{figure}

	We shall see that the homotopy colimit of such a lifting $\widetilde B$ is, in fact, the solution we have been seeking.
	
	\begin{prop}
		In the presence of an ambient group $G$ acting on $X$ and giving rise to $\fX$,
		\[
		\hocolim_{\cO^{c\fX}}\widetilde B_-X\simeq_ p B_GX.
		\]
		\begin{proof}
			Follows immediately from Problem \ref{abstract linking problem}, whose solution is detailed in Subsection \ref{solving the problem}.
		\end{proof}
	\end{prop}
	
  	We are now in the position to fully state the problem we want to solve:
	
	\begin{problem}\label{abstract linking problem}
		Given an abstract fusion action system $\fX$, construct an \emph{abstract linking action system} associated to $\fX$.  This should be a category $\cL^\fX$ together with functors
		\[
		\pi:\cL^\fX\to\fX^{cX}\qquad\textrm{and}\qquad\XX:\cL^\fX\to\TOP
		\]
		such that
		\begin{itemize}
			\item  In the presence of an ambient group $G$, this $\cL^{cX}_G$ and the functor $\overline\XX$ described in Subsection \ref{ambient action systems} give an abstract linking action system associated to $\fX_G$.
			\item  If $\overline\pi:\cL^\fX\to\cO^{c\fX}$ is the composite $\cL^\fX\stackrel\pi\to\fX^{cX}\to\cO^{c\fX}$, then the left homotopy Kan extension of $\XX$ over $\overline\pi$ is a homotopy lifting of $\overline B_-X:\cO^{c\fX}\to ho\TOP$, which will be denoted $\widetilde B_{\cL^\fX}:\cO^{c\fX}\to\TOP$.  Consequently,
			\[
			\hocolim_{\cL^\fX}\XX\simeq\hocolim_{\cO^{c\fX}}\widetilde B_{\cL^\fX}.
			\]
		\end{itemize}
	\end{problem}
	
	\subsubsection{Solving the problem}\label{solving the problem}

	\begin{definition}
	An \emph{($X$-centric) linking action system} associated to the fusion action system $\fX$ is a category $\cL^\fX$ whose objects are the $X$-centric subgroups of $S$, together with a pair of functors
	\[
	\xymatrix{
		\cT_S^{cX}\ar[r]^-\delta&
		\cL^\fX\ar[r]^-\pi&
		\fX^{cX}.
	}
	\]
	For each $s\in N_S(P,Q)$, let $\widehat s$ be the corresponding morphism $\delta_{P,Q}(s)\in\cL^\fX(P,Q)$.  For $\fg\in\cL^\fX(P,Q)$ denote the components of $\pi_{P,Q}(\fg)$ by the pair $(c_\fg,\ell_\fg)$.  
	
	The following conditions must be satisfied:
	\begin{enumerate}
	\item[(A)] The functors $\delta$ and $\pi$ are the identity on objects.  $\delta$ is injective and $\pi$ is surjective on morphisms.  Moreover, $Z(P;X)$ acts (via $\delta$) right-freely on $\cL^\fX(P,Q)$, and $\pi_{P,Q}:\cL^\fX(P,Q)\to\fX^{cX}(P,Q)$ is the orbit map of this action.
	\item[(B)] For each $p\in P$, we have $c_{\widehat p}=c_p:P\to P$ and $\ell_{\widehat p}=\ell_p:X\to X$.
	\item[(C)] For each $\fg\in\cL^\fX(P,Q)$ and $p\in P$, the following commutes in $\cL^\fX$:
		\[
			\xymatrix
			{
				P\ar[r]^{\fg}\ar[d]_{\widehat p}&
				Q\ar[d]^{\widehat{c_\fg(p)}}\\
				P\ar[r]_\fg&
				Q
			}
		\]
	\end{enumerate}
	
	Finally, every linking action system comes naturally equipped with a functor $\XX:~\cL^\fX\to\TOP$ defined by
	\[
	\xymatrix{
		P\ar@{|..>}[r]\ar[d]_\fg&	X\ar[d]^{\ell_\fg}\\
		Q\ar@{|..>}[r]&			X
		}
	\]
	where $X$ is regarded as a discrete topological space.
	\end{definition}

	\begin{remark}  It is easy to see that:
	\begin{itemize}
		\item  If $X=*$ is the trivial $S$-set and $\fX=\cF$ is a saturated fusion system, then an abstract linking action system is the same as a centric linking system as defined in \cite{BLO2}.  In this case $\XX$ is just the trivial functor $\cL^*\to\TOP$ whose homotopy colimit is $\left|\cL^*\right|$.
		\item  In the presence of an ambient group $G$, the category $\cL_G^{cX}$ of Subsection \ref{ambient action systems} is an example of an abstract linking action system associated to $\fX_G$.  In particular, this definition satisfies the first condition of Problem \ref{abstract linking problem}.
	\end{itemize}
	\end{remark}

	Unless it is necessary to emphasize the fact that we are looking at the $X$-centric subcategories, we shall omit explicit notational reference to them.

	The remainder of this section is devoted to proving that these abstract linking action systems satisfy the second condition of Problem \ref{abstract linking problem}.  We begin with some basic properties of abstract linking action systems:

	\begin{prop}\label{unique right lifting}
		Let
		$
		\xymatrix
		{
			P\ar[r]^{(\varphi,\sigma)}&
			Q\ar[r]^{(\psi,\tau)}&
			R
		}
		$
	be a sequence of morphisms in $\fX$.  Then for any
	\[
	\fg\in\pi^{-1}_{Q,R}((\psi,\tau))\subseteq\cL^\fX(Q,R)\qquad\textrm{and}
	\qquad\widetilde{\fg\fh}\in\pi^{-1}_{P,R}((\psi\varphi,\tau\sigma))\subseteq\cL^\fX(P,R)
	\]
	there is a unique $\fh\in\pi^{-1}_{P,Q}((\varphi,\sigma))\subseteq\cL^\fX(P,Q)$ such that $\fg\fh=\widetilde{\fg\fh}$.
	\begin{proof}
		Pick any $\fh'\in\pi_{P,Q}^{-1}((\varphi,\sigma))\subseteq\cL^\fX(P,Q)$.  By Axiom (A), the $X$-center $Z(P;X)$ acts freely and transitively on $\pi_{P,R}^{-1}((\psi\varphi,\tau\sigma))\subseteq\cL^\fX(P,R)$, so there is a unique $z\in Z(P;X)$ such that
	$
		\widetilde{\fg\fh}=\fg\fh'\widehat z
	$.
	Therefore setting $\fh=\fh'\widehat z$ gives a morphism in $\pi_{P,Q}^{-1}((\varphi,\sigma))$ such that $\fg\fh=\widetilde{\fg\fh}$.
	
	To prove uniqueness, suppose that we have $\fh_1,\fh_2\in\cL^\fX(P,Q)$ such that $\fg\fh_1=\fg\fh_2$, $\fh_i$ lifts $(\varphi,\psi)$, and $\fg\fh_i=\widetilde{\fg\fh}$.  Then $c_\fg c_{\fh_1}=c_{\fg}c_{\fh_2}\in\fX(P,R)$, and since every morphism of a fusion action system is mono (the first coordinate is an injective group map and the second coordinate is invertible), we conclude that $c_{\fh_1}=c_{\fh_2}$.  Therefore by Axiom (A) there is a unique $z\in Z(P;X)$ so that $\fh_1=\fh_2\circ\widehat z$.  Then the facts that
	\[
		\widetilde{\fg\fh}=\fg\fh_1=\fg\fh_2\circ\widehat z=\widetilde{\fg\fh}\circ\widehat z
	\]
	and $Z(P;X)$ acts right-freely on $\cL^\fX(P,R)$ implies that $z=1$, so $\fh_1=\fh_2$.
	\end{proof}
	\end{prop}
	
	\begin{remark}\label{remark following proposition}
		Note that the full strength of the Axioms of a linking action system is not used to prove this proposition.  In particular, the fact that $\cL^\fX$ comes equipped with a functor $\cT_S^{cX}\to\cL^\fX$ could have been replaced with the (apparently) weaker assumption that for every $X$-centric $P\leq S$, there is an injection $P\to\cL^\fX(P)$ such that Axiom (A) holds.  This observation will sometimes be useful when trying to construct linking action systems.
	\end{remark}
	
	Proposition \ref{unique right lifting} has a number of immediate consequences:
 
	\begin{cor}
		Every morphism in $\cL^\fX$ is categorically mono.
		\begin{proof}
			Suppose that we're given morphisms $\fg\in\cL^\fX(Q,R)$ and $\fh,\fh'\in\cL^\fX(P,Q)$ such that $\fg\fh=\fg\fh'$.  Then $\pi(\fg\fh)=\pi(\fg\fh')$, and $\fg$ is a lifting of $\pi(\fg)$, so the uniqueness statement of Proposition \ref{unique right lifting} forces $\fh=\fh'$.
		\end{proof}
	\end{cor}
	
	\begin{cor}\label{isomorphisms detected in action system}
	If $\fg\in\cL^\fX(P,Q)$ is such that $(c_\fg,\ell_\fg)\in\fX(P,Q)$ is an isomorphism, then $\fg$ is itself an isomorphism.
	\begin{proof}
		Apply Proposition \ref{unique right lifting} to the sequence $\xymatrix
		{
		Q\ar[rr]^{(c_\fg,\ell_\fg)^{-1}}&&	P\ar[rr]^{(c_\fg,\ell_\fg)}&& Q
		}$
		with $\fg$ lifting $(c_\fg,\ell_\fg)$ and $\id_Q^{\cL^\fX}$ lifting the composite $\id_Q^{\fX}$.  We obtain a unique $\fh\in\cL^\fX(Q,P)$ such that $\fg\fh=\id_Q^{\cL^\fX}$.  It follows that $\fg\fh\fg=\fg=\fg\circ\id_P^{\cL^\fX}$ and thus $\fh\fg=\id_P^{\cL^\fX}$ as $\fg$ is categorically mono.  Therefore $\fh=\fg^{-1}$.
	\end{proof}
	\end{cor}
	
	\begin{notation}  We denote by $\fri_P^Q$ the morphism $\delta_{P,Q}(1)\in\cL^\fX(P,Q)$, and call this the ``inclusion'' of $P$ in $Q$ in $\cL^\fX$.
	\end{notation}
	
	\begin{cor}\label{restrictions exist}
	For any $\fg\in\cL^\fX(P,Q)$ and $X$-centric subgroups $P^*\leq P$ and $Q^*\leq Q$ such that $c_\fg(P^*)\leq Q^*$, there is a unique morphism $\res_{P^*}^{Q^*}(\fg)\in\cL^\fX(P^*,Q^*)$ such that the following diagram commutes in $\cL^\fX$:
	\[
	\xymatrix
	{
	P\ar[rr]^\fg&&	Q\\
	{P^*}\ar[u]^{\fri_{P^*}^{P}}\ar[rr]_{\res_{P^*}^{Q^*}(\fg)}&&	{Q^*}\ar[u]_{\fri_{Q^*}^Q}
	}
	\]
	\begin{proof}
		Apply Proposition \ref{unique right lifting} to the diagram
		$
		\xymatrix
		{
		P^*\ar[r]^{c_\fg|_{P^*}}&	Q^*\ar[r]^{\iota_{Q^*}^Q}& Q
		}
		$
		with $\fri_{Q^*}^Q$ lifting $\iota_{Q^*}^Q$ and $\fg\circ\fri_{P^*}^P$ lifting the composite.
	\end{proof}
	\end{cor}
	
	\begin{notation}
	The morphism $\res_{P^*}^{Q^*}(\fg)$ is called the ``restriction'' of $\fg$, and will sometimes be denoted $\fg\big|_{P^*}^{Q^*}$ or just $\fg|_{P^*}$ if the target is clear.
	\end{notation}
	
	\begin{cor}\label{factorization in linking action systems}
	Every morphism in $\cL^\fX$ factors uniquely as an isomorphism followed by an inclusion.
	\begin{proof}
	For $\fg\in\cL^\fX(P,Q)$, the morphism $(c_\fg,\ell_\fg)\in\fX(P,c_\fg(P))$ is an isomorphism in the underlying action system, so $\res_P^{c_\fg(P)}(\fg)$ is an isomorphism in $\cL^\fX$ by Corollary \ref{isomorphisms detected in action system}.  Clearly $\fg=\fri_{c_\fg(P)}^Q\circ\res_P^{c_\fg(P)}(\fg)$, and uniqueness now follows by another application of Proposition \ref{unique right lifting}.
	\end{proof}
	\end{cor}
	
	\begin{prop}\label{linking action systems epi}
	Every morphism in $\cL^\fX$ is categorically epi.
	\begin{proof}
		By Corollary \ref{factorization in linking action systems}, it suffices to show that $\fri:=\fri_P^Q$ is epi for all $X$-centric $P\leq Q$.   Moreover, it suffices to assume that $P\trianglelefteq Q$, as any inclusion of $ p$-groups can be refined to a sequence of normal inclusions. Let $\fg_1,\fg_2\in\cL^\fX(Q,R)$ be two morphisms such that $\fg_1\circ\fri=\fg_2\circ\fri$; we want to show that $\fg_1=\fg_2$.
		
		The image of $\fg_i$ in $\fX(Q,R)$ is $(c_{\fg_i},\sigma_i)$ and the assumption that $\fg_1\circ\fri=\fg_2\circ\fri$ implies that $\sigma_1=\sigma_2$.  We first show that it also follows that $c_{\fg_1}=c_{\fg_2}$.  
		
		Note that the assumption that $P\trianglelefteq Q$ implies that every $q\in Q$ determines a morphism in $\widehat q\in\cL^\fX(P)$, which is easily seen to be the restriction of $\widehat q\in\cL^\fX(Q)$.  Then the diagram 
		\[
		\xymatrix
		{
			Q\ar[r]^{\fg_i}\ar[d]_{\widehat q}&	R\ar[d]^{\widehat{c_{\fg_i}(q)}}\\
			Q\ar[r]_{\fg_i}&		R
		}
		\]
		commutes in $\cL^\fX$ by Axiom (C) for $i=1,2$.  The assumption on the $\fg_i$ also implies that their restrictions to $P$ in $\cL^\fX$ are equal; denote this common morphism by $\fh\in\cL^\fX(P,c_\fh(P))$.  We can therefore form the restriction of this diagram, which gives us
		\[
		\xymatrix
		{
			P\ar[r]^\fh\ar[d]_{\widehat q}&	c_\fh(P)\ar[d]^{\widehat{c_{\fg_i}(q)}}\\
			P\ar[r]_\fh&		c_\fg(P)
		}
		\]
		for $i=1,2$.  As all the morphisms in this restriction diagram are iso, and three of the four do not depend on choice of $i$, we conclude that $\widehat{c_{\fg_1}(q)}=\widehat{c_{\fg_2}(q)}$ for all $q\in Q$.  The assignment $q\mapsto\widehat q$ is injective, so we conclude that $c_{\fg_1}=c_{\fg_2}$ on $Q$.
		
		Thus we have $\pi(\fg_1)=\pi(\fg_2)$, so by Axiom (2), there is some $z\in Z(Q;X)$ such that $\fg_2=\fg_1\circ\widehat z$.  We then compute
		$
		\fg_1\circ\widehat z\circ \fri=\fg_2\circ\fri=\fg_1\circ\fri
		$.  The fact that $\fg_1$ is mono implies that $\widehat z\circ\fri=\fri$.  
		
		Finally, the fact that $P\leq Q$ implies that $z\in Z(P;X)$, and Axiom (3) again shows that $\widehat z\circ\fri=\fri\circ\widehat z$ (for $\widehat z$ respectively an isomorphism in $\cL^\fX$ of $Q$ and $P$).  The freeness of the right action of $Z(P;X)$ on $\cL^\fX(P,Q)$ forces $z=1$, and the result is proved.
	\end{proof}
	\end{prop}
	
	\begin{cor}
	Extensions  of morphisms in $\cL^\fX$ are unique when they exist.  In other words, for any $\fg^*\in\cL^\fX(P^*,Q^*)$ and subgroups $P^*\leq P$ and $Q^*\leq Q$, there is at most one morphism $\fg\in\cL^\fX(P,Q)$ such that the diagram
	\[
	\xymatrix
	{
		P\ar[r]^\fg&	Q\\
		P^*\ar[u]^{\fri_{P^*}^P}\ar[r]_{\fg^*}&	Q^*\ar[u]_{\fri_{Q^*}^Q}
	}
	\]
	commutes in $\cL^\fX$.
	\begin{proof}
 		If $\fg'$ is another such extension, the equalities
		$
			\fg\circ\fri_{P^*}^P=\fri_{Q^*}^Q\circ\fg^*=\fg'\circ\fri_{P^*}^P
		$
		and the fact that $\fri_{P^*}^P$ is epi force $\fg=\fg'$.
	\end{proof}
	\end{cor}
	
	\begin{prop}\label{fully normalized Sylow in linking action systems}
	If $P\leq S$ is fully normalized in $\cF$, then $\widehat{N_S(P)}\big|_P^P\in\Syl_ p\big(\cL^\fX(P)\big)$.
	\begin{proof}
		By Axiom (A) of linking action systems, $\left|\cL^\fX(P)\right|=\left|\fX(P)\right|\cdot\left|Z(P;X)\right|$ and by definition of $\Aut_S(P;X)$,
		\[
			\left|N_S(P)\right|=\left|\Aut_S(P;X)\right|\cdot\left|Z_S(P;X)\right|
			=\left|\Aut_S(P;X)\right|\cdot\left|Z(P;X)\right|.
		\]
		As $P$ is fully normalized, the Saturation Axioms imply that $\Aut_S(P;X)\in\Syl_ p\left(\fX(P)\right)$, and the result easily follows.
	\end{proof}
	\end{prop}
	
	\begin{notation}
		Let $\overline\pi$ be the composite $\cL^\fX\stackrel\pi\to\fX\to\cO^{c\fX}$.
	\end{notation}
	
	\begin{prop}\label{free target action in linking action system}
	Suppose that $\fg,\fh\in\cL^\fX(P,Q)$ are such that $\overline\pi(\fg)=\overline\pi(\fh)$.  Then there is a unique element $q\in Q$ such that $\fh=\widehat q\circ\fg$.  In other words, the map $\overline\pi_{P,Q}:\cL^\fX(P,Q)\to\cO^{c\fX}(P,Q)$ is the orbit map of the free left action of $Q$ on $\cL^\fX(P,Q)$.
	\begin{proof}
	The condition on $\fg$ and $\fh$ implies that there is some $q'\in Q$ such that
	\[
		(c_\fh,\ell_\fh)=(c_{q'},\ell_{q'})\circ(c_\fg,\ell_\fg).
	\]
	Condition (A) then implies that there is a unique $z\in Z(P;X)$ such that
	\[
		\fh=\widehat{q'}\circ\fg\circ\widehat z=\widehat{q'}\circ\widehat{c_\fg(z)}\circ\fg
	\]
	where the second equality follows from Condition (C).  Setting $q=q'\circ c_\fg(z)$ gives the desired element.  
	
	The uniqueness of $q$ is a direct consequence of the fact that $\fg$ is epi and the assignment $q\mapsto\widehat q$ is injective.
	\end{proof}
	\end{prop}

\begin{figure}[h!]
	\[
	\xymatrix
	{
		\cL^\fX\ar[d]_{\overline\pi}\ar[r]^-\XX&	\TOP\\
		\cO^{c\fX}\ar[ur]_{\widetilde B_-X}
	}
	\qquad\qquad
	\xymatrix
	{
	&	\TOP\ar[d]\\
	\cO^{c\fX}\ar@{-->}[ur]^?\ar[r]_-{\overline B_-X}&	ho\TOP
	}
	\]
	\caption{What we have vs. What we want}\label{have vs. want}
	\end{figure}
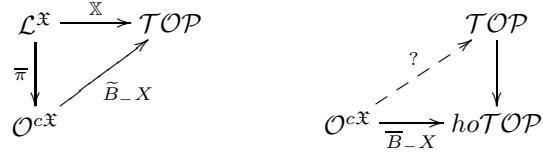

	We now have the necessary results to show that we have solved the second point of Problem \ref{abstract linking problem}.  Recall that $\widetilde B_-X:\cO^{c\fX}\to\TOP$ is the left homotopy Kan extension of $\XX:\cL^\fX\to\TOP$ over $\overline\pi:\cL^\fX\to\cO^{c\fX}$, and that we are looking for a homotopy lifting of $\overline B_-X:\cO^{c\fX}\to ho\TOP$, as depicted in Firgure \ref{have vs. want}.
	
	For ease of notation, write $\widetilde B$ for $\widetilde B_-X$ and $B$ for $\overline B_-X$.
\\\\
	\begin{prop}
		The functor $\widetilde B$ is a homotopy lifting of $B$.
	\begin{proof}
		{\bf{On objects:}}  We want to show that $\widetilde B(P)\simeq B(P)=B{}_PX$.  Recall that
		\[
			\widetilde B(P)\simeq\hocolim_{(Q,\alpha)\in(\overline \pi\downarrow P)}\XX(Q)
		\]
		where $(\overline\pi\downarrow P)$ is the overcategory of $\overline\pi$ over $P$.  The objects are pairs $(Q,\alpha)$, where $Q$ is an $X$-centric subgroup of $S$, viewed as either an object of $\cL^\fX$ or $\cO^{c\fX}$ as appropriate, and $\alpha$ is a morphism in $\cO^{c\fX}(Q,P)$.  A morphism from $(Q,\alpha)$ to $(R,\beta)$ is $\fg\in\cL^\fX(Q,R)$ such that the following 	commutes in $\cO^{c\fX}$:
	\[
	\xymatrix
	{
		Q\ar[rr]^{\overline \pi(\fg)}\ar[dr]_\alpha&
		&
		R\ar[dl]^\beta\\
		&
		P
	}
	\]

	Thomason's theorem \cite{ThomasonHomotopyColim} shows that this homotopy colimit expression of $\overline B(P)$ has the homotopy type of $|\cG(P)|$, where $\cG(P)$ is the Grothendieck category associated to $(\overline \pi\downarrow P)$ and the functor $\XX$.

 	Explicitly, $\cG(P)$ is the category whose objects are pairs $((Q,\alpha),x)$ where $(Q,\alpha)$ is an object of $(\overline\pi\downarrow P)$ and $x\in X$.  A morphism $((Q,\alpha),x)\to((R,\beta),y))$ is a morphism $\fg\in(\overline\pi\downarrow P)((Q,\alpha),(R,\beta))$---and therefore  $\fg\in\cL^\fX(Q,R)$---such that $\ell_\fg(x)=y$.  

	Let $\check\cB(P)\subseteq\cG(P)$ be the subcategory whose objects are pairs $((P,\mathrm{id}_P),x)$ for all $x\in X$, and where
	\[
		\check\cB(P)(((P,\id_P),x),((P,\id_P),y))=\left\{\widehat p\big|p\in P, p\cdot x=y\right\}.
	\]
	This category is isomorphic to $\cB{}_PX$; the claim is that the inclusion $\check\cB(P)\subseteq\cG(P)$ induces a deformation retract after realization.  This will follow if we can find a functor $\Psi:\cG(P)\to\check\cB(P)$ that is the identity on $\check\cB(P)$ together with a natural transformation $F:\Id_{\cG(P)}\Rightarrow\iota_{\check\cB(P)}^{\cG(P)}\circ\Psi$, where $\iota_{\check\cB(P)}^{\cG(P)}:\check\cB(P)\to\cG(P)$ is the inclusion functor.

	Pick some lifting of morphisms $\xi:\Mor(\cO^{c\fX})\to\Mor(\cL^\fX)$, and assume that $\xi$ sends identities to identities.  Let $\Psi:\cG(P)\to\check\cB(P)$ be the functor 
	\[
	\xymatrix{
	((Q,\alpha),x)\ar@{|..>}[r]\ar[d]_\fg&	\left((P,\id_P),\ell_{\xi(\alpha)}(x)\right)\ar[d]^{\Psi(\fg)}\\
	((R,\beta),y)\ar@{|..>}[r]&			\left((P,\id_P),\ell_{\xi(\beta)}(y)\right)
	}
	\]
	where $\Psi(\fg)$ is defined to be $\widehat p$ for the unique $p\in P$ such that
	\[
		\xymatrix{
		Q\ar[r]^\fg\ar[d]_{\xi(\alpha)}&
		R\ar[d]^{\xi(\beta)}\\
		P\ar[r]_{\widehat p}&
		P
		}
	\]
	The existence and uniqueness of $p$ follow from Proposition \ref{free target action in linking action system} and the fact that 
	$
	\overline\pi(\xi (\beta)\circ\fg)=\beta\circ\overline\pi(\fg)=\alpha=\overline\pi(\xi(\alpha))
	$
	because $\fg$ is a morphism in  $(\overline\pi\downarrow P)$.  
	
	We must check that $\widehat p$ indeed defines a morphism in $\cG(P)$:
	\[
		((P,\id_P),\ell_{\xi(\alpha)}(x))\to((P,\id_P),\ell_{\xi(\beta)} (y)).
	\]
	  To do this, we must show that
	\[
	\xymatrix
	{
		P\ar[rr]^{(c_p,\ell_p)}\ar[dr]_{\id_P}&
		&
		P\ar[dl]^{\id_P}\\
		&
		P
	}
	\]
	commutes in $\cO^{c\fX}$ and that
	$
		\ell_p\circ\ell_{\xi(\alpha)}(x)=\ell_{\xi(\beta)} (y)
	$.
	The first is obvious from the definition of $\cO^{c\fX}$, and the second follows from the fact that $\widehat p\circ \xi(\alpha)=\xi(\beta)\circ\fg$, so that
	$
		\ell_p\circ \ell_{\xi(\alpha)}(x)=\ell_{\xi(\beta)}\circ \ell_\fg(x)=\ell_{\xi(\beta)}(y)
	$
	by the assumption on $\fg$.

	Observe that $\Psi|_{\check\cB(P)}$ is the identity functor.

	Define the natural transformation $\Theta:\Id_{\cG(P)}\Rightarrow\iota_{\check\cB(P)}^{\cG(P)}\circ\Psi$ by
	\[
		\Theta(((Q,\alpha),x))=\xi(\alpha):((Q,\alpha),x)\to((P,\id_P),\ell_{\xi(\alpha)}(x)).
	\]
	That this is a morphism in $\cG(P)$ follows easily from the definition of $\cG(P)$ above and the fact that $\xi(\alpha)$ is a lift of $\alpha$ to $\cL$.  If
	$
		\fg:((Q,\alpha),x)\to((R,\beta),y)
	$
	is a morphism in $\cG(P)$, we want the diagram
	\[
	\xymatrix
	{
		((Q,\alpha),x)\ar[r]^{\fg}\ar[d]_{\xi(\alpha)}&	((R,\beta),y)\ar[d]^{\xi(\beta)}\\
		((P,\id_P),\ell_{\xi(\alpha)}(x))\ar[r]_{\Psi(\fg)}&	((P,\id_P),\ell_{\xi(\beta)}(y))		
	}
	\]
	to commute in $\cL^\fX$, which follows from the definition of $\Psi(\fg)$.
	
	Finally, note that for any object of $\check x\in\check\cB(P)$, it is easy to see that $F(\check x)$ is the identity.  Thus the geometric realization of $\Theta$ gives a homotopy that shows that $\left|\Psi\right|$ realizes $\left|\check\cB(P)\right|$ as a deformation retract of $\left|\cG(P)\right|$, and we have the result on objects.

	{\bf On morphisms:}  First observe that each $[\varphi,\sigma]\in\cO^{c\fX}(P,P')$ induces a functor $(\overline\pi\downarrow P)\to(\overline\pi\downarrow P')$:
	\[
	\xymatrix
	{
	(Q,\alpha)\ar@{|..>}[r]\ar[d]_\fg&		(Q,[\varphi,\sigma]\circ\alpha)\ar[d]^\fg\\
	(R,\beta)\ar@{|..>}[r]&				(R,[\varphi,\sigma]\circ\beta)
	}
	\]
	This in turn induces a functor $[\varphi,f]_*:\cG(P)\to\cG(P')$, defined by 
	\[
	\xymatrix{
		((Q,\alpha),x)\ar@{|..>}[r]\ar[d]_{\fg}&
		((Q,[\varphi,\sigma]\circ\alpha),x)\ar[d]^{\fg}\\
		((R,\beta),\ell_\fg x)\ar@{|..>}[r]&
		((R,[\varphi,\sigma]\circ\beta),\ell_\fg x)
	}
	\]
	Similarly, for $(\varphi,\sigma)\in\fX(P,P')$ a lift of $[\varphi,\sigma]\in\cO^{c\fX}(P,P')$, there is a functor $\check\cB_{(\varphi,\sigma)}: \check\cB(P)\to\check\cB(P')$:
	\[
	\xymatrix
	{
		((P,\id_P),x)\ar@{|..>}[r]\ar[d]_{\widehat p}&
		((P',\id_{P'}),\sigma(x)\ar[d]^{\widehat{\varphi p}})\\
		((P,\id_P),p\cdot x)\ar@{|..>}[r]&
		((P',\id_{P'}),\varphi p\cdot \sigma(x))
	}
	\]
	
	To finish the proof that $\widetilde B$ is a homotopy lifting of $B$ we must show that the diagram of functors
	\[
	\xymatrix
	{
		\check\cB(P)\ar[r]^{\iota}\ar[d]_{\check\cB_{(\varphi,\sigma)}}&
		\cG(P)\ar[d]^{[\varphi,\sigma]_*}\\
		\check\cB(P')\ar[r]_\iota&
		\cG(P')
	}
	\]
	commutes up to natural transformation.  
	
	Let $F_1$ be the functor given by the top path of the diagram.  Explicitly, $F_1$ is the functor
	\[
	\xymatrix
	{
		((P,\id_P),x)\ar@{|..>}[r]\ar[d]_{\widehat p}		&		((P,[\varphi,\sigma]),x)\ar[d]^{\widehat p}\\
		((P,\id_P),p\cdot x)\ar@{|..>}[r]				&		((P,[\varphi,\sigma]),p\cdot x)
	}
	\]
	Let the bottom path be the functor $F_2$:
	\[
	\xymatrix
	{
		((P,\id_P),x)\ar@{|..>}[r]\ar[d]_{\widehat p}&
		((P',\id_{P'}),\sigma(x))\ar[d]^{\widehat{\varphi p}}\\
		((P,\id_P),p\cdot x)\ar@{|..>}[r]&
		((P',\id_{P'}),\varphi p\cdot \sigma(x))
	}
	\]
	Choose some $\fh\in\cL^\fX(P,P')$ that lifts $(\varphi,\sigma)$, and define the natural transformation
	\[
	\Phi:F_1\Rightarrow F_2:((P,id),x)\mapsto\fh.
	\]
	To see that $\fh$ is indeed a morphism in $\cG(P')$ from $((P,[\varphi,\sigma]),x)$ to $((P',\id_{P'}),\sigma(x))$, note that 
	\[
	\xymatrix
	{
		P\ar[rr]^{\overline\pi(\fh)}\ar[dr]_{[\varphi,f]}&
		&
		P'\ar[dl]^{\id_{P'}}\\
		&
		P'
	}
	\]
	commutes because $\fh$ is a lifting of $[\varphi,\sigma]$.  Moreover, $\ell_\fh x=\sigma(x)$ because $\fh$ is a lifting of $(\varphi,\sigma)$.
	
	Finally, if $\widehat p$ is a morphism in $\check\cB(P)$ from $((P,\id_P),x)$ to $((P,\id_P),y)$, the diagram
	\[
	\xymatrix
	{
		((P,[\varphi,\sigma]),x)\ar[r]^\fh\ar[d]_{\widehat p}&
		((P',\id_{P'}),\sigma(x))\ar[d]^{\widehat{\varphi p}}\\
		((P,[\varphi,\sigma]),y)\ar[r]_\fh&
		((P',\id_{P'}),\sigma(y))
	}
	\]
	commutes in $\cG(P')$ because Axiom (C) and the fact that $c_\fh=\varphi$ ($\fh$ lifts $(\varphi,\sigma)$) together imply that
	$
		\widehat{\varphi p}\circ\fh=\fh\circ\widehat p
	$
	in $\cL^\fX(P,P')$.

   This completes the proof that $\widetilde B$ is a homotopy lifting of $B$.
  \end{proof}
 \end{prop}

	\subsection{Obstruction theory}\label{fusion action obstruction theory}
	
	Fix an abstract fusion action system $\fX$.  The material in this section does not depend on $\fX$'s being saturated.
	
	We can solve the homotopy lifting problem for $\overline\cB_-X:\cO^{c\fX}\to ho\TOP$, and thus create a space that we might call a ``classifying space'' for $\fX$, if we have an associated linking action system $\cL^\fX$.  In this section we describe how we can construct $\cL^\fX$ from the data of $\fX$, or more accurately, we describe the difficulties in doing so.  We also describe the obstructions to constructing $\cL^\fX$ uniquely.  
	
	\begin{definition}
	If there is a unique $\cL^\fX$ associated to $\fX$, the space $\hocolim_{\cL^\fX}\XX$ is the \emph{classifying space} for the fusion action system.  We denote this space by $B\fX$.
	\end{definition}
	
	\begin{notation}
		We denote by $\left(\widetilde\varphi,\widetilde\sigma\right)\in\fX(P,Q)$ a lift of $[\varphi,\sigma]\in\cO^{c\fX}(P,Q)$.
	\end{notation}
	
	\begin{definition}
	Let $\cZ^\fX:\left(\cO^{c\fX}\right)^{op}\to\cA b$ be the functor
	\[
	\xymatrix
	{
		P\ar@{|..>}[r]\ar[d]_{[\varphi,\sigma]}&
		Z(P;X)&
		Z(\widetilde\varphi(P);X)\ar[l]_{\widetilde\varphi^{-1}}&
		Z_S(\widetilde\varphi(P);X)\ar@{=}[l]\\
		Q\ar@{|..>}[r]&
		Z(Q;X)\ar[u]_{\cZ^\fX([\varphi,\sigma])}\ar@{=}[rr]&&
		Z_S(Q;X)\ar@{^(->}[u]_{\textrm{incl}}
	}
	\]
	where $\cZ^\fX([\varphi,\sigma])$ is the unique map of groups that makes the rectangle commute.
	\end{definition}
	
	\begin{prop}
		The functor $\cZ^\fX$ is well-defined.
		\begin{proof}
			The equalities come from the assumption that $P$ is $X$-centric.  We must simply show that for $(\widetilde\varphi',\widetilde\sigma')$ any another lift of $[\varphi,\sigma]$, the maps $\widetilde\varphi^{-1}$ and $(\widetilde\varphi')^{-1}$ are equal on $Z(Q;X)$.  Since that $\widetilde\varphi$ and $\widetilde\varphi'$ are both lifts of $\varphi$, there is some $q\in Q$ such that $\widetilde\varphi'=c_q\circ\widetilde\varphi$.  Therefore $c_q^{-1}$ is the identity on $Z(Q;X)$, and the result follows.
		\end{proof}
	\end{prop}
	
	Just as in the theory of $p$-local finite groups, there is a well-defined obstruction theory that governs the existence and uniqueness of linking action systems associated to a fusion action system.
	
	\begin{theorem}\label{fusion action system obstructions to p-local finite group actions}
	The data of the abstract fusion action $\fX$ determines an element $u\in\lim^3_{\cO^{c\fX}}\cZ^\fX$ that vanishes precisely when there is a linking action system $\cL^\fX$ associated to $\fX$.  Moreover, if there is a linking action system, the group $\lim^2_{\cO^x}\cZ^\fX$ acts transitively on the set of linking action systems viewed as categories over $\cO^{c\fX}$.
	\begin{proof}
	Just as the proof of Theorem \ref{K-normalizer saturation} looks almost exactly like that of the analogue for fusion systems, so too can the proof of the obstruction theory governing $p$-local finite group from \cite{BLO2} be carried over with only the obvious changes to the world of $p$-local finite group actions.  Again, we omit reproducing what amounts to identical symbols with a different interpretation, and instead refer the reader to the original source.
	\end{proof}
	\end{theorem}
	
\subsection{Linking action systems as transporter systems}\label{Linking action systems as transporter systems}

	Let $\fX$ be a saturated fusion action system and $\cL^\fX$ an associated linking action system.

	Proposition \ref{unique right lifting} can be thought of a ``unique right lifting'' lemma, and we saw that it implied many useful properties for a linking action system $\cL^\fX$.  In particular, we learned from it that all morphisms are categorically mono, that restriction is a well-defined notion (given the existence of specified ``inclusion morphisms''), and ultimately that all morphisms are epi.  This last result ``should'' have been derived from a ``unique left lifting'' lemma, but instead we used the additional structure of the inclusion morphisms to derive it.  Indeed, there is no direct left lifting analogue of Proposition \ref{unique right lifting}; instead, we have to settle with the following, which turns out to have its own uses:
	
	\begin{prop}\label{unique left pseudolifting}
			Let
		$
		\xymatrix
		{
			P\ar[r]^{(\varphi,\sigma)}&
			Q\ar[r]^{(\psi,\tau)}&
			R
		}
		$
	be a sequence of morphisms in $\fX$.  For any
	\[
	\fh\in\pi^{-1}_{P,Q}((\varphi,\sigma))\subseteq\cL^\fX(P,Q)\qquad\textrm{and}
	\qquad\widetilde{\fg\fh}\in\pi^{-1}_{P,R}((\psi\varphi,\tau\sigma))\subseteq\cL^\fX(P,R)
	\]
	there is a unique $\fg\in\cL^\fX(P,Q)$ such that $\fg\fh=\widetilde{\fg\fh}$.  Moreover, there is a unique $z\in \varphi(Z(P;X))$ such that
	$
		(c_\fg,\ell_\fg)=(\varphi c_z,\sigma\ell_z)
	$.
	\begin{proof}
		The morphism $\fh$ is epi by Proposition \ref{linking action systems epi}, so there is at most one $\fg$ such that $\fg\fh=\widetilde{\fg\fh}$.
		
		Pick any $\fg'$ lifting $(\psi,\tau)$. Thus $\fg'\circ\fh$ and $\widetilde{\fg\fh}$ have the same image in $\fX(P,R)$, and by Axioms (A) and (C), there is a unique $z\in Z(P;X)$ such that
		$
			\widetilde{\fg\fh}=\fg'\circ\fh\circ\widehat z=\fg'\circ\widehat{\varphi(z)}\circ\fh
		$.
		Then $\fg=\fg'\circ\widehat{\varphi(z)}$ gives us the existence statement. 
	\end{proof}
	\end{prop}
	
	\begin{remark}
		The main difference between the left ``lifting'' of Proposition \ref{unique left pseudolifting} and the right lifting of Proposition \ref{unique right lifting} is that this more recent result cannot actually lift the morphism $(\psi,\tau)\in\fX(Q,R)$, but only some $\varphi(Z(P;X))$-translate of it.  The difference between the two stems from the possibility that $Q\gneq\varphi(P)$, in which case $\varphi$ need not take central elements to central elements; as $Q$ can be bigger, it is possible that $\varphi(z)$ acts nontrivially on $Q$, even though $z$ acts trivially on $P$.
	\end{remark}
	
	Recall that an \emph{extension} of $\fg\in\cL^\fX(P,Q)$ is $\widetilde\fg\in\cL^\fX(\widetilde P,\widetilde Q)$ for $P\leq\widetilde P$ and $Q\leq \widetilde Q$ such that
	\[
	\xymatrix
	{
		\widetilde P\ar[r]^{\widetilde\fg}&	\widetilde Q\\
		P\ar[u]^{\fri_P^{\widetilde P}}\ar[r]_\fg&	Q\ar[u]_{\fri_Q^{\widetilde Q}}
	}
	\]
	commutes in $\cL^\fX$.  We have already seen (thanks to the fact that all morphisms of $\cL^\fX$ are both epi and mono) that extensions are unique if they exist, and we are now in the position to say when exactly they do exist:
	
	\begin{prop}\label{linking action extensions}
	Let $\fg\in\cL^\fX(P,Q)$ be an isomorphism and let $\widetilde P,\widetilde Q\leq S$ be such that $P\trianglelefteq
\widetilde P$, $Q\trianglelefteq\widetilde Q$, and
	$
		\fg\circ\delta_{P,P}\left(\widetilde P\right)\circ\fg^{-1}\leq\delta_{Q,Q}\left(\widetilde Q\right)
	$.
	Then there is a unique extension $\widetilde\fg\in\cL^\fX(\widetilde P,\widetilde Q)$ of $\fg$.
	\begin{proof}
		First suppose that $Q$ is fully $X$-centralized.   By the Extension Axiom for saturated fusion actions, the morphism $(c_\fg,\ell_\fg)\in\fX(P,Q)$ extends to some $(\varphi,\ell_\fg)\in\fX\left(N_{(c_\fg,\ell_\fg)},S\right)$.  The condition on $\widetilde P$ implies that $\widetilde P\leq N_{(c_\fg,\ell_\fg)}$ (project the condition down to $\fX$) and the condition on $\widetilde Q$ implies that $\varphi(p')\in \widetilde Q$ for all $p'\in\widetilde P$.  We can thus rename $(\varphi,\ell_\fg)$ to be its restriction in $\fX(\widetilde P,\widetilde Q)$.  For the sequence
		$
		\xymatrix
		{
			P\ar[rr]^{\left(\iota_P^{\widetilde P},\id_X\right)}&&\widetilde P\ar[rr]^{\left(\varphi,\ell_\fg\right)}&&
			\widetilde Q
		}
		$
		in $\fX$, let $\fri_P^{\widetilde P}\in\cL^\fX(P,\widetilde P)$ lift the first map and $\fri_Q^{\widetilde Q}\circ\fg$ lift the composite.  Then Proposition \ref{unique left pseudolifting} applies to give a unique extension $\widetilde\fg$ of $\fg$ as desired (though note that it need not be the case that $\widetilde\fg$ is a lift of $(\varphi,\ell_\fg)$, only a lift of a $\varphi(Z(P;X))$-translate of it).
		
		Now consider the general case, where $Q$ need not be fully $X$-centralized.  Let $R$ be fully normalized and $\cF$-conjugate to $P$ and $Q$.  For any $\fh\in\cL^\fX(Q,R)$, we have that 
		$
			\fh\circ\widehat{N_S(Q)}\big|_Q^Q\circ\fh^{-1}\leq\cL^\fX(R)
		$
		is an inclusion of a $ p$-subgroup.  As $R$ is fully normalized, Proposition \ref{fully normalized Sylow in linking action systems} states that $\delta_{R,R}\left(N_S(R)\right)$ is Sylow in $\cL^\fX(R)$, so $\fh$ can be chosen so that
		$
			\fh\circ\widehat{N_S(Q)}\big|_Q^Q\circ\fh^{-1}\leq\widehat{N_S(R)}\big|_R^R
		$.
		The subgroup $R$ is fully $X$-centralized, so the first part of this proof implies that there are morphisms $\widetilde\fh$ extending $\fh$ to $N_S(Q)$ and $\widetilde{\fh\fg}$ extending $\fh\fg$ to $\widetilde P$.  Let $\overline\fh$ be the restricted isomorphism of $\widetilde\fh$ with source $\widetilde Q$, and similarly $\overline{\fh\fg}$ the restricted isomorphism of $\widetilde{\fh\fg}$ with source $\widetilde P$.  The situation can be represented as:
		\[
		\xymatrix
		{
			N_S(Q)\ar[r]^{\widetilde\fh}&	N_S(R)&&&	\overline P\ar[r]^-{\widetilde{\fh\fg}}&	N_S(R)\\
			\widetilde Q\ar[u]^\fri\ar[r]^-\cong_-{\overline\fh}&	c_{\widetilde\fh}\left(\widetilde Q\right)\ar[u]_\fri&&&
			\widetilde P\ar[u]^=\ar[r]_-{\overline{\fh\fg}}^-\cong&	c_{\widetilde{\fh\fg}}\left(\widetilde P\right)\ar[u]_\fri\\
			Q\ar[u]^\fri\ar[r]^\cong_\fh&	R\ar[u]_\fri&&&	P\ar[u]^\fri\ar[r]_{\fh\circ\fg}&	R\ar[u]_\fri
		}
		\]
		The claim is that $c_{\widetilde{\fh\fg}}\left(\widetilde P\right)\leq c_{\widetilde\fh}\left(\widetilde Q\right)$.  Axiom (C) implies that for $p'\in\widetilde P$ and $q'\in\widetilde Q$,
		\[
			\widehat{c_{\widetilde{\fh\fg}}(p')}=\overline{\fh\fg}\circ\widehat{p'}\circ\overline{\fh\fg}^{-1}\qquad\textrm{and}\qquad
			\widehat{c_{\widetilde{\fh}}(q')}=\overline{\fh}\circ\widehat{q'}\circ\overline{\fh}^{-1}.
		\]
		Observing that every $p'\in\widetilde P$ defines a morphism in $\cL^\fX(P)$, we can restrict this to get
		\[
			\delta_{R,R}\left(c_{\widetilde{fg}}\left(\widetilde P\right)\right)=\left(\fh\circ\fg\right)\circ\delta_{P,P}
			\left(\widetilde P\right)\circ\left(\fh\circ\fg\right)^{-1}
			\leq\fh\circ\delta_{Q,Q}\left(\widetilde Q\right)\circ\fh^{-1}=
			\delta_{R,R}\left(c_\fh\left(\widetilde Q\right)\right)
		\]
		where the inequality comes from the initial assumption on $\widetilde P$ and $\widetilde Q$, and the claim is proved.
		
		Therefore the Divisibility Axiom of fusion action systems implies that there is some $(\psi,\tau)\in\fX(\widetilde P,\widetilde Q)$ such that
		$
			c_{\widetilde{\fh\fg}}=(\iota,\id_X)\circ c_{\widetilde\fh}\circ(\psi,\tau)\in\fX\left(\widetilde P,N_S(R)\right)
		$.
		Now Proposition \ref{unique right lifting} implies that there is a unique $\widetilde\fg\in\cL^\fX\left(\widetilde P,\widetilde Q\right)$ such that $\widetilde{\fh\fg}=\widetilde\fh\widetilde\fg$.  Restricting this to $P$ we get $\fh\circ\fg=\fh\circ\res_P^Q(\widetilde\fg)$, and the fact that $\fh$ is mono implies that $\res_P^Q(\widetilde\fg)=\fg$, as desired.
	\end{proof}
	\end{prop}
	
	\begin{remark}
	Note that the condition on the overgroups can be restated as follows:  For every $p'\in\widetilde P$, there is some $q'\in\widetilde Q$ such that $\fg\circ\widehat {p'}\circ\fg^{-1}=\widehat {q'}$.  This condition is morally the same as the definition of the extender $N_{(\varphi,\sigma)}$ used to state the Extension Axiom for saturated fusion actions.  The key difference is that, when stated in terms of the linking action system, the extension condition on the source is sharper, which allows us to relax the assumption that the target of $\fg$ be fully $X$-centralized.
	\end{remark}
	
	Recall the notion of an \emph{abstract transporter system} associated to a fusion system $\cF$, as introduced in \cite{OliverVenturaTransporterSystems} and briefly described in Subsection \ref{group theory without groups}.
	
	\begin{cor}
		For a linking action system $\cL^\fX$ associated to $\fX$ the composite $\cL^\fX\to\fX\to\cF$, together with $\delta:\cT_S^{cX}\to\cL^\fX$, give $\cL^\fX$ the structure of a transporter system associated to $\cF$.
		\begin{proof}
			The only difficult part of the proof is the extension condition, which Proposition \ref{linking action extensions} shows to be true.
		\end{proof}
	\end{cor}
	
	We can interpret this result as saying that a fusion action system $\fX$ and an associated linking system $\cL^\fX$ give rise to a transporter system on the underlying fusion system $\cF$ together with a map $\Mor(\cL^\fX)\to\Sigma_X$ that takes composition to multiplication and inclusions to the identity.  In \cite{OliverVenturaTransporterSystems} it is shown that this map is equivalent to the data of a group map $\pi_1(|\cL^\fX|)\to\Sigma_X$.  We now set out to reverse this process:
	
	Fix a (saturated) fusion system $\cF$ on the $ p$-group $S$ and a transporter system $\cT$ associated to $\cF$.  Assume that we are given a group map $\theta:\pi_1(|\cT|)\to\Sigma_X$, or equivalently a map $\Mor(\cT)\to\Sigma_X$ that sends composition to multiplication and inclusions to the identity.
	
	\begin{definition}
		Let $\fX^\theta$ be the category with $\Ob(\fX^\theta)=\Ob(\cT)$ and morphisms given by
		\[
			\fX^\theta(P,Q)=\left\{(\varphi,\sigma)\in\Inj(P,Q)\times\Sigma_X\big|\exists\fg\in\cT(P,Q)\textrm{ such that }(\varphi,\sigma)=(c_\fg,\theta(\fg))\right\}.
		\]
	\end{definition}
	
	This allows us to define an action of $S$ on $X$ as follows:  The $ p$-group $S$ embeds as a subgroup of $\cT(S)$ via the structure map $\delta:\cT_S^{\Ob(\cT)}\to\cT$.  We denote by $\widehat S$ the image of $S$.  Moreover, $\theta$ defines a $\cT(S)$-action on $X$, and thus an $S$-action by restriction.
	
	Clearly $\fX^\theta$ is a fusion action system, or at least generates one once we allow for restrictions of morphisms to subgroups not in $\Ob(\cT)$.  Let $\cF^\theta$ be the underlying fusion system.
	
	We would for $\fX^\theta$ to be saturated, but for now we must settle for a weaker condition.
	
	\begin{definition}
		For $\cC$ a collection of subgroups of $S$ closed under $\cF^\theta$-conjugacy and overgroups, and $\fX$ a fusion action system on $S$, we say that $\fX$ is \emph{$\Ob(\cC)$-saturated} if the Saturation Axioms hold for all $P\in\cC$.
	\end{definition}

	We need a little terminology to prove object-saturation of $\cF^\theta$:
	
	\begin{notation}
		In the above situation, for any $P\in\Ob(\cT)$ we define $E(P)=\ker[\cT(P)\to\cF^\theta(P)]$.  We also denote by $K(P)$ the kernel of the action map $\theta:\cT(P)\to\Sigma(X)$.  Finally, let $C$ be the core of the $S$-action on $X$, so that $\widehat C=\widehat S\cap K(S)$.  We can therefore define the notions of $X$-normalizers and $X$-centralizers of objects of $\cT$ in the obvious way.
	\end{notation}

	\begin{prop}\label{X-whatever Sylow in transporter systems}
		For each $P\in\Ob(\cT)$,
		\begin{itemize}
			\item $P$ is fully normalized in $\cF^\theta$ if and only if $\widehat{N_S(P)}\in\Syl_ p(\cT(P))$.
			\item $P$ is fully centralized in $\cF^\theta$ if and only if $\widehat{Z_S(P)}\in\Syl_ p(E(P))$.
			\item $P$ is fully $X$-normalized in $\cF^\theta$ if and only if $\widehat{N_S(P;X)}\in\Syl_ p(K(P))$.
			\item $P$ is fully $X$-centralized in $\cF^\theta$ if and only if $\widehat{Z_S(P;X)}\in\Syl_ p(E(P)\cap K(P))$.
		\end{itemize}
	\begin{proof}
		The first two points are proved in \cite{OliverVenturaTransporterSystems}, Proposition 3.4.  The proofs of the remaining two points follow basically the same argument as the second.
		
		Proof of third point:  For $P\in\Ob(\cT)$, let $Q$ be $\cF^\theta$-conjugate to $P$ and fully normalized, so by the first point $\widehat{N_S(Q)}\in\Syl_ p(\cT(Q))$.   Therefore 
		\[
		\widehat{N_S(Q;X)}=\widehat{N_S(Q)\cap K(Q)}\in\Syl_ p\cT(Q).
		\]
		Now, $\cF^\theta$ is $\Ob(\cT)$-saturated by \cite{OliverVenturaTransporterSystems}, so the proof of Proposition \ref{X-whatever implications} applies here to give us that $P$ is also fully $X$-normalized.  Since $K(Q)\cong K(P)$, we have $\widehat{N_S(P;X)}\in\Syl_ p(K(P))$ if and only if $|N_S(P;X)|=|N_S(Q;X)|$, or equivalently, if and only if $P$ is fully $X$-normalized.
		
		The proof of the fourth point is the same as that of the third, replacing every instance of $K(-)$ with $E(-)\cap K(-)$.
	\end{proof}
	\end{prop}
	
	The notation in the following Corollary is a direct analogy with that introduced to describe fusion action systems.
	
	\begin{cor}
		For all $P\in\Ob(\cT)$
		\begin{itemize}
		\item If $P$ is fully normalized then
		\begin{itemize}
			\item$\fX^\theta_S(P)\in\Syl_ p(\fX^\theta(P))$.
			\item$\cF^\theta_S(P)\in\Syl_ p(\cF^\theta(P))$.
			\item$\Sigma_\theta^S(P)\in\Syl_ p(\Sigma_\theta(P))$.
		\end{itemize}
		\item If $P$ is fully $X$-normalized, then $\cF^\theta_S(P)_0\in\Syl_ p(\cF^\theta(P)_0)$.
		\item If $P$ is fully centralized, then $\Sigma_\theta^S(P)_0\in\Syl_ p(\Sigma_\theta(P)_0)$.
		\end{itemize}
		\begin{proof}
			Each of these follows from the observation that the image of a Sylow is Sylow in the quotient.
		\end{proof}
	\end{cor}
	
	\begin{prop}\label{extensions for transporter actions}
	Let $Q\in\Ob(\cT)$ be fully $X$-centralized and $(\varphi,\sigma)\in\Iso_{\fX^\theta}(P,Q)$.  Then there is some $(\widetilde\varphi,\sigma)\in\fX^\theta(N_{(\varphi,\sigma)},S)$ that extends $(\varphi,\sigma)$.
	\begin{proof}
		We first claim that  $\widehat{N_S(Q)}\in\Syl_ p\left(\widehat{N_S(Q)}\cdot E(Q)\cap K(Q)\right)$:
		\begin{eqnarray*}
			\left[\widehat{N_S(Q)}:\widehat{N_S(Q)}\cdot E(Q)\cap K(Q)\right]&=&
			\frac{\left|\widehat{N_S(Q)}\right|\cdot|E(Q)\cap K(Q)|}{\left|N_S(Q)\right|\cdot
			\left|\widehat{N_S(Q)}\cap E(Q)\cap K(Q)\right|}\\
			&=&\left[\widehat{N_S(Q)}\cap E(Q)\cap K(Q):E(Q)\cap K(Q)\right].
		\end{eqnarray*}
		The fact that $\widehat{N_S(Q)}\cap E(Q)\cap K(Q)=Z_S(Q;X)$ and the final point of Proposition~\ref{X-whatever Sylow in transporter systems} gives the claim.
		
		Now, pick $\fg\in\cT(P,Q)$ such that $(c_\fg,\theta(\fg))=(\varphi,\sigma)$.  By definition of $N_{(\varphi,\sigma)}$,
		\[
			\fg\circ\widehat{N_{(\varphi,\sigma)}}\big|_P^P\circ\fg^{-1}\leq\widehat{N_S(Q)}\big|_Q^Q\cdot E(Q)\cap K(Q)
		\]
		and so by the Sylow result just proved, there is some $\fh\in E(Q)\cap K(Q)$ such that
		\[
			(\fh\fg)\circ\widehat{N_{(\varphi,\sigma)}}\big|_P^P\circ(\fh\fg)^{-1}\leq\widehat{N_S(Q)}\big|_Q^Q.
		\]
		Then $\fh\fg\in\cT(P,Q)_\iso$, $N_{(\varphi,\sigma)}\trianglerighteq P$, and $N_S(Q)\triangleright Q$, so the conditions for the Extension Axiom (II) of abstract transporter systems are satisfied.  Therefore there is some $\widetilde{\fh\fg}\in\cT(N_{(\varphi,\sigma)},N_S(Q))$ that extends $\fh\fg$.  This implies that $(c_{\fh\fg},\theta(\fh\fg))\in\fX^\theta(N_{(\varphi,\sigma)},S)$ extends $(c_\fg,\theta(\fg))=(\varphi,\sigma)$ in $\fX^\theta$, and the result is proved.
	\end{proof}
	\end{prop}
	
	\begin{cor}\label{object saturation}
		The fusion action system $\fX^\theta$ is $\Ob(\cT)$-saturated.
		\begin{proof}
			All the Axioms have been verified in Propositions \ref{X-whatever Sylow in transporter systems} and \ref{extensions for transporter actions}.
		\end{proof}
	\end{cor}
	
	\begin{remark}
	Corollary \ref{object saturation} suggests that we should consider under what circumstances a collection of subgroups $\cH$ of $S$ has the property that if a fusion action system $\fX$ is $\cH$-saturated, then $\fX$ is itself saturated.  Such a result for fusion systems can be found in \cite{BCGLO1}, suggesting that it is not unreasonable to expect that the same would apply for fusion action systems.  Hopefully this point will be addressed in future research.
	\end{remark}
	
	We have seen that a transporter system $\cT$ together with a map $\theta:\pi_1(\cT)\to\Sigma_X$ determine a saturated fusion action system, at least so far as the objects of $\cT$ are aware.  We are given natural functors
	$
	\xymatrix
	{
	\cT_S^{\Ob(\cT)}\ar[r]&		\cT\ar[r]&	\fX^\theta
	}
	$, 
	and we can ask how close this is to being the data of an $X$-centric linking action system associated to $\fX^\theta$.  
	
	The following result states that, so long as all the $X$-centric subgroups are accounted for,  $\cT$ fails to be a linking action system ``in a $ p'$-way,'' and moreover that it contains enough data to construct a linking action system $\cL^\theta$:
	
	\begin{prop}
	Suppose that in the above situation $\Ob(\cT)$ contains all $X$-centric subgroups of $S$, and let $\cT^{cX}$ denote the fully subcategory with these as the objects.  Then for any $X$-centric $P\leq S$, there is a unique $ p'$-group $EK'(P)$ such that
	$
	E(P)\cap K(P)=\widehat{Z(P;X)}\times EK'(P)
	$.
	Furthermore, $EK'(P)$ is the subgroup of all $ p'$-elements of $E(P)\cap K(P)$. 
	
	Consequently, if $\cL^\theta$ is the category whose objects are the $X$-centric subgroups of $S$ and whose morphisms are given by
	\[
		\cL^\theta(P,Q)=\cT(P,Q)/EK'(P),
	\]
	then $\cL^\theta$ is an $X$-centric linking action system associated to $\fX^\theta$.
	\begin{proof}
		Axiom (C) of transporter systems implies that $E(P)$ commutes with $\widehat P$, so in particular $E(P)\cap K(P)$ does as well.  Therefore the fact that $P$ is $X$-centric implies that
		\[
			\widehat S\cap E(P)\cap K(P)=\widehat{Z_S(P;X)}=\widehat{Z(P;X)}\trianglelefteq E(P)\cap K(P),
		\]
		and the fact that $P$ is fully $X$-centralized implies that $\widehat{Z(P;X)}$ is a normal abelian Sylow subgroup of $E(P)\cap K(P)$.  The Schur-Zassenhaus theorem then implies the existence and uniqueness of $EK'(P)$, from which it easily follows that $\cL^\theta$ is an $X$-centric linking action system associated to $\fX^\theta$.
	\end{proof}
	\end{prop}
	
\subsection{Stabilizers of $ p$-local finite group actions}\label{stablizer p-local finite groups}

		For $H\leq G$ finite groups, it is a basic result that
	\[
		BH\simeq EG\times_H *\simeq EG\times_G G/H.
	\]

	In this section we prove that the analogous statement for linking actions systems is true. Let $\fX$ be a saturated fusion action system, and let $\cL^\fX$ be an associated linking action system.  Recall that $C$ denotes the core of the $S$ action on $X$.
	
	Recall that the fusion action system $\fX$ is \emph{transitive} if $\fX(1)$ acts transitively on $X$.
	
	\begin{lemma}
		If the fusion action system $\fX$ is saturated, $\pi_\Sigma\left(\fX(1)\right)=\pi_\Sigma\left(\fX(C)\right)$.
		\begin{proof}
			The non-obvious inclusion is $\pi_\Sigma\left(\fX(1)\right)\subseteq\pi_\Sigma\left(\fX(C)\right)$, which follows from the Extension Axiom for fusion action systems and from the easy calculation that $C\leq N_{(\id_1,\sigma)}$ for any $(\id_1,\sigma)\in\fX(1)$.
		\end{proof}
	\end{lemma}
	
	\begin{remark}
		We could therefore have defined transitivity of fusion actions in terms of the group $\pi_\Sigma\left(\fX(C)\right)\leq\Sigma_X$.  If we want to concentrate on linking action systems, this alternate characterization has the advantage that $C$ is always $X$-centric, and therefore  we can define transitivity of a linking action system in terms of subgroups of $S$ that are witnessed by $\cL^\fX$.
	\end{remark}
	
	We wish to introduce the notion of the ``stabilizer'' of a point $x\in X$ from the point of view of fusion and linking actions.  Recall that we denote by $\fX_x$ the stabilizer fusion action subsystem of $\fX$.  Moreover, we shall denote by $\cF_x$ the underlying fusion system of $\fX_x$.
		
	\begin{definition}
		Given a linking action system $\cL^\fX$ associated to $\fX$, the \emph{stabilizer linking action system of $x$}, $\cL^\fX_x=\cS tab_{\cL^\fX}(x)$, is the category whose objects are those $X$-centric subgroups that are contained in $S_x$ and whose morphisms are given by
		$
			\cL^\fX_x(P,Q)=~\left\{\fg\in\cL^\fX(P,Q)\big|\ell_\fg(x)=x\right\}
		$.
	\end{definition}

	\begin{remark}
	We can think of $\fX_x$ and $\cL^\fX_x$ as the preimages under the natural maps $\Mor(\fX)\to\Sigma_X$ and $\Mor(\cL^\fX)\to\Sigma_X$, respectively, of the subgroup $\Sigma_{X-\{x\}}$.
	\end{remark}
		
	Again, we must restrict our attention to the stabilizers of \emph{fully stabilized} points of $X$, namely those points whose stabilizers are of maximal order.  The following is the linking action system analogue of Lemma \ref{full stabilization in fusion action systems}.
	
	\begin{lemma}\label{fully stabilized characterization}  
		The point $x\in X$ is fully stabilized if and only if
		$
			\widehat{S_x}\big|_C^C\in\Syl_ p\left(\cL^\fX_x(C)\right)
		$.
		\begin{proof}
			The group $\cL^\fX(C)$ naturally acts on $X$ by the composition
			\[
			\xymatrix
			{
				\cL^\fX(C)\ar[rr]^-{\pi_{C,C}}&&	\fX(C)\ar[rr]^-{\pi_\Sigma}&&	\Sigma_X.
			}
			\]
			Axiom (B) of linking action systems implies that 
			\[
				\widehat{S_x}\big|_C^C=\left(\widehat{S}\big|_C^C\right)_x,
			\]
			and it follows easily from the definitions that
			\[
				\widehat{S_x}\big|_C^C=\left(\widehat S\big|_C^C\right)\cap\left(\cL^\fX_x(C)\right).
			\]
			The core $C$ is strongly closed in $\cF$, so in particular it is fully normalized and $N_S(C)=S$.  Proposition \ref{X-whatever Sylow in transporter systems} then implies that that
			$
				\widehat S\big|_C^C\in\Syl_ p\left(\cL^\fX(C)\right)
			$.
			
			Thus we have again reduced the problem to the case of actual finite groups, as in the proof of Lemma \ref{full stabilization in fusion action systems}, and the result follows.
		\end{proof}
	\end{lemma}
	
	With this interpretation of the stabilizer fusion, action, and linking systems, we find ourselves in the situation examined in \cite{OliverVenturaTransporterSystems}, and we recall the following result:
	
	\begin{prop}\label{Oliver Ventura 4.1}
		Let $\cT$ be an abstract transporter system associated to the fusion system $\cF$ on the $ p$-group $S$.  Fix a finite group $\Gamma$ and a group homomorphism $\Phi:\pi_1(|\cT|)\to\Gamma$, or equivalently, a map $\Mor(\cT)\to\Gamma$ that takes composition to multiplication and inclusions to the identity.  For any subgroup $H\leq\Gamma$ let $S_H\leq S$ be the maximal subgroup whose elements (viewed as morphisms of $\cT$) are sent to $H$, and assume that that $S_1\in\Ob(\cT)$. 
		
		Let $\cT_H\subseteq\cT$ be the subcategory whose objects are those of $\cT$ that are contained in $H$ and whose morphisms are given by
		$
		\cT_H(P,Q)=\left\{\fg\in\cT(P,Q)\big|\Phi(\fg)\in H\right\}
		$.
		
		Let $\cF_H\subseteq\cF$ be the fusion system on $S_H$ generated by $\pi(\cT_H)$, and let
		\[
		\xymatrix
		{
		\cT_{S_H}^{\Ob(\cT_H)}(S_H)\ar[rr]^-{\delta_{H,H}}&&
		\cT_H\ar[rr]^-{\pi_H}&&
		\cF_H
		}
		\]
		be the restrictions of the structure maps for the transporter system $\cT$.	Then:
		\begin{itemize}
		\item[(a)]  $\Phi(\Mor(\cT))=\Phi(\cT(S_1))$.
		\item[(b)]  $\cT_H$ is a transporter system associated to $\cF_H$ if and only if $\delta_{S_1,S_1}(S_H)\in\Syl_ p(\cT_H(S_1))$.
		\item[(c)]  If the condition of Item (b) is satisfied and all fully centralized $P\leq S$ have the property that $Z_{S_1}(P)\leq P$ implies $P\in\Ob(\cT)$ then $\cF_H$ is a saturated fusion system.
		\item[(d)]  If for all $P\in\Ob(\cT)$ we have $P\cap S_1\in\Ob(\cT)$, then $|\cT_H|$ has the homotopy type of the covering space of $|\cT|$ with fundamental group $\Phi^{-1}(H)$.
		\end{itemize}
		\begin{proof}
		\cite[Proposition 4.1]{OliverVenturaTransporterSystems}.
		\end{proof}
	\end{prop}
	
	\begin{prop}
	If $x\in X$ is fully stabilized, the stabilizer fusion system $\cF_x\subseteq\cF$ is saturated and $\cL^\fX_x$ is a transporter system associated to it.
	\begin{proof}
		From the natural map $\Mor(\cL^\fX)\to\Sigma_X$ and $H=\Sigma_{X-\{x\}}\leq\Sigma_X$, in the notation of Proposition \ref{Oliver Ventura 4.1} we have $C=S_1$, $S_x=S_H$, $\cL^\fX_x=\cT_H$, and $\cF_x=\cF_H$.  That $\cL^\fX_x$ is a transporter system associated to $\cF_x$ is then simply an application of Lemma \ref{fully stabilized characterization} (b) to Proposition \ref{Oliver Ventura 4.1} (b).  
		
		That $\cF_x$ is saturated follows from Theorem \ref{stabilizer saturation}, which implies that the stabilizer fusion action subsystem $\fX_x$ is saturated.  However, with the additional structure of a linking action system, the result is easier to prove, thanks to the work of Oliver and Ventura.  Item (c)  of Proposition \ref{Oliver Ventura 4.1} shows that if $P\leq S$ is fully centralized and $Z_C(P)\leq P$, then $P$ is $X$-centric, i.e.,  that $Z_C(Q)=Z_S(Q;X)=Z( Q;X)$ for all $Q$ $\cF$-conjugate to $P$.  If $\varphi\in\Iso_\cF(P,Q)$, then $\varphi(Z(P;X))=Z(Q;X)\leq Z_S(Q;X)$.  Because $P$ is fully centralized, it is also fully $X$-centralized, so the assumption that $Z(P;X)=Z_S(P;X)$ and comparison of orders implies that $Z(Q;X)=Z_S(Q;X)$.  Thus the conditions of Item (c) are satisfied and $\cF_x$ is saturated.
	\end{proof}
	\end{prop}
	
	We now find ourselves in the following situation:  Let $\cL^\fX$ be a linking action system associated to the transitive saturated fusion action system $\fX$, and $x$ a fully stabilized point of $X$.  We would like to understand the topological information of the stabilizer linking action system $\cL^\fX_x$ as it relates to that of $\cL^\fX$, and indeed Proposition \ref{Oliver Ventura 4.1} gives us some relevant information in terms of subgroups of $\pi_1(|\cL^\fX|)$.  We can also calculate the homotopy type of $|\cL^\fX_x|$ directly, as follows:
	
	Let  $\iota:\cL^\fX_x\to\cL^\fX$ be the inclusion functor, and let $F:\cL^\fX\to\TOP$ be the left homotopy Kan extension of the trivial functor $*:\cL^\fX_x\to\TOP$ over $\iota$.  We already have the functor $\XX:\cL^\fX\to\TOP$ defined as part of the data of $\cL^\fX$, and we would like to relate these.
	
	\begin{prop}
		$F$ is equivalent to $\XX$ as functors $\cL^\fX\to\TOP$.
		\begin{proof}
			Thomason's theorem \cite{ThomasonHomotopyColim} tells us that we have a homotopy equivalence
			\[
			F(P)\simeq\hocolim_{(\iota\downarrow P)}*=|(\iota\downarrow P)|
			\]
			for all $P\in\Ob(\cL^\fX)$.  We first prove that every component of $|(\iota\downarrow P)|$ is contractible and that the components can be put in natural correspondence with the points of $X$.
			
			The equivalence of $F$ with $\XX$ is not natural:  For every $y\in X$, Corollary \ref{stabilizer subconjugacy} implies that $\cL^\fX(S_y,S_x)$ is nonempty as $x$ is fully stabilized.  Let $\fg_y\in\cL^\fX(S_y,S_x)$ be a choice of a morphism in this hom-set for each $y$, and assume $\fg_x=\id_{S_x}$.  Also let $\overline{\fg_y}=(\fg_y)_{\mathrm{iso}}$ be the restricted isomorphism of $\fg_y$ (cf. Corollary \ref{factorization in linking action systems}), $\fg_y^P$ the restriction of $\fg_y$ to $P_y\leq S_y$, and $\overline{\fg_y^P}=\left(\fg_y^P\right)_{\mathrm{iso}}$ the restricted isomorphism of $\fg_y^P$.
			
			For any $(Q,\alpha)\in(\iota\downarrow P)$, so that $Q\in\Ob(\cL^\fX_x)$ and $\alpha\in\cL^\fX(Q,P)$, we have $\ell_\alpha(x)=y$ for some $y\in X$.  Since the pair $(c_\alpha,\ell_\alpha)$ is intertwined and $Q\leq S_x$, we have $c_\alpha(Q)\leq P_y\leq S_y$.  Thus we have the following commutative diagram in $\cL^\fX$, where all morphisms labeled $\fri$ are the obvious inclusions:
			\[
			\xymatrix{
			S_x\\
			c_{\fg_y(S_y)}\ar[u]^\fri&&S_y\ar[ll]^{\overline{\fg_y}}_\cong\ar[ull]_{\fg_y}\\
			c_{{\fg_y}(P_y)}\ar[u]^\fri&&P_y\ar[ll]^{\overline{\fg_y^P}}_\cong\ar[u]^\fri\ar[rr]^\fri&&P\\
			Q\ar[u]^\fri\ar[rr]_{\alpha_\iso}^\cong&&c_\alpha(Q)\ar[u]^\fri\ar[rru]_\fri
			}
			\]
			The top composition is $\alpha$ by definition.  By the choices made above, we have $\ell_\fg\circ\ell_\alpha(x)=x$, and therefore there is a  factorization
			\[
			\alpha=\underbrace{\left(\fri_{P_y}^P\circ\overline{\fg_y^P}^{-1}\right)}_{\cL^\fX(c_{\fg_y}(P_y),P)}\circ\underbrace{\left(\overline{\fg_y^P}\circ\fri_{c_\alpha(Q)}^{P_y}\circ\alpha_{\mathrm{iso}}\right)}_{\cL^\fX_x(Q,c_{\fg_y}(P_y))}.
			\]
			This is the unique (as all morphisms of $\cL^\fX$ are categorically mono and epi) factorization of $\alpha:Q\to P$ as a composite $Q\to c_{\fg_y}(P_y)\to P$, so what we have really proved is the following:
			
			For any object $(Q,\alpha)\in(\iota\downarrow P)$ such that $\ell_\alpha(x)=y$, there is a unique morphism from $(Q,\alpha)$ to $\left(c_{\fg_y}(P_y),\iota_{P_y}^P\circ\overline{\fg_{y}^P}^{-1}\right)$ in $(\iota\downarrow P)$, namely $\overline{\fg_y^P}\circ\fri_{c_\alpha(Q)}^{P_y}\circ\alpha_{\mathrm{iso}}$.  In other words, we have found a terminal object in the component containing $(Q,\alpha)$ that depends only on our choice of the $\fg_y$ and where $\alpha$ sends $x$.  This shows that $F(P)$ is homotopically discrete, and as we have assumed that $\cL^\fX$ is transitive, the components are naturally identified with the points of $X$.
			
			All we have to do is see how $F$ acts on morphisms, compared to the functor $\XX$.  Recall that $\XX$ sends the morphism $\fh\in\cL^\fX(P,P')$ to the map of spaces $\ell_\fh$.  On the other hand,  $F(\fh)$ is induced by the functor $(\iota\downarrow P)\to(\iota\downarrow P')$ that sends $(Q,\alpha)$ to $(Q,\fh\circ\alpha)$.  If $(Q,\alpha)$ is in the component we have identified with $y$ we have $\ell_\alpha(x)=y$, and then $(Q,\fh\circ\alpha)$ is in the component corresponding to $\ell_\fh\circ\ell_\alpha(x)=\ell_\fh(y)$.  This is just to say that $F(\fh)$ permutes the space homotopy equivalent $X$ by the permutation $\ell_\fh$, so the result is proved.
		\end{proof}
	\end{prop}
	
	\begin{theorem}\label{stabilizer classifying spaces}
		In the above situation,
		\[
			\hocolim_{\cL^\fX}\XX\simeq\left|\cL^\fX_x\right|.
		\]
	\end{theorem}
	
	\begin{remark}
	The final piece of interpretation comes from thinking of the left hand side as the $ p$-local finite group action-theoretic version of $EG\times_G G/H$, and the right hand side as $EH\times_H*$. 
	\end{remark}

\bibliography{Sources}
\bibliographystyle{alphanum}

\end{document}